\documentclass[12pt]{article}
\usepackage{amsmath}
\usepackage{amssymb}
\usepackage{amsthm}
\usepackage[totalwidth=17cm, totalheight=25cm]{geometry}
\usepackage{graphicx}
\usepackage{mathabx}
\usepackage{tabularx}
	\newcolumntype{C}[1]{>{\centering\arraybackslash}m{#1}} 
	\newcolumntype{R}[1]{>{\raggedleft\arraybackslash}m{#1}} 

\newtheoremstyle{boldplain}
{9pt}
{9pt}
{\itshape}
{}
{\bfseries}
{.}
{.5em}
{\thmname{#1}\thmnumber{ #2}\thmnote{ (#3)}}%

\newtheoremstyle{bolddefinition}
{9pt}
{9pt}
{}
{}
{\bfseries}
{.}
{.5em}
{\thmname{#1}\thmnumber{ #2}\thmnote{ (#3)}}%

\theoremstyle{boldplain}

\newtheorem{cor}[equation]{Corollary}

\newtheorem{lem}[equation]{Lemma}
\newtheorem{lemma}[equation]{Lemma}
\newtheorem{prop}[equation]{Proposition}
\newtheorem{sublem}[equation]{Sublemma}
\newtheorem{thm}[equation]{Theorem}

\theoremstyle{bolddefinition}
\newtheorem{dfn}[equation]{Definition}
\newtheorem{definition}[equation]{Definition}
\newtheorem{defn}[equation]{Definition}

\newtheorem{example}[equation]{Example}

\newtheorem{rem}[equation]{Remark}

\newfont{\bigbf}{cmbx10 scaled\magstep1}
\normalbaselines
\numberwithin{equation}{section}

\def\no{\noindent}

\def\C{{\mathbb C}}
\def\R{{\mathbb R}}
\def\H{{\mathbb H}}

\def\N{{\mathbb N}}

\def\Z{{\mathbb Z}}

\def\al{\alpha}
\def\be{\beta}
\def\ga{\gamma}
\def\Ga{\Gamma}
\def\de{\delta}
\def\De{\Delta}
\def\eps{\epsilon}
\def\la{\lambda}
\def\La{\Lambda}
\def\si{\sigma}
\def\Si{\Sigma}
\def\ups{\upsilon}

\def\Om{\Omega}

\def\3{\ss}

\def\acts{\curvearrowright}

\def\D{\partial}
\def\Diag{\mathop{\hbox{Diag}}}
\def\diam{\mathop{\hbox{diam}}}
\def\diamo{\diamondsuit}
\def\embed{\hookrightarrow}
\def\Homeo{\operatorname{Homeo}}
\def\B{\operatorname{B}}
\def\Fix{\operatorname{Fix}}
\def\Flag{\operatorname{Flag}}

\def\Flagt{\Flag(\tau_{mod})}

\def\geo{\partial_{\infty}}

\def\Hom{\operatorname{Hom}}
\def\id{\mathop{\hbox{id}}}

\def\interior{\operatorname{int}}
\def\Isom{\mathop{\hbox{Isom}}}
\def\LasGa{\La_{\si_{mod}}(\Ga)}
\def\LatGa{\La_{\tau_{mod}}(\Ga)}
\def\lra{\longrightarrow}

\def\midp{\operatorname{mid}}

\def\ol{\overline}

\def\pihalf{\frac{\pi}{2}}

\def\2pithird{\frac{2\pi}{3}}
\def\piforth{\frac{\pi}{4}}

\def\pos{\mathop{\hbox{pos}}}

\def\Ra{\Rightarrow}
\def\st{\operatorname{st}}
\def\ost{\mathop{\hbox{ost}}}
\def\Stab{\operatorname{Stab}}
\def\tangle{\angle_{Tits}}
\def\tauto{\buildrel f\over\lra}
\def\tits{\partial_{Tits}}

\def\Trip{\mathop{\hbox{Trip}}\nolimits}

\def\Vert{\operatorname{Vert}}

\def\8{\infty}
\def\<{\langle}
\def\>{\rangle}

\def\BI{\begin{itemize}}
\def\EI{\end{itemize}}

\title{Morse actions of discrete groups on symmetric spaces}
\author{Michael Kapovich, Bernhard Leeb, Joan Porti}
\date{March 29, 2014}

\begin{document}

\maketitle

\begin{abstract}
\no
We study the geometry and dynamics of discrete infinite covolume subgroups
of higher rank semisimple Lie groups.
We introduce and prove the equivalence of several conditions, 
capturing ``rank one behavior'' 
of discrete subgroups of higher rank Lie groups. 
They are direct generalizations 
of rank one equivalents to convex cocompactness.
We also prove that our notions are equivalent 
to the notion of Anosov subgroup,
for which we provide a closely related, 
but simplified and more accessible reformulation,
avoiding the geodesic flow of the group.
We show moreover that the Anosov condition can be relaxed further 
by requiring only non-uniform unbounded expansion 
along the (quasi)geodesics in the group.

A substantial part of the paper is devoted to the coarse geometry
of these discrete subgroups. 
A key concept which emerges from our analysis 
is that of {\em Morse quasigeodesics}
in higher rank symmetric spaces, 
generalizing the Morse property for quasigeodesics 
in Gromov hyperbolic spaces. 
It leads to the notion of {\em Morse actions} 
of word hyperbolic groups on symmetric spaces,
i.e.\ actions for which the orbit maps are Morse quasiisometric embeddings,
and thus provides a coarse geometric characterization 
for the class of subgroups considered in this paper.
A basic result is a local-to-global principle 
for Morse quasigeodesics and actions. 
As an application of our techniques 
we show algorithmic recognizability 
of Morse actions and construct Morse ``Schottky subgroups'' 
of higher rank semisimple Lie groups via arguments 
not based on Tits' ping-pong.
Our argument is purely geometric 
and proceeds by constructing    
equivariant Morse quasiisometric embeddings of trees 
into higher rank symmetric spaces. 
\end{abstract}

\tableofcontents

\section{Introduction}

This paper is motivated by the search for 
``geometrically nice" infinite covolume discrete groups of isometries 
of higher rank symmetric spaces, 
respectively, discrete subgroups of semisimple Lie groups with finite center. 
A reasonable class of such groups should be broad enough and at the same time manageable. 

The class of groups considered in this paper can be viewed as groups with 
{\em rank one actions on higher rank symmetric spaces}, 
i.e.\ discrete subgroups of semisimple Lie groups 
which exhibit some rank one behavior.
They are a direct generalization of {\em convex cocompact} subgroups of rank one Lie groups. 
The strength of the notion of convex cocompactness in rank one relies on a number of different 
characterizations in terms of 
geometry/dynamics (conical limit set), 
dynamics (expansion at the limit set, uniform convergence action on the limit set), 
coarse geometry (undistorted)
and topology (existence of a natural compactification of the quotient locally symmetric space). 
Furthermore, 
for these subgroups one can prove results which are inaccessible or unavailable otherwise 
(e.g.\ topological or algebraic finiteness properties). 
In chapter~\ref{sec:rank1}, 
we will go through different characterizations of rank one convex cocompact groups. 

In higher rank,
some of these characterizations turn out to be too restrictive, 
and others too weak. 
For example,
it was shown by Kleiner and Leeb \cite{convcoco}
that convex cocompactness is too restrictive in higher rank,
as it is satisfied only by few subgroups. 
On the other hand,
undistortion is way too weak;
for instance, undistorted subgroups can fail to be finitely presented.
Thus, one is forced to look for alternative generalizations of convex cocompactness in higher rank. 

In this paper, 
we will consider, in the context of {\em weakly regular} discrete subgroups $\Gamma\subset G$ 
of semisimple Lie groups (with finitely many connected components and finite center),  
four notions generalizing convex cocompactness to higher rank: 

(i) conical (RCA) subgroups

(ii) subgroups expanding at a (suitable) limit set

(iii) asymptotically embedded subgroups

(iv) Morse subgroups

Whereas (i)-(iii) are asymptotic (mostly dynamical) conditions, 
condition (iv) is coarse geometric. 
The Morse condition is a suitable strengthening of undistortedness in higher rank; 
the orbit quasigeodesics must satisfy an additional restriction
which is a higher rank analogue of the Morse property for quasigeodesics in hyperbolic spaces. 
See below for the definitions. 

A consequence of the equivalence (iii)$\Leftrightarrow$(iv) 
is the {\em structural stability} of asymptotically embedded subgroups, 
generalizing Sullivan's Structural Stability Theorem in rank one \cite{Sullivan}. 
Another consequence of this equivalence is 
that asymptotically embedded subgroups are undistorted and uniformly regular. 

Furthermore, 
an important feature of the Morse condition (iv) is that it admits a {\em local} characterization. 
The localness of the Morse property implies that the space of Morse actions is open, 
and that Morse subgroups are algorithmically recognizable.

We illustrate our techniques by constructing Morse-Schottky actions of free groups on higher rank symmetric spaces.
Unlike all previously known constructions, our proof does not rely on ping-pong, but
is purely geometric and proceeds by constructing  equivariant quasi-isometric embeddings of trees. 
The key step is the observation that a certain local {\em straightness} property for 
sufficiently spaced sequences of points in the symmetric space implies the global Morse property. 
This observation is also at the heart of the proof 
of the local-to-global principle for Morse actions. 

An advantage of the notion of conicality/expansivity
over the notions of asymptotic embeddedness and Morse 
is that it does not a priori assume that the subgroup is word hyperbolic as an abstract group. 
For example, uniform lattices in semisimple Lie groups satisfy this property. 
Thus, the conicality condition can serve as an ingredient in a general theory of ``geometric niceness" 
in higher rank which includes groups which are not (relatively) hyperbolic. 
 
In our paper \cite{coco}, 
we study the dynamics of weakly regular antipodal discrete subgroups acting 
on the partial flag manifolds associated to semisimple Lie groups.
We construct domains of discontinuity in general. 
Furthermore, 
we prove cocompactness of actions under the additional assumption of expansivity, equivalently, conicality 
using some of the discussion in this paper. 
The latter results can be regarded as a weak form of cocompactness 
for the action on the domain of discontinuity at infinity in rank one. 

We prove (section \ref{sec:weak Anosov}, Theorem \ref{thm:anosov iff morse}) 
that our characterizations (i)-(iv) are equivalent 
to the notion of {\em Anosov representation} (subgroup) 
introduced by Labourie \cite{Labourie} and developed further by Guichard and Wienhard 
\cite{GW}. 
In section~\ref{sec:anosov}, we also give a closely related, 
but simplified and more accessible reformulation of the Anosov property.
Our definition involves only an expansion property for the group action on the flag manifold  
and avoids the notion of expansion/contraction of flows on bundles.     
In particular, 
it avoids using the geodesic flow for hyperbolic groups 
(whose construction is highly involved for groups which do not arise as the fundamental group of a closed negatively curved Riemannian manifold). 
In section \ref{sec:weak Anosov},
we present a further relaxation of our version of the Anosov condition,
by requiring only non-uni\-form unbounded expansion along 
quasigeodesics in $\Ga$.
Nevertheless, we show that the resulting class of subgroups
remains the same (Theorem \ref{thm:anosov iff morse}). 

While our methods are independent of the ones in \cite{Labourie,GW},
it was known before that Anosov subgroups are undistorted (\cite[Theorem 5.3]{GW}), uniformly regular (\cite[Proposition 3.16]{GW}), 
and that Anosov representations are structurally stable (\cite[Theorem 5.13]{GW}). 

In addition to the equivalent conditions mentioned above, 
we also introduce a notion of {\em boundary embeddedness} for discrete subgroups 
which are word hyperbolic as abstract groups. 
It is a priori a weakening of asymptotic embeddedness, 
as it only requires an antipodal equivariant embedding of the Gromov boundary 
into the flag manifold. 
Here, antipodality means that distinct points go to opposite simplices. 
Boundary embeddedness constitutes the non-dynamical part of the definition of Anosov subgroup, 
omitting the expansivity. 
We show that boundary embeddedness is equivalent to asymptotic embeddedness in the regular case
(see Proposition~\ref{prop:homeochlim}). 

\medskip
We now describe in more detail some of the concepts and results discussed above.

Let $X=G/K$ be a symmetric space of noncompact type. 
Of primary importance to us will be the visual boundary $\geo X$ 
and the fact that it carries a natural structure as a spherical building $\tits X$. 
We recall that the visual limit set $\La(\Ga)\subset\geo X$ of a discrete subgroup $\Ga\subset G$ 
is the set of accumulation points of an(y) $\Ga$-orbit $\Ga x\subset X$. 

Weak regularity of discrete subgroups is an asymptotic condition 
regarding the directions of segments connecting orbit points. 
It is defined with respect to a fixed face $\tau_{mod}$ of the model Weyl chamber $\si_{mod}$ 
associated to $G$ and $X$. 
The condition of $\tau_{mod}$-regularity is a relaxation of uniform $\tau_{mod}$-regularity 
which requires that all limit points of $\Ga$ are $\tau_{mod}$-regular, 
see Definition~\ref{def:tauregsubs}. 
In the case when $\tau_{mod}=\si_{mod}$, 
uniform regularity simply means that all limit points of $\Ga$ are interior points of chambers 
in the Tits building. 

The $\tau_{mod}$-regularity of a discrete subgroup $\Ga\subset G$
can be read off its dynamics at infinity, 
namely it is equivalent to a contraction-expansion property 
of the $\Ga$-action on the partial flag manifold 
$$
\Flag(\tau_{mod})=G/P_{\tau_{mod}}
$$ 
associated to $G$ and $\tau_{mod}$, 
generalizing convergence dynamics in rank one, 
see section~\ref{sec:contrexp}.

For $\tau_{mod}$-regular subgroups and, more generally, for subgroups containing diverging 
$\tau_{mod}$-regular sequences, we have a {\em $\tau_{mod}$-limit set } 
$$
\La_{\tau_{mod}}(\Ga)\subset\Flag(\tau_{mod}) .
$$ 
In the uniformly $\tau_{mod}$-regular case, the $\tau_{mod}$-limit set is the natural projection of 
$\La(\Ga)$ to $\Flag(\tau_{mod})$. 
Furthermore, in the uniformly regular case, 
$\La_{\si_{mod}}(\Ga)$ is the set of chambers which contain ordinary limit points.
We will refer to points in $\La_{\tau_{mod}}(\Ga)$ as {\em limit simplices} or {\em limit flags}. 
We call $\Ga$ {\em $\tau_{mod}$-non-elementary} if 
it has at least three limit flags of type $\tau_{mod}$.

\begin{dfn}[Antipodal]
A $\tau_{mod}$-regular discrete subgroup $\Ga\subset G$ is called {\em $\tau_{mod}$-antipodal},
if any two distinct limit simplices in $\La_{\tau_{mod}}(\Ga)$ are opposite. 
\end{dfn}

The notion of conicality of limit simplices is due to Albuquerque \cite[Def.\ 5.2]{Albuquerque}.
For simplicity, we give it here only in the regular case 
and refer the reader to Definition \ref{dfn:conicaltau} for the general case. 
\begin{dfn}[Conicality]
\label{dfn:conicalintro}
A $\si_{mod}$-regular discrete subgroup $\Ga\subset G$ is called {\em $\si_{mod}$-conical}
if for every limit chamber $\si\in\La_{\si_{mod}}(\Ga)$
there exists a sequence $\ga_n\to\infty$ in $\Ga$ 
such that for a(ny) point $x\in X$ the sequence of orbit points $\ga_nx$ 
is contained in a tubular neighborhood of the euclidean Weyl chamber $V(x,\si)$
with tip $x$ and asymptotic to $\si$. 
\end{dfn}

\begin{dfn}[RCA]
A discrete subgroup $\Ga\subset G$ is called {\em $\tau_{mod}$-RCA} 
if it is $\tau_{mod}$-regular, $\tau_{mod}$-conical and $\tau_{mod}$-antipodal. 
\end{dfn}

Following Sullivan \cite{Sullivan},
we call a subgroup expanding at infinity
if its action on the appropriate partial flag manifold is expanding at the limit set.
More precisely:
\begin{dfn}[Expanding]
We call a $\tau_{mod}$-regular discrete subgroup $\Ga\subset G$ 
{\em $\tau_{mod}$-expanding at the limit set}
if for every limit flag in $\La_{\tau_{mod}} (\Ga)$ 
there exists a neighborhood $U$ in $\Flag(\tau_{mod})$ and an element $\ga\in\Ga$ 
which is uniformly expanding on $U$,
i.e.\ for some constant $c>1$ and all $\tau_1,\tau_2\in U$ we have:
\begin{equation*}
d(\ga\tau_1,\ga\tau_2)\geq c\cdot d(\tau_1,\tau_2) 
\end{equation*}
Here, the distance $d$ is induced by a Riemannian background metric on the flag manifold. 
\end{dfn}

\begin{dfn}[Asymptotically embedded]
We call a $\tau_{mod}$-antipodal $\tau_{mod}$-regular discrete subgroup $\Ga\subset G$ 
{\em $\tau_{mod}$-asymptotically embedded} 
if $\Gamma$ is word hyperbolic 
and there exists a $\Ga$-equivariant homeomorphism 
\begin{equation*}
\label{eq:mapalphatauintro}
\alpha: \geo \Gamma \buildrel\cong\over\to 
\La_{\tau_{mod}} (\Ga)\subset \Flag(\tau_{mod})
\end{equation*}
of its Gromov boundary onto its $\tau_{mod}$-limit set. 
\end{dfn}

For simplicity, 
we define the Morse property only in the regular case, 
see Definitions~\ref{dfn:mqg}, \ref{dfn:morseemb} and~\ref{dfn:morseact} for the general case.
\begin{dfn}[Morse]
(i) A uniformly regular quasigeodesic ray in $X$ is called {\em $\si_{mod}$-Morse} 
if it converges to a chamber at infinity in uniform conical fashion,
compare Definition~\ref{dfn:conicalintro}. 

(ii) 
An isometric action $\rho:\Ga\acts X$ of a 
word hyperbolic 
group $\Ga$ is called 
{\em $\si_{mod}$-Morse} if its orbit map sends uniform quasigeodesic rays in $\Ga$ 
to uniformly $\si_{mod}$-Morse quasigeodesic rays in $X$. 
In this case the image $\rho(\Ga)$ will be called a {\em $\si_{mod}$-Morse subgroup}. 
\end{dfn}
It follows immediately from the definition that Morse actions are properly discontinuous 
and Morse subgroups are discrete and undistorted. 

The first main result of this paper is
(see Theorems~\ref{thm:regopphyp}, \ref{thm:anosov iff morse} and~\ref{thm:asymorse}):
\begin{thm}[Equivalence]
\label{thm:equiv}
For $\tau_{mod}$-non-elementary $\tau_{mod}$-regular discrete subgroups $\Ga\subset G$ 
the following properties are equivalent:

1. $\Ga$ is $\tau_{mod}$-RCA. 

2. $\Ga$ is $\tau_{mod}$-antipodal and its action on $\Flag(\tau_{mod})$ is expanding at 
$\La_{\tau_{mod}}(\Ga)$. 

3. $\Ga$ is $\tau_{mod}$-asymptotically embedded. 

4. $\Ga$ is $\tau_{mod}$-Morse. 

5. $\Ga$ is $\tau_{mod}$-Anosov. 
\end{thm}
Whereas properties 3-5 include that $\Ga$ is word hyperbolic,
it follows that subgroups having property 1 or 2 must be word hyperbolic.

Among the consequences of this theorem are (see Theorem~\ref{thm:qiembrtaureg}):
\begin{thm}
[Undistortedness and uniform regularity] 
A subgroup $\Ga\subset G$ satisfying one of the equivalent properties  1-5
is $\tau_{mod}$-uniformly regular and undistorted. 
\end{thm}

We recall that 
quasigeodesics in Gromov hyperbolic spaces can be recognized locally 
by looking at sufficiently large finite pieces. 
Our second main theorem is 
an analogous result for Morse quasigeodesics in symmetric spaces and, as a consequence,
for Morse actions.
Morse actions can be recognized locally by looking at the image of sufficiently large balls 
in the group under the orbit map (see Theorem~\ref{thm:locglobmqiembhyp}):
\begin{thm}[Local-to-global]
For a word hyperbolic group $\Ga$, 
locally $\tau_{mod}$-Morse actions $\Ga\acts X$ 
(for suitable parameters) 
are $\tau_{mod}$-Morse actions.
\end{thm}
As a consequence of the localness and equivalence theorems, we obtain
(see Theorems~\ref{thm:morsestab} and~\ref{thm:strcstab}): 
\begin{thm}
[Openness of the space of Morse actions] 
For a word hyperbolic group $\Ga$, the subset of $\tau_{mod}$-Morse actions is open in $\Hom(\Ga, G)$. 
\end{thm}

\begin{thm}
[Structural stability] 
Let $\Ga$ be word hyperbolic.
Then for $\tau_{mod}$-Morse actions $\rho:\Ga\acts X$, 
the boundary embedding $\al_{\rho}$ depends continuously on the action $\rho$. 
\end{thm}
Thus, actions sufficiently close to a faithful Morse action 
are again discrete and faithful. 

The localness implies furthermore, that Morse actions are algorithmically recognizable
(see section~\ref{sec:algrec}):
\begin{thm}
[Algorithmic recognition of Morse actions] 
Let $\Ga$ be word hyperbolic.
Then there exists an algorithm whose inputs are homomorphisms $\rho: \Ga\to G$ (defined on generators of $\Ga$) 
and which terminates if and only if $\rho$ defines a $\tau_{mod}$-Morse action $\Ga\acts X$. 
\end{thm}

If the action is not Morse, the algorithm runs forever. 
We do not know if there is an algorithm which recognizes non-Morse actions. 
Note, that even in hyperbolic 3-space no algorithm is known which recognizes 
that a finitely generated group is not convex cocompact.

\medskip
{\bf Organization of the paper.} 
In section 
\ref{sec:prelim} we review basic definitions from the theory of symmetric spaces and spherical buildings. We also prove some results on geometry of such spaces, including geometry of parallel sets and associated decompositions of $X$, cones in $X$ over certain subsets of $\geo X$ and dynamics 
of transvections of $X$. 
In section \ref{sec:dyn prelim} we discuss several standard concepts of topological dynamics, namely, convergence actions, expansivity and conical limit points. 
In section \ref{sec:rank1} we give a list of equivalent definitions of convex cocompact subgroups of rank one Lie groups. 
In section \ref{sec:regularity and limit sets} we introduce several key asymptotic notions describing geometry of discrete groups, such as regularity and uniform regularity, various limit sets as well as antipodality and conicality of limit sets in partial flag manifolds. In section \ref{sec:Asymptotic properties} we prove the equivalence of the first three items in Theorem \ref{thm:equiv} and establish some fundamental properties of asymptotically embedded groups which lead to the concept of Morse quasigeodesics and Morse actions. 
We conclude the section by discussing the Anosov condition. 
The notions of Morse quasigeodesics and actions 
are discussed in detail in section \ref{sec:morse}. 
In that section, among other things, 
we establish local-to-global principles, 
prove Theorem \ref{thm:equiv}
and show that Morse actions are structurally  stable and algorithmically recognizable. We also construct 
Morse-Schottky actions of free groups on symmetric spaces.

\medskip 
{\bf Acknowledgements.} 
The first author was supported by NSF grants DMS-09-05802 and DMS-12-05312. The last author was supported by grants  Mineco MTM2012-34834 and AGAUR  SGR2009-1207. The three authors are also grateful to 
the GEAR grant which partially supported the IHP trimester in Winter of 2012 (DMS 1107452,
1107263, 1107367 ``RNMS: Geometric structures and representation varieties'' 
(the GEAR Network)), and to the Max Planck Institute for Mathematics in Bonn, where some of this work was done.

\section{Geometric preliminaries}
\label{sec:prelim}

In this section we collect some standard material on  
Coxeter complexes, the geometry of nonpositively curved symmetric spaces and associated 
spherical Tits buildings; we refer the reader to \cite{qirigid} and \cite{habil} 
for more detailed discussion of symmetric spaces and buildings.
We also prove some new results on the geometry of parallel sets and cones in symmetric spaces.

We start with some general notations:

A {\em tubular neighborhood} of a subset $A$ of a metric space $X$ is an open $r$-neighborhood of $A$ in $X$
for some $r>0$, i.e.\ the set 
$$
N_r(A)= \{x\in X: d(x, A)<r\}. 
$$
We will use the notation
$$
{B}(a,r)=\{x\in X: d(x,a)< r\}
$$
and
$$
\bar{B}(a,r)=\{x\in X: d(x,a)\le r\}
$$
for the open and, respectively, closed $r$-ball, centered at $a$. 

A {\em geodesic} in a metric space $X$ is an isometric embedding from a (possibly infinite) interval $I\subset\R$ into $X$. In the context of finitely generated groups $\Ga$ equipped with word metrics we will also sometimes use the notion of {\em discrete geodesics}, which are isometric maps 
$$
\ga: I\to G 
$$
where $I$ is an interval (possibly infinite) in $\Z$. Similarly, we will be talking about {\em discrete quasigeodesics} whose domains are intervals in $\Z$.

\begin{definition}\label{def:expansion factor}
For a diffeomorphism $\Phi: M\to M$ of a Riemannian manifold $M$ 
we define the {\em expansion factor} $\eps(\Phi,x)$ at a point $x\in M$ as
\begin{equation*}
\eps(\Phi,x) = \inf_{v\in T_xM-\{0\}} \frac{\|d\Phi(v)\|}{\|v\|} = \| (d\Phi_x)^{-1} \|^{-1} .
\end{equation*}
The {\em exponential expansion factor} of $\Phi$ at $x$ is defined as $\log (\eps(\Phi,x))$. 
\end{definition}

\subsection{Coxeter complexes}
\label{sec:cox}

A {\em spherical Coxeter complex} is  a pair $(S,W)$ 
consisting of a unit sphere $S$ in a Euclidean vector space $V$ and 
a finite group $W$ which acts isometrically on $S$ 
and is generated by reflections at hyperplanes. 
We will use the notation $\angle$ for the angular metric on $S$. 
Throughout the paper, we assume that $W$ does not fix a point in $S$ 
and is associated with a root system $R$. 
Spherical Coxeter complexes will occur as {\em model apartments} 
of spherical buildings, 
mostly of Tits boundaries of symmetric spaces, 
and will in this context usually be denoted by $a_{mod}$. 

A {\em wall} $m_{\rho}$ in $S$ 
is the fixed point set of a hyperplane reflection $\rho$ in $W$. 
An {\em half-apartment} in $S$ is 
a closed hemisphere  bounded by a wall. 
A point $\xi\in S$ is called {\em singular} 
if it belongs to a wall and {\em regular} otherwise. 

The action $W\acts S$ determines on $S$ a structure 
as a simplicial complex 
whose facets, called {\em chambers}, 
are the closures of the connected components of
\begin{equation*}
S-\bigcup_{\rho} m_{\rho}
\end{equation*}
where the union is taken over all reflections $\rho$ in $W$. 
We will refer to the simplices in this complex as {\em faces}. 
(If one allows fixed points for $W$ on $S$, 
then $S$ carries only a structure as a cell complex.) 
Codimension one faces of this complex are called {\em panels}. The {\em interior} $\interior(\tau)$ 
of a face $\tau$ is the complement in $\tau$ to the union of walls not containing $\tau$. The interiors  
$\interior(\tau)$ are called {\em open simplices} in $(S,W)$. 
A geodesic sphere in $S$ is called simplicial if it equals an intersection of walls.

The chambers are fundamental domains for the action $W\acts S$. 
We define the {\em spherical model Weyl chamber} 
as the quotient $\si_{mod}=S/W$. 
The natural projection $\theta:S\to\si_{mod}$ 
restricts to an isometry on every chamber.

It is convenient, and we will frequently do so, to identify  $\si_{mod}$ with a chamber in $S$ (traditionally called the {\em positive chamber}). An important elementary property of the chamber $\si_{mod}$ is that its diameter (with respect to the spherical metric) is $\le \pihalf$.

Given a face $\tau_{mod}$ of $\si_{mod}$, we define the subgroup $W_{\tau_{mod}}\subset W$ as the stabilizer of $\tau_{mod}$ in 
$W$. An identification of $\si_{mod}$ with a chamber $\si\subset S$ 
determines a generating set of $W$, 
namely the reflections at the walls bounding $\si_{mod}$, 
and hence a word metric on $W$; 
the longest element with respect to this metric is denoted $w_0$. 
This element sends 
$\si_{mod}$ to the opposite chamber in $S$. 
We say that two points $\xi, \hat\xi\in S$ are {\em Weyl antipodes} if $\hat\xi=w_0 \xi$. 
We define the {\em standard involution} (also known as the {\em Cartan involution}) 
$$
\iota=\iota_S: S\to S$$ 
as the composition $-w_0$. This involution preserves $\si_{mod}$ 
and equals the identity if and only if $-\id_S\in W$ 
because then $w_0=-\id_S$. 

A point $\xi$ in $S$ is called a {\em root} 
if the hemisphere centered at $\xi$ is simplicial, 
equivalently, is bounded by a wall. Every root point $\xi$ defines a certain linear functional $\alpha_\xi\in V^*$ on the Euclidean vector space $V$ containing $S$ as the unit sphere; this linear functional is also called root. The  kernel of $\al_\xi$ is the hyperplane in $V$ which intersects $S$ along the wall defined by $\xi$.  The set of roots in $V^*$  is denoted by $R$ and called the {\em root system} associated 
with the Coxeter complex $(S,W)$.  We refer the reader to \cite{Bourbaki} for details.

\medskip 
Suppose that $S$ is identified with the 
sphere at infinity of a Euclidean space $F$, 
$S\cong\geo F$, 
where $\geo F$ is equipped with the angular metric. 
For a chamber $\si\subset S$ and a point $x\in F$ 
we define the {\em Weyl sector} 
$V(x, \si)\subset F$ 
as the cone over $\si$ with tip $x$, 
that is, 
as the union of rays emanating from $x$ and asymptotic to $\si$. 

After fixing an origin $o\in F$ 
the group $W$ lifts to a group of isometries of $F$ fixing $o$. 
The sectors $V(o,\si)$ are then fundamental domains 
for the action of $W\acts F$. 

We define the {\em euclidean model Weyl chamber} 
as the quotient 
$V_{mod}=F/W$; 
we will also denote it by $\Delta$ or $\De_{euc}$. 
It is canonically isometric to 
the complete euclidean cone over $\si_{mod}$. 
The natural projection $F\to V_{mod}$ 
restricts to an isometry on the sector $V(o,\si)$ 
for every chamber $\si\subset S$. 
Furthermore, 
for a closed subset $\bar A\subset\si_{mod}$
we define $V(0,\bar A)\subset V_{mod}$
as the complete cone over $\bar A$ with tip $0$.
In particular,
a face $\tau_{mod}\subset\si_{mod}$ corresponds
to a face $V(0,\tau_{mod})$ of $V_{mod}$.

We define the {\em $\Delta$-valued distance function} 
or {\em $\Delta$-distance} 
$d_{\De}$ on $F$ by 
\begin{equation*}
d_\Delta(x,y)= proj(y-x)\in\De
\end{equation*}
where $proj: F\to F/W\cong\Delta$ is the quotient map. 
Note the symmetry property
\begin{equation}
\label{eq:symprop}
d_{\De}(x,y)= \iota_{\De} d_{\De}(y,x) 
\end{equation}
where $\iota_{\De}$ denotes the isometric involution of $\De$ 
induced by $\iota_S$. 
The Weyl group is precisely the group of isometries  
for the $\Delta$-valued distance on $F_{mod}$ 
which fix the origin.

\subsection{Hadamard manifolds}
\label{sec:had}

In this section only,
$X$ denotes a Hadamard manifold, 
i.e.\ a simply connected complete Riemannian manifold 
with nonpositive sectional curvature.  We will use the notation $\Isom(X)$ for the full isometry group of $X$. 

Any two points in $X$ are connected by a unique geodesic segment. 
We will use the notation $xy$ for the oriented geodesic 
segment connecting $x$ to $y$ and parameterized with unit speed. 
We will be treating geodesic segments, geodesic rays and complete geodesics as isometric maps of intervals 
to $X$; we sometime will abuse the terminology and identify geodesics and their images. 

We will denote by $\angle_x(y,z)$ 
the angle between the geodesic segments $xy$ and $xz$ at the point $x$. 
For $x\in X$ we let $\Si_xX$ denote 
the {\em space of directions} of $X$ at $x$, 
i.e.\ the unit sphere in the tangent space $T_xX$, 
equipped with the angle metric. 

The {\em ideal} or {\em visual boundary} of $X$, 
denoted $\geo X$, is the set of asymptote classes of geodesic rays in $X$, 
where two rays are {\em asymptotic} 
if and only if they have finite Hausdorff distance. 
Points in $\geo X$ are called {\em ideal points}. 
For $\xi\in \geo X$ and $x\in X$ we denote by $x\xi$ the geodesic ray 
emanating from $x$ and asymptotic to $\xi$, 
i.e.\ representing the ideal point $\xi$. For $x\in X$ we have a natural map
$$
\log_x: \geo X\to \Si_xX
$$
sending $\xi\in \geo X$ to the velocity vector at $x$ of the geodesic ray $x\xi$. 
The {\em cone} or {\em visual topology} on $\geo X$ 
is characterized by the property that all  the maps $\log_x$ are homeomorphisms;  
with respect to this topology, $\geo X$ is homeomorphic to 
the sphere of dimension $\dim(X)-1$. The visual topology extends to $\bar{X}=X\cup \geo X$ as follows: 
A sequence $(x_n)$ converges to an ideal point $\xi\in \geo X$ 
if the sequence of geodesic segments $xx_n$ 
emanating from some (any) base point $x$ 
converges to the ray $x\xi$ pointwise (equivalently, uniformly on compacts in $\R$). 
This topology makes $\bar{X}$ into a closed ball. 
We define the visual boundary of a subset $A\subset X$ 
as the set $\geo A=\bar A\cap\geo X$ 
of its accumulation points at infinity. 

The ideal boundary $\geo X$ carries the natural 
{\em Tits (angle) metric} $\tangle$, 
defined as 
\begin{equation*}
\tangle(\xi,\eta)= \sup_{x\in X} \angle_{x}(\xi,\eta)
\end{equation*}
where $\angle_{x}(\xi,\eta)$ 
is the angle between the geodesic rays $x\xi$ and $x\eta$. 
The {\em Tits boundary} $\tits X$ is the metric space $(\geo X,\tangle)$. 
The Tits metric is lower semicontinuous with respect to the visual topology 
and, accordingly, 
the {\em Tits topology} induced by the Tits metric 
is finer than the visual topology. 
It is discrete if there is an upper negative curvature bound, 
and becomes nontrivial if flat directions occur. 
For instance, 
the Tits boundary of flat $r$-space is the unit $(r-1)$-sphere, 
$\tits\R^r\cong S^{r-1}(1)$. 
An isometric embedding $X\to Y$ of Hadamard spaces 
induces an isometric embedding $\tits X\to\tits Y$ of Tits boundaries. 

A subset $A$ of $\tits X$ is called {\em convex} if for any two points $\xi, \eta\in A$ with $\tangle(\xi,\eta)<\pi$, the (unique) geodesic 
$\xi\eta$ connecting $\xi$ and $\eta$ in $\tits X$ is entirely contained in $A$.

\medskip
Let $\xi\in \geo X$ be an ideal point. 
For a unit speed geodesic ray $\rho:[0,+\infty)\to X$ 
asymptotic to $\xi$ 
one defines the {\em Busemann function} $b_\xi$ on $X$ as
the uniform monotonic limit
\begin{equation*}
b_\xi(x)=\lim_{t\to+\infty} (d(x, \rho(t))-t).
\end{equation*}
Altering the ray $\rho$ changes $b_\xi$ by an additive constant. 
Along the ray we have  
\begin{equation*}
b_\xi(\rho(t))=-t. 
\end{equation*}
The Busemann function $b_\xi$ is convex, 1-Lipschitz 
and measures the relative distance from $\xi$. 
The sublevel sets 
\begin{equation*}
Hb_{\xi,x}=\{b_{\xi}\leq b_{\xi}(x)\}\subset X
\end{equation*}
are called (closed) {\em horoballs centered at $\xi$}. 
Horoballs are convex. 
The ideal boundaries of horoballs are 
the closed $\pihalf$-balls at infinity 
with respect to the Tits metric, 
\begin{equation*}
\geo Hb_{\xi,x}= \{\tangle(\xi, \cdot)\le \pi/2\}\subset\geo X .
\end{equation*}
Busemann functions are {\em asymptotically linear} along rays; 
if $\gamma:[0,+\infty)\to X$ is a unit speed geodesic ray 
asymptotic to $\eta\in\geo X$, $\ga(+\infty)=\eta$,
then 
\begin{equation*}
\lim_{t\to\infty} \frac{b_\xi(\gamma(t))}{t}=-\cos \tangle(\xi, \eta). 
\end{equation*}
This limit is called the {\em asymptotic slope} of $b_\xi$ at $\eta$. 
In particular, 
rays asymptotic to $\eta$ enter horoballs centered at ideal points $\xi$ 
with $\tangle(\xi,\eta)<\pihalf$. 

\subsection{Symmetric spaces of noncompact type}
\label{sec:symm}

The standard references for this and the following section are \cite{Eberlein} and \cite{Helgason}. Our treatment of this standard material is more geometric than the one presented in these books. 

A symmetric space, denoted by $X$ throughout this paper, 
is said to be of {\em noncompact type} 
if it is nonpositively curved, simply connected 
and has no Euclidean factor. 
In particular, it is a Hadamard manifold. 
We will identify $X$ with the quotient $G/K$ 
where $G$ is a semisimple Lie group 
acting isometrically and transitively on $X$, 
and $K$ is a maximal compact subgroup of $G$.  
We will assume that $G$ is commensurable with the isometry group $\Isom(X)$ 
in the sense that we allow finite kernel and cokernel 
for the natural map $G\to\Isom(X)$. 
In particular, 
the image of $G$ in $\Isom(X)$ contains the identity component $\Isom(X)_o$. 
The Lie group $G$ carries a natural structure as a real algebraic group. 

A {\em point reflection} (confusingly, also known as a {\em Cartan involution}) at a point $x\in X$ 
is an isometry $\si_x$ which fixes $x$ and has differential $-\id_{T_xX}$ in $x$. 
In a symmetric space,
point reflections exist in all points (by definition).
A {\em transvection} of $X$ is an isometry 
which is the product $\si_x\si_{x'}$ of two point reflections; 
it preserves the oriented geodesic through $x$ and $x'$ 
and the parallel vector fields along it. 
The transvections preserving a unit speed geodesic $c(t)$ 
form a one parameter subgroup $(T^c_t)$ of $\Isom(X)_o$ 
where $T^c_t$ denotes the transvection mapping 
$c(s)\mapsto c(s+t)$. 
A nontrivial isometry $\phi$ of $X$ is called {\em axial} 
if it preserves a geodesic $l$ and shifts along it. 
(It does not have to be a transvection.)
The geodesic $l$ is called an {\em axis} of $\phi$. 
Axes are in general not unique. 
They are parallel to each other. 

A {\em flat} in $X$ is a totally geodesic flat submanifold, 
equivalently, 
a convex subset isometric to a Euclidean space. 
A maximal flat in $X$ is a flat 
which is not contained in any larger flat; 
we will use the notation $F$ for maximal flats. 
The group $\Isom(X)_o$ acts transitively on the set of maximal flats; 
the common dimension of maximal  flats is called the {\em rank} of $X$. 
The space $X$ has rank one if and only if it has 
strictly negative sectional curvature. 

A maximal flat $F$ 
is preserved by all transvections along geodesic lines contained in it. 
In general, there exist nontrivial isometries of $X$ fixing $F$ pointwise. 
The subgroup of isometries of $F$ 
which are induced by elements of $G$ 
is isomorphic to a semidirect product $ \R^r \rtimes W$, 
where $r$ is the rank of $X$. 
The subgroup $\R^r$ acts simply transitively on $F$ by translations. 
The linear part $W$ 
is a finite reflection group, called the {\em Weyl group} of $G$ and $X$. 
Since maximal flats are equivalent modulo $G$, 
the action $W\acts F$ is well-defined up to isometric conjugacy. 

We will think of the Weyl group as 
acting on a {\em model flat} $F_{mod}\cong \R^r$ 
and on its ideal boundary sphere at infinity, 
the {\em model apartment} $a_{mod}=\tits F_{mod}\cong S^{r-1}$. 
The pair $(a_{mod},W)$ is the {\em spherical Coxeter complex} 
associated with $X$. 
We identify the spherical model Weyl chamber $\si_{mod}$ 
with a (fundamental) chamber in the model apartment, 
$\si_{mod}\subset a_{mod}$. 
Accordingly, 
we identify the {\em euclidean model Weyl chamber} $V_{mod}$ 
with the sector in $F_{mod}$ 
with tip in the origin and ideal boundary $\si_{mod}$, 
$V_{mod}\subset F_{mod}$. 

The $\Delta$-valued distance naturally extends from $F_{mod}$ to $X$ 
because every pair of points lies in a maximal flat. 
In order to define the distance $d_\Delta(x,y)$ of two points $x,y\in X$ 
one chooses a maximal flat $F$ containing $x,y$ 
and identifies it isometrically with $F_{mod}$ 
while preserving the types of points at infinity. 
The resulting quantity $d_\Delta(x,y)$ is independent of the choices. 
We refer the reader to \cite{KLM} for the detailed discussion of 
{\em metric properties} of $d_\Delta$.

For every maximal flat $F\subset X$, 
we have a Tits isometric embedding $\geo F\subset\geo X$ 
of its ideal boundary sphere. 
There is an identification $\geo F\cong a_{mod}$ 
with the model apartment, 
unique up to composition with elements in $W$. 
The Coxeter complex structure on $a_{mod}$ 
induces a simplicial structure on $\geo F$. 
The ideal boundaries of maximal flats cover $\geo X$ 
because every geodesic ray in $X$ is contained in a maximal flat.
Moreover, their intersections are simplicial. 
One thus obtains a $G$-invariant 
piecewise spherical {\em simplicial structure} on $\geo X$ 
which makes $\geo X$ into a {\em spherical building} and, 
also taking into account the visual topology, 
into a topological spherical building. 
It is called the {\em spherical} or {\em Tits building} 
associated to $X$. 
The Tits metric is the path metric 
with respect to the piecewise spherical structure. 
We will refer to the simplices as {\em faces}.

The ideal boundaries $\geo F\subset\geo X$ 
of the maximal flats $F\subset X$ 
are precisely the {\em apartments} 
with respect to the spherical building structure at infinity, 
which in turn are precisely the convex subsets 
isometric to the unit $(r-1)$-sphere 
with respect to the Tits metric. 
Any two points in $\geo X$ lie in a common apartment.

The action $G\acts\geo X$ on ideal points is not transitive if $X$ has rank $\ge 2$. 
Every $G$-orbit meets every chamber exactly once. 
The quotient can be identified with the spherical model chamber, 
$\geo X/G\cong\si_{mod}$. 
We call the projection 
\begin{equation*}
\theta:\geo X\to\geo X/G \cong\si_{mod}
\end{equation*}
the {\em type} map. 
It restricts to an isometry on every chamber $\si\subset\geo X$. 
We call the inverse 
$\kappa_{\si}=(\theta|_{\si})^{-1}: \si_{mod}\to\si$
the {\em (chamber) chart} for $\si$. 
Consequently, 
$\theta$ restricts to an isometry on every face $\tau\subset\geo X$. 
We call $\theta(\tau)\subset\si_{mod}$ 
the {\em type} of the face $\tau$ and 
$\kappa_{\tau}=(\theta|_{\tau})^{-1}:\theta(\tau)\to\tau$
its {\em chart}. 
We define the {\em type} of an ideal point $\xi\in\geo X$ 
as its image $\theta(\xi)\in\si_{mod}$. 
A point $\xi\in \geo X$ is called {\em regular} 
if its type is an interior point of $\si_{mod}$, 
and {\em singular} otherwise. 
We denote by $\geo^{reg}X\subset\geo X$ 
the set of regular ideal boundary points. 
A point $\rho\in \tits X$ is said to be of {\em root type} if $\theta(\rho)$ is a root in $\si_{mod}\subset S$. 
Equivalently, the closed $\pihalf$-ball centered at $\rho$ (with respect to the Tits metric) 
is simplicial, i.e.\ is a simplicial subcomplex of $\tits X$. 

A geodesic segment $xy$ in $X$ is called {\em regular} if $x\ne y$ and for the unique geodesic ray $x\xi$ containing $xy$ the point $\xi\in \tits X$ is regular. Equivalently, the vector $d_\Delta(x,y)$ belongs  to the interior of $V_{mod}$. 

Two ideal points $\xi,\eta\in\geo X$ are called {\em antipodal} 
if $\tangle(\xi, \eta)=\pi$. 

We say that two simplices $\tau_1,\tau_2\subset\geo X$ 
are {\em opposite} (or {\em antipodal}) {\em with respect to a point $x\in X$} if $\tau_2=\si_x\tau_1$, 
where $\si_x$ is the reflection at the point $x$. 
We say that two simplices $\tau_1,\tau_2\subset\geo X$ 
are {\em opposite} (or {\em antipodal}) 
if they are opposite simplices in the apartments containing both of them,
equivalently, if every interior point of $\tau_1$ has an antipode in the interior of $\tau_2$ and vice versa,
equivalently, if they are opposite with respect to some point $x\in X$.
Their types are then related by $\theta(\tau_2)=\iota(\theta(\tau_1))$.
We will frequently use the notation $\tau, \hat\tau$ and $\tau_+, \tau_-$ for antipodal simplices. 

One can quantify the antipodality of simplices 
of $\iota$-invariant type $\tau_{mod}\subset\si_{mod}$ as follows: 
Pick an $\iota$-invariant type $\bar\zeta$ in the interior of $\tau_{mod}$, $\iota(\zeta)=\zeta$. 
Given two simplices $\tau_+, \tau_-$ in $\tits X$ of type $\tau_{mod}$ and a point $x\in X$ define 
the {\em $\zeta$-angle} 
\begin{equation}\label{eq:zeta-angle}
{\angle_x}^{\zeta}(\tau_-, \tau_+)= \angle_x(\xi_-, \xi_+),
\end{equation}
where $\xi_\pm\in \tau_\pm$ are such that $\theta(\xi_\pm)=\zeta$. Similarly, define the {\em $\zeta$-Tits angle} 
\begin{equation}\label{eq:zeta-tangle}
{\tangle}^\zeta(\tau_-, \tau_+)= {\angle_x}^{\zeta}(\tau_-, \tau_+),  
\end{equation}
where $x$ belongs to a flat $F\subset X$ such that $\tau_-, \tau_+\subset \tits F$. Then simplices $\tau_\pm$ 
(of the same type) are antipodal iff
$$
\tangle^\zeta(\tau_-, \tau_+)=\pi
$$
for some, equivalently, every, choice of $\zeta$ as above.

A pair of opposite chambers $\si_+$ and $\si_-$ 
is contained in a unique apartment,  
which we will denote by $a(\si_+,\si_-)$; the apartment  $a(\si_+,\si_-)$  is the ideal boundary of a unique maximal flat $F(\si_+, \si_-)$ in $X$.

For a point $x\in X$ and a simplex $\tau\subset\geo X$ 
we define the {\em (Weyl) sector} $V=V(x,\tau)\subset X$ 
as the union of rays $x\xi$ for all ideal points $\xi\in\tau$. 
It is contained in a flat. (Thus, Weyl sectors in $X$ are isometric images of Weyl sectors $V(0, \tau_{mod})\subset V_{mod}$ under isometric embeddings $F_{mod}\to X$.)  More generally, 
for a point $x\in X$ and a closed subset $A\subset\geo X$,
we define the {\em Weyl cone} $V(x,A)$
as the union of all rays $x\xi$ for $\xi\in A$. 
It is, in general, not flat.

\medskip
The stabilizers $B_{\si}\subset G$ of the chambers $\sigma\subset \geo X$ 
are the {\em Borel subgroups} of $G$. 
After identifying the model chamber with a chamber in $\geo X$, 
$\si_{mod}\subset\geo X$, 
we call $B=B_{\si_{mod}}$ the {\em positive} Borel subgroup. 
The group $G$ acts transitively on the set of chambers in $\geo X$, 
which we will then identify with $G/B$, 
the {\em  full flag manifold} of $G$. 
The Borel subgroups are algebraic subgroups of $G$, 
and $G/B$ is a real projective variety. 
The set $\D_FX\cong G/B$ of chambers in $\geo X$ 
is called the {\em F\"urstenberg boundary} of $X$; 
we will equip it with the visual topology 
(as opposed to the Zariski topology coming from $G/B$) 
which coincides with its manifold topology 
as a compact homogeneous $G$-space. 
After picking an interior point $\xi \in\interior(\si_{mod})$, 
we can identify the $G$-orbit $G\xi\subset\geo X$ 
$G$-equivariantly and homeomorphically with $\D_FX$ 
by assigning to the (regular) point $g\xi$ 
the (unique) chamber $g\si_{mod}$ containing it. 

The stabilizers $P_{\tau}\subset G$ of simplices $\tau\subset\geo X$ 
are the {\em parabolic subgroups} of $G$. 
The group $G$ acts transitively on simplices of the same type. 
The set 
$\Flag(\tau_{mod})\cong G/P_{\tau_{mod}}$ 
of the simplices $\tau$ of type $\theta(\tau)=\tau_{mod}\subset\si_{mod}$ 
is called the {\em  partial flag manifold}. 
In particular, 
$\Flag(\si_{mod})=\D_FX$. 
Again, we equip the flag manifolds with the visual topology; 
it agrees with their topology as compact homogeneous $G$-spaces. 

We also fix a Riemannian metric on each flag manifold $\Flag(\tau_{mod})$; the particular choice of the metric will be irrelevant.

For a flag manifold $\Flag(\tau_{mod})$ and a simplex $\hat\tau$ of type $\iota(\tau_{mod})$ 
we define the {\em open Schubert stratum} $C(\hat\tau)\subset \Flag(\tau_{mod})$ as the subset of simplices opposite to $\hat\tau$. It follows from semicontinuity of the Tits distance that 
the subset $C(\hat\tau)\subset \Flag(\tau_{mod})$ is indeed open. Furthermore, this subset is also 
dense in $\Flag(\tau_{mod})$. We note that for rank 1 symmetric spaces, 
the only flag manifold associated to $G$ is $\geo X$ 
and the open Schubert strata are the complements of points.

\subsection{Parallel sets, cones, horocycles and decompositions of symmetric spaces}

\subsubsection{Parallel sets}
\label{sec:parset}

Let $s\subset\tits X$ be an isometrically embedded unit sphere. 
We denote by $P(s)\subset X$ the {\em parallel set} associated to $s$. 
It can be defined as the 
union of maximal flats $F\subset X$ asymptotic to $s$, 
$s\subset\geo F$. 
Alternatively, 
one can define it as the 
union of flats $f\subset X$ with ideal boundary $\geo f=s$. 
The parallel set is a totally geodesic subspace 
and splits metrically as the product 
\begin{equation}
\label{eq:parsplit}
P(s)\cong f\times CS(s)
\end{equation}
of any of these flats and a symmetric space $CS(s)$ 
called its {\em cross section}. 
Accordingly, 
the ideal boundary of the parallel set is a metric suspension 
\begin{equation}
\label{eq:bdparsersusp}
\tits P(s)\cong\tits f\circ\tits CS(s) .
\end{equation}
It coincides with the subbuilding $\B(s)\subset\geo X$ 
consisting of the union of all apartments $a\subset\geo X$ containing $s$,
\begin{equation*}
\B(s)=\geo P(s).
\end{equation*}

It is immediate that parallel sets are nonpositively curved symmetric spaces.
However, they do not have noncompact type as their Euclidean de Rham factors 
are nontrivial. 
The factor $f$ in the splitting (\ref{eq:bdparsersusp})
of the parallel set is then the Euclidean de Rham factor
and the cross section $CS(s)$ has trivial euclidean de Rham factor,
i.e.\ it is a symmetric space of noncompact type.

For a pair of antipodal simplices $\tau_+, \tau_-\subset\geo X$
there exists a unique minimal singular sphere $s(\tau_-, \tau_+)\subset \geo X$ containing them.
We denote $P(\tau_-, \tau_+):=P(s(\tau_-, \tau_+))$; this parallel set is the union of (maximal) flats $F\subset X$ whose ideal boundaries contain $\tau_-\cup \tau_+$. In order to simplify the notation, we will denote $\B(s(\tau_-, \tau_+))$ simply by $\B(\tau_-, \tau_+)$.

\subsubsection{Stars, Weyl cones and diamonds}
\label{sec:star}

\begin{dfn}[Stars]
\label{dfn:star}
Let $\tau\subset\tits X$ be a simplex. 
We define the {\em star} $\st(\tau)$ of the open simplex $\interior(\tau)$ as the subcomplex of $\tits X$ 
consisting of all simplices intersecting the open simplex $\interior(\tau)$ nontrivially (i.e., containing $\tau$). In other words, $\st(\tau)$ is the smallest subcomplex of $\tits X$ containing all chambers $\si$ such that $\tau\subset \si$. 
We note that $\st(\tau)$ is also known as the {\em residue} of $\tau$, see e.g. \cite{AB}. 

We define the {\em open star} $\ost(\tau)\subset\geo X$ 
as the union of all open simplices  whose closure intersects $\interior(\tau) $ nontrivially. For the model simplex $\tau_{mod}$, we will use the notation 
$\ost(\tau_{mod})$ to denote its open star in the simplicial complex consisting of faces of $\si_{mod}$. 
\end{dfn}

Note that $\ost(\tau)$ is an open subset of the simplex $\si_{mod}$; it 
does not include any open faces of $\tau$ except for the interior of $\tau$. Furthermore, $\D\st(\tau)=\st(\tau)-\ost(\tau)$ is the union of all panels $\pi$ of type $\theta(\pi)\not\supset\tau_{mod}$ which are contained in a chamber with face $\tau$. 

\medskip 
{\bf $\tau_{mod}$-regularity.} Using the notion of open stars we now generalize the 
standard notion of regularity in $\tits X$ to {\em regularity relative to faces $\tau_{mod}$ of $\si_{mod}$}. 
Let $\tau_{mod}$ be a face of $\si_{mod}$. A point $\xi\in \tits X$ such that  $\theta(\xi)\in \ost(\tau_{mod})$,  is called $\tau_{mod}$-regular. 
If $\Theta\subset \ost(\tau_{mod})$ is a compact subset, then we will refer to points $\xi\in \tits X$ such that $\theta(\xi)\in \Theta$, as being {\em $\Theta$-regular}. Analogously to the definition of regular geodesic segments, a nondegenerate geodesic segment $xy$ is called 
{\em $\tau_{mod}$-regular}, resp. {\em $\Theta$-regular} if it is contained in a geodesic ray $x\xi$ with $\xi$ 
$\tau_{mod}$-regular, resp. {\em $\Theta$-regular}. 

We note that each subset $\Theta$ as above determines the simplex $\tau_{mod}$, namely, $\tau_{mod}$ is the smallest face of $\si_{mod}$ 
such that  $\Theta\subset \ost(\tau_{mod})$.

The fact that the diameter of $\si_{mod}$ is $\le \pihalf$ immediately implies: 

\begin{lem}\label{lem:half}
For every interior point $\xi\in \tau_{mod}$ and every compact subset $\Theta\subset \ost(\tau_{mod})$ we have
$$
\Theta\subset B(\xi,\pihalf). 
$$
\end{lem}

{\bf An isolation property in subbuildings $\B(s)$.}

\begin{lem}
\label{lem:poleisol}
Let $s\subset\geo X$ be a singular sphere.
Then the following hold:

(i)
The simplices contained in $s$ are isolated among the simplices contained in $\B(s)=\geo P(s)$.
In other words: If $\tau_n\to\tau$ is a converging sequence of simplices $\tau_n\subset \B(s)$ with limit simplex $\tau\subset s$, 
then $\tau_n=\tau$ for large $n$.

(ii)
For a simplex $\tau\subset s$, the subset $\ost(\tau)\cap \B(s)$ is open in $\B(s)$ with respect to the visual topology.
In particular,
if $\tau$ is top-dimensional in $s$, then $\ost(\tau)$ is an open subset of $\B(s)$.
\end{lem}
\proof
In view of the decomposition (\ref{eq:bdparsersusp}) of $\B(s)$ as a metric suspension, 
a point $\xi\in s$ has a unique antipode $\hat\xi$ in $\B(s)$, and this antipode is contained in $s$, $\hat\xi\in s$.
Furthermore, 
\begin{equation}
\label{eq:sumdistpole}
\tangle(\xi,\cdot)+\tangle(\cdot,\hat\xi)\equiv\pi
\end{equation} 
on $\B(s)$.
We recall that,
due to the lower semicontinuity of the Tits metric, 
Tits closed balls $\bar B(\eta,r)$ in $\geo X$ are closed also with respect to the visual topology. 
Therefore 
(\ref{eq:sumdistpole}) implies that 
Tits open balls $B(\xi,r)\cap \B(s)$ in $\B(s)$ centered at $\xi\in s$
are open in $\B(s)$ also with respect to the visual topology,
because $B(\xi,r)\cap \B(s)=\B(s)-\bar B(\hat\xi,\pi-r)$.

We recall furthermore that,
due to the finiteness of possible Tits distances between points of a fixed type $\bar\xi\in\si_{mod}$,
the subset $\theta^{-1}\theta(\xi)$ is Tits discrete, 
i.e.\  $B(\xi,\eps)\cap\theta^{-1}\theta(\xi)=\{\xi\}$ for sufficiently small $\eps>0$,
depending on $\theta(\xi)$.
It follows that $\xi$ is an {\em isolated} point of $\theta^{-1}\theta(\xi)\cap \B(s)$
with respect to the visual topology,
i.e.\ it is isolated among points in $\B(s)$ of the same type.

Now we prove (i).
Since the face type of the $\tau_n$ must stabilize,
we may assume that it is constant, i.e.\ $\theta(\tau_n)=\tau_{mod}\subset\si_{mod}$ for all $n$
and $\theta(\tau)=\tau_{mod}$.
We choose a type $\bar\xi\in\interior(\tau_{mod})$ 
and consider the points $\xi_n\in\interior(\tau_n)$ and $\xi\in\interior(\tau)$ 
of type $\theta(\xi_n)=\theta(\xi)=\bar\xi$.
Then $\xi_n\to\xi$.
Since $\xi$ is isolated in $\theta^{-1}\theta(\xi)\cap \B(s)$, it follows that $\xi_n=\xi$ and hence $\tau_n=\tau$ for large $n$.

To verify (ii), we note that 
the subset $\st(\tau)\cap \B(s)$ is the union of all chambers $\si$ with $\tau\subset\si\subset \B(s)$. 
This family of chambers is obviously closed and thus compact in the space of all chambers contained in $\B(s)$.
It is also {\em open} in it as a consequence of part (i).
Indeed, if $\si_n\to\si$ is a converging sequence of chambers $\si_n\subset \B(s)$,
then $\tau_n\to\tau$ for the faces $\tau_n\subset\si_n$ of type $\theta(\tau_n)=\theta(\tau)$,
and it follows that $\tau_n=\tau$ and hence $\tau\subset\si_n$ for large $n$.

Consider now a sequence of points $\eta_n\in \B(s)$ 
such that $\eta_n\to\eta\in\ost(\tau)$.
We must show that $\eta_n\in\ost(\tau)$ for large $n$.
Let $\si_n\subset \B(s)$ be chambers containing them, $\eta_n\in\si_n$. 
It suffices to check that $\si_n\subset\st(\tau)$ for large $n$,
because $\theta(\eta_n)\to\theta(\eta)$, and hence $\theta(\eta_n)\in\ost(\theta(\tau))\subset\si_{mod}$ for large $n$.
Suppose the contrary, i.e.\ that one can pass to a subsequence such that 
$\si_n\not\subset\st(\tau)$ for all $n$.
We may assume that $\si_n\to\si\subset \B(s)$.
By the openness above, it follows that $\si\not\subset\st(\tau)$, equivalently $\tau\not\subset\si$.
On the other hand, $\eta\in\si$.
Together with $\eta\in\ost(\tau)$ this implies that $\tau\subset\si$, 
a contradiction. 
\qed

\medskip 
{\bf Convexity of stars.} 

\begin{lem}
\label{lem:starintapt}The star $\st(\bar\tau)$ of a simplex $\bar\tau\subset a_{mod}$ 
is a convex subset of $a_{mod}$. Furthermore,  $\st(\bar\tau)$ equals the intersection of the simplicial hemispheres $\bar h\subset a_{mod}$ 
such that  $\interior(\bar\tau)\subset\interior{\bar h}$. 
\end{lem}
\proof
If a hemisphere $\bar h$ contains a simplex $\bar\tau$, but does not contain it in its boundary, 
then all chambers containing this simplex as a face belong to the (closed) hemisphere. 
Vice versa, if a chamber $\bar\si$ does not contain $\bar\tau$ as a face, 
then there exists a wall which separates $\bar\si$ from $\bar\tau$. 
\qed

\medskip
Similarly, the star $\st(\tau)$ of a simplex $\tau\subset\tits X$ 
is a convex subset of $\tits X$. 
One can represent it as the intersection of all simplicial $\pihalf$-balls 
which contain $\interior(\tau)$ in their interior. One can represent $\st(\tau)$ also as the intersection of fewer balls:

\begin{lem}[Convexity of stars]
\label{lem:starconv}
(i) 
Let $\tau\subset\tits X$ be a simplex. 
Then $\st(\tau)$ equals the intersection of the simplicial $\pihalf$-balls 
whose interior contains $\interior(\tau)$. 

(ii)
For any simplex $\hat\tau$ opposite to $\tau$, 
the star $\st(\tau)$ equals the intersection 
of the subbuilding 
$\B(\tau,\hat\tau)=\geo P(\tau,\hat\tau)$ 
with all simplicial $\pihalf$-balls 
whose interior contains $\interior(\tau)$ 
and whose center lies in this subbuilding. 
\end{lem}
\proof
(i) If a simplicial $\pihalf$-ball contains a simplex $\tau$, but does not contain it in its boundary, 
then all chambers containing this simplex as a face belong to this ball. 
Vice versa, let $\si$ be a chamber which does not contain $\tau$ as a face. 
There exists an apartment $a\subset\tits X$ which contains $\si$ and $\tau$. 
As before in the proof of Lemma~\ref{lem:starintapt}, 
there exists a simplicial hemisphere $h\subset a$ containing $\tau$ but not $\si$. 
Then the simplicial $\pihalf$-ball with the same center as $h$ contains $\tau$ but not $\si$. 

(ii) Note first that $\st(\tau)\subset \B(\tau,\hat\tau)$. 
Then we argue as in part (i), observing that 
if $\si\subset \B(\tau,\hat\tau)$ then $a$ can be chosen inside $\B(\tau,\hat\tau)$. 
\qed

\medskip
{\bf Convexity of cones.} 

Now we prove a corresponding convexity statement in the symmetric space, 
namely that the Weyl cones $V(x,\st(\tau))$ are convex. We begin with 

\begin{lemma}
For every $x\in P(\tau,\hat\tau)$,  
the Weyl cone $V(x,\st(\tau))$ is contained in the parallel set 
$P(\tau,\hat\tau)$.
\end{lemma}
\proof Consider a chamber $\si$ in $\tits X$ containing $\tau$.  The Weyl sector $V(x,\si)$ is contained in a (unique maximal) flat $F\subset X$. 
Since $\tau, \hat \tau$ are antipodal with respect to $x$, $\tau\cup \hat\tau\subset \geo F$. Therefore, $F\subset P(\tau,\hat\tau)$. \qed

\begin{prop}[Convexity of Weyl cones]
\label{prop:wconeconv}
Let $\hat\tau$ be the simplex opposite to $\tau$ with respect to $x$. 
Then the Weyl cone $V(x,\st(\tau))$ is the intersection of the parallel set $P(\tau,\hat\tau)$
with the horoballs which are centered at $\geo P(\tau,\hat\tau)$ and contain $V(x,\st(\tau))$. 
\end{prop}
\proof
One inclusion is clear. 
We must prove that each point $y\in P(\tau,\hat\tau) \setminus V(x,\st(\tau))$ 
is not contained in one of these horoballs. 
There exists a maximal flat $F\subset P(\tau,\hat\tau)$ containing $x$ and $y$. 
(Any two points in a parallel set lie in a common maximal flat.)
We extend the oriented segment $xy$ to a ray $x\eta$ inside $F$. 

As in the proof of Lemma~\ref{lem:starconv}, 
there exists $\zeta\in\geo F$ such that 
$\bar B(\zeta,\pihalf)$ contains $\st(\tau)$ 
but does not contain $\eta$. 
Then the horoball $Hb_{\zeta,x}$ intersects $F$ in a half-space 
which contains $x$ in its boundary hyperplane but does not contain $\eta$ in its ideal boundary. 
Therefore does not contain $y$.
By convexity, 
$V(x,\st(\tau))\subset Hb_{\zeta,x}$. 
\qed

\medskip
We will need a version of the above convexity results 
for more general stars and cones. 

\begin{dfn}[Weyl convex]
A subset $\Theta\subseteq\si_{mod}$ is called {\em $\tau_{mod}$-Weyl convex}, 
if its symmetrization $W_{\tau_{mod}}\Theta\subseteq\st(\si_{mod})$ 
is a convex subset of $a_{mod}$. 
\end{dfn}
Below, 
$\Theta$ will always denote a $\tau_{mod}$-Weyl convex subset of $\si_{mod}$. 

\medskip 
{\bf Weyl convexity implies convexity.} 

Note that $\si_{mod}$ itself is Weyl convex, since $W_{\tau_{mod}}\si_{mod}=\st(\tau_{mod})$ is convex. 
A point $\bar\xi\in\si_{mod}$ is a Weyl convex subset if and only if it belongs to $\tau_{mod}$. 
\begin{dfn}[$\Theta$-stars]
We define the {\em $\Theta$-star} of a simplex $\tau\subset\geo X$ of type $\tau_{mod}$ as 
$\st_{\Theta}(\tau)=\st(\tau)\cap\theta^{-1}(\Theta)$. 
\end{dfn}
The symmetrization $W_{\tau_{mod}}\Theta\subset a_{mod}$ 
equals the $\Theta$-star of $\tau_{mod}$ inside the model apartment, 
$\st_{\Theta}^{a_{mod}}(\tau_{mod})=W_{\tau_{mod}}\Theta$. 

Our next result establishes convexity of $\Theta$-stars. More precisely:

\begin{lem}[Convexity of $\Theta$-stars]
\label{lem:thetastarconv}
For every simplex $\tau\subset\geo X$ of type $\tau_{mod}$, 
the $\Theta$-star $\st_{\Theta}(\tau)\subset\geo X$ 
equals the intersection of all $\pihalf$-balls containing it. 
\end{lem}
\proof
Let $\zeta\in\geo X$. 
For an apartment $a\subset\geo X$ containing $\tau$ and $\zeta$, 
we have $a\cap\st_{\Theta}(\tau)=\st_{\Theta}^a(\tau)$, 
the intrinsic star of $\tau$ inside the apartment $a$
(viewed as a thin building).
Moreover, 
$\st_{\Theta}(\tau)$ equals the union of all these intersections with apartments. 

By definition of spherical buildings, 
for any two apartments $a,a'\subset\geo X$ containing $\tau$ and $\zeta$, 
there exists an isometry $a\to a'$ fixing (the convex hull of) $\tau$ and $\zeta$. 
This identification carries $\st_{\Theta}(\tau)\cap a$ to $\st_{\Theta}(\tau)\cap a'$. 
Hence $\st_{\Theta}(\tau)\subset B(\zeta,\pihalf)$ 
if and only if $\st_{\Theta}(\tau)\cap a\subset B(\zeta,\pihalf)$ for one (any) of these apartments $a$. 
This shows that the intersection of all $\pihalf$-balls containing $\st(\tau)$ 
is not strictly bigger than $\st(\tau)$, 
since any point in $\geo X$ lies in a common apartment with $\tau$. 
\qed

\begin{prop}[Convexity of $\Theta$-cones]
\label{prop:thconeconv}
For every point $x\in X$ and every simplex $\tau\subset\geo X$ of type $\tau_{mod}$ 
the $\Theta$-cone $V(x,\st_{\Theta}(\tau))$ is convex of $X$. 
\end{prop}
\proof
In view of Lemma~\ref{lem:thetastarconv}, 
the proof of Proposition~\ref{prop:wconeconv} goes through. 
\qed

\medskip
The following consequence will be important for us.
\begin{cor}[Nested $\Theta$-cones]
\label{cor:nestcone}
If $x'\in V(x,\st_{\Theta}(\tau))$, then 
$V(x',\st_{\Theta}(\tau))\subset V(x,\st_{\Theta}(\tau))$. 
\end{cor}

Let $xy\subset X$ be an oriented $\tau_{mod}$-regular geodesic segment. 
Then we define the simplex $\tau=\tau(xy)\subset\geo X$ as follows: 
Forward extend the segment $xy$ to the geodesic ray $x\xi$, 
and let $\tau$ be the unique face of type $\tau_{mod}$ of $\tits X$ such that $\xi\in \st(\tau)$.

\begin{dfn}[Diamond]
We define the {\em $\Theta$-diamond} of a $\Theta$-regular segment $x_-x_+$ as 
\begin{equation*}
\diamo_{\Theta}(x_-,x_+)=
V(x_-,\st_{\Theta}(\tau_+))\cap V(x_+,\st_{\Theta}(\tau_-))
\subset P(\tau_-,\tau_+)
\end{equation*}
where $\tau_{\pm}=\tau(x_{\mp}x_{\pm})$. 
\end{dfn}

The next result follows immediately from Corollary~\ref{cor:nestcone}.
\begin{lem}
If $x'_-x'_+$ is a $\Theta$-regular segment contained in 
$\diamo_{\Theta}(x'_-,x'_+)$
and if $\tau(x'_{\mp}x'_{\pm})=\tau(x_{\mp}x_{\pm})$, 
then 
$\diamo_{\Theta}(x'_-,x'_+)\subset\diamo_{\Theta}(x_-,x_+)$. 
\end{lem}

\subsubsection{Strong asymptote classes}
\label{sec:strongasy}

Recall that two unit speed geodesic rays $\rho_1,\rho_2:[0,+\infty)\to X$ in a Hadamard manifold 
are called {\em asymptotic} if the convex function 
$t\mapsto d(\rho_1(t),\rho_2(t))$ on $[0,+\infty)$ is bounded, 
and they are called {\em strongly asymptotic} if 
$d(\rho_1(t),\rho_2(t))\to 0$ as $t\to+\infty$. 
In the case when $X$ is a symmetric space, one verifies using Jacobi fields on $X$ 
that the decay is exponential with rate depending on the type of the ideal boundary point 
$\rho_1(+\infty)=\rho_2(+\infty)$ (see \cite{Eberlein}). 

We generalize these two notions to sectors. 
A sector $V(x,\tau)$ has a canonical isometric parametrization 
by a {\em sector chart} 
$\kappa_{V(x,\tau)}: V(0,\theta(\tau))\to V(x,\tau)$ 
preserving types at infinity; 
here, $\theta(\tau)\subset\si_{mod}$ is the type of the simplex $\tau$ and  $ V(0,\theta(\tau))$ is a face of the model sector $V_{mod}$. 

For two sectors $V(x_1,\tau_1)$ and $V(x_2,\tau_2)$ of the same type, 
$\theta(\tau_1)=\theta(\tau_2)$, 
the distance function from points in $\kappa_{V(x_1,\tau_1)}$ to the sector $\kappa_{V(x_2,\tau_2)}$ 
\begin{equation}
\label{eq:distsect}
d(\kappa_{V(x_1,\tau_1)},\kappa_{V(x_2,\tau_2)}):
V(0,\theta(\tau_1))\to[0,+\infty)
\end{equation}
is convex. 
The two sectors are called {\em asymptotic}, 
if this distance function is bounded, 
equivalently, if they coincide at infinity, $\tau_1=\tau_2$. 
For two asymptotic sectors $V(x_1,\tau)$ and $V(x_2,\tau)$, 
we define 
\begin{equation*}
d_{\tau}(x_1,x_2):= 
\inf d(\kappa_{V(x_1,\tau)},\kappa_{V(x_2,\tau)});
\end{equation*}
this defines a pseudo-metric $d_{\tau}$ on $X$, 
viewed as the set of (tips of) sectors asymptotic to $\tau$. 
\begin{dfn}
The sectors $V(x_1,\tau)$ and $V(x_2,\tau)$ are 
{\em strongly asymptotic} 
if $d_{\tau}(x_1,x_2)=0$. 
\end{dfn}
This is equivalent to the property 
that for some (any) $\xi\in\interior(\tau)$ 
the rays $x_1\xi$ and $x_2\xi$ are strongly asymptotic. 

We denote by 
\begin{equation*}
X_{\tau}^{par}=X/\sim_{d_{\tau}}
\end{equation*}
the {\em space of (parametrized) strong asymptote classes}
of sectors asymptotic to $\tau$. 

We show next that parallel sets represent strong asymptote classes.
For a simplex $\hat\tau$ opposite to $\tau$ 
we consider the restriction 
\begin{equation}
\label{eq:parasy}
P(\tau,\hat\tau)\to X_{\tau}^{par}
\end{equation}
of the projection $X\to X_{\tau}^{par}$. 
We observe that for points $x_1,x_2\in P(\tau,\hat\tau)$
the distance function (\ref{eq:distsect}) is constant
$\equiv d(x_1,x_2)$. 
Hence (\ref{eq:parasy}) is an isometric embedding. 
To see that it is also surjective, 
we need to verify that 
every sector $V(x,\tau)$ is strongly asymptotic to a sector 
$V(x',\tau)\subset P(\tau,\hat\tau)$. 
This follows from the corresponding fact for geodesics:
\begin{lem}
\label{lem:exstrasygeo}
Let $\xi,\hat\xi$ be antipodal ideal points. 
Then every geodesic asymptotic to $\xi$ is strongly asymptotic to 
a geodesic whose other end is asymptotic to $\hat\xi$.
\end{lem}
\proof
Let $c_1(t)$ be a unit speed geodesic forward asymptotic to $\xi$. 
Then $t\mapsto d(c_1(t),P(\xi,\hat\xi))$ 
is convex and bounded on $[0,+\infty)$, 
and hence non-increasing. We claim that the limit 
$$
d:=\lim_{t\to\infty} d(c_1(t),P(\xi,\hat\xi))
$$
equals zero. To see this, 
we choose a unit speed geodesic $c_2(t)$ 
in $P(\xi,\hat\xi)$ forward asymptotic to $\xi$ 
and use the transvections along it to ``pull back'' $c_1$: 
The geodesics 
$c_1^s:=T^{c_2}_{-s}c_1(\cdot+s)$ 
form a bounded family as $s\to+\infty$ 
and subconverge to a geodesic $c_1^{+\infty}$. 
Since the transvections $T^{c_2}_s$ 
preserve the parallel set $P(\xi,\hat\xi)$, 
the distance functions 
$d(c_1^s(\cdot),P(\xi,\hat\xi))=d(c_1(\cdot+s),P(\xi,\hat\xi))$ 
converge locally uniformly on $\R$ and uniformly on $[0,+\infty)$ 
to the constant $d$. 
It follows that the limit geodesic $c_1^{+\infty}$ 
has distance $\equiv d$ from $P(\xi,\hat\xi)$. 
The same argument, applied to $c_2$ instead of the parallel set, 
implies that $c_1^{+\infty}$ is parallel to $c_2$. 
Thus, $d=0$. 

Now we find a geodesic in $P(\xi,\hat\xi)$ strongly asymptotic to $c_1$ 
as follows. 
Let $t_n\to+\infty$. 
Then there exist unit speed geodesics $c'_n(t)$ in $P(\xi,\hat\xi)$ 
forward asymptotic to $\xi$ and 
with $d(c'_n(t_n),c_1(t_n))=d(c_1(t_n),P(\xi,\hat\xi))\to0$. 
This family of parallel geodesics $c'_n$ is a Cauchy family 
and converges to a geodesic in $P(\xi,\hat\xi)$ 
which is strongly asymptotic to $c_1$. 
\qed

\medskip
We conclude: 
\begin{cor}
\label{cor:strasspariso}
The map (\ref{eq:parasy}) is an isometry. 
\end{cor}

Let $\xi\in\geo X$, 
and let $c(t)$ be a geodesic asymptotic to it, $c(+\infty)=\xi$.
We observe that for every $\eta\in\geo P(c)$ the restriction $b_{\eta}\circ c$ is linear,
because there exists a flat $f$ containing $c$ and asymptotic to $\eta$, $\eta\in\geo f$. 

As a consequence, for any two strongly asymptotic geodesics $c_1(t)$ and $c_2(t)$ at $\eta$,
the restricted Busemann functions $b_{\eta}\circ c_i$ coincide for every $\eta\in\st(\xi)\subset\geo P(c_1)\cap\geo P(c_2)$. 
In fact, there is the following criterion for strong asymptoticity:
\begin{lem}
\label{lem:strasycrit}
For geodesics $c_1(t)$ and $c_2(t)$ asymptotic to $\xi\in\geo X$ 
the following are equivalent:

(i) $c_1(t)$ and $c_2(t)$ are strongly asymptotic. 

(ii) $b_{\eta}\circ c_1=b_{\eta}\circ c_2$ for every $\eta\in\st(\xi)$. 

(ii') $b_{\eta}\circ c_1=b_{\eta}\circ c_2$ for every $\eta\in B(\xi,\eps)$ for some $\eps>0$.
\end{lem}
\proof
We may replace the geodesics $c_i$ by parallel ones without changing their strong asymptote classes 
(Lemma~\ref{lem:exstrasygeo})
and hence assume that they lie in one flat.
Then the assertion is clear.
\qed

\subsubsection{Horocycles and horocyclic subgroups}
\label{sec:horocycles}

We fix a simplex $\tau\subset\geo X$, respectively, 
a parabolic subgroup $P_{\tau}\subset G$ 
and discuss various foliations of $X$ naturally associated to it.

We begin with foliations by flats and parallel sets: 
First, 
we denote by ${\mathcal F}_{\tau}$ the partition of $X$ 
into minimal flats asymptotic to $\tau$, 
i.e.\ singular flats $f\subset X$ such that $\tau$ is a top-dimensional simplex in $\geo f$. 
Second, 
any simplex $\hat\tau$ opposite to $\tau$ 
spans together with $\tau$ a singular sphere 
$s(\tau,\hat\tau)\subset\geo X$
and determines the parallel set 
$P(\tau,\hat\tau):=P(s(\tau,\hat\tau))$ 
with cross section 
$CS(\tau,\hat\tau):=CS(s(\tau,\hat\tau))$. 
The parallel sets $P(\tau,\hat\tau)$ for all $\hat\tau$ 
form a partition ${\mathcal P}_{\tau}$ of $X$, 
which is a coarsening of ${\mathcal F}_{\tau}$. 
The parabolic subgroup $P_{\tau}$ preserves both partitions 
and acts transitively on their leaves, 
because it acts transitively on $X$. 
(This in turn can be derived e.g.\ from the transitivity of $G$ 
on maximal flats.)
This implies that these partitions are smooth foliations. 

We describe now preferred identifications of the leaves of these foliations 
by the actions of certain subgroups of $P_{\tau}$. 
Their orbits will be the submanifolds orthogonal and complementary to the foliations,
i.e.\ the integral submanifolds of the distributions normal to the foliations.

The tuple $(b_{\xi})_{\xi\in\Vert(\tau)}$ 
of Busemann functions for the vertices $\xi$ of $\tau$ 
(well-defined up to additive constants) 
provides affine coordinates simultaneously for each of the flats 
in ${\mathcal F}_{\tau}$. 
The Busemann functions at the other ideal points in $\tau$ are linear combinations of these.
The normal subgroup 
\begin{equation}
\label{eq:largehorocgp}
\bigcap_{\xi\in\Vert(\tau)}\Stab(b_{\xi})
=\bigcap_{\xi\in\tau}\Stab(b_{\xi})
\subset P_{\tau}
\end{equation}
acts transitively on the set of these flats and preserves the coordinates; 
it thus provides consistent preferred identifications between them. 
The level sets of $(b_{\xi})_{\xi\in\Vert(\tau)}$ 
are submanifolds orthogonal to the flats in ${\mathcal F}_{\tau}$, 
because the gradient directions of the Busemann functions $b_{\xi}$
at a point $x\in f\in{\mathcal F}_{\tau}$ 
span the tangent space $T_xf$. 
They form a smooth foliation ${\mathcal F}_{\tau}^{\perp}$ 
and are the orbits of the group (\ref{eq:largehorocgp}). 

We define the {\em horocyclic subgroup} at $\tau$ as the {\em normal subgroup} 
\begin{equation}
\label{eq:charhoro}
N_{\tau}=\bigcap_{\xi\in\st(\tau)}\Stab(b_{\xi})
\subset\Fix(\st(\tau))\triangleleft P_{\tau} .
\end{equation}
Note that 
as a consequence of Lemma~\ref{lem:strasycrit},
$N_{\tau}$ preserves the strong asymptote classes of geodesics at all ideal points $\xi\in\ost(\tau)$.

We now give a method for constructing isometries in $N_{\tau}$.

Let $\xi\in\geo X$ be an ideal point,
and let $c(t)$ be a unit speed geodesic forward asymptotic to it,
$c(+\infty)=\xi$.
Consider the one parameter group
$(T_t^c)$ 
of transvections along $c$. 
The $T_t^c$ fix $\geo P(c)$ pointwise 
and shift the Busemann functions $b_{\eta}$ for $\eta\in\geo P(c)$
by additive constants:
\begin{equation*}
b_{\eta}\circ T_t^c-b_{\eta} \equiv -t\cdot\cos\tangle(\eta,\xi)
\end{equation*}

\begin{lem}
\label{lem:constrhiso}
Let $c_1(t)$ and $c_2(t)$ be geodesics forward asymptotic to $\xi$,
which are strongly asymptotic.
Then there exists an isometry $n\in G$ with the properties:

(i) $n\circ c_1=c_2$.

(ii) $n$ fixes $\geo P(c_1)\cap\geo P(c_2)$ pointwise.

(iii) $b_{\eta}\circ n\equiv b_{\eta}$ for all $\eta\in\geo P(c_1)\cap\geo P(c_2)$. 

\no
In particular, $n\in N_{\tau}$ for the simplex $\tau$ spanned by $\xi$, $\xi\in\interior(\tau)$.
\end{lem}
\proof
By our observation above, 
the isometries $T_{-t}^{c_2}\circ T_t^{c_1}$ 
fix $\geo P(c_1)\cap\geo P(c_2)$ pointwise
and preserve the Busemann functions $b_{\eta}$ for all $\eta\in\geo P(c_1)\cap\geo P(c_2)$.
(The geodesics $c_i$ need not have unit speed;
they have the same positive speed since they are strongly asymptotic.)
Moreover, 
they form a bounded family and, 
as $t\to+\infty$, 
subconverge to an isometry $n\in G$ with the same properties 
and which maps $c_1$ to $c_2$ while preserving parametrizations, 
compare the proof of Lemma~\ref{lem:exstrasygeo}.  
The last assertion follows because $\st(\tau)\subset\geo P(c_1)\cap\geo P(c_2)$.
\qed

\begin{cor}
\label{cor:horotranspar}
$N_{\tau}$ acts transitively on 

(i) every strong asymptote class of geodesics at every ideal point $\xi\in\interior(\tau)$;

(ii) the set of leaves of ${\mathcal P}_{\tau}$.
\end{cor}
\proof
Part (i) is a direct consequence. 

To verify (ii), 
we choose $\xi\in\interior(\tau)$. 
Given leaves $P(\tau,\hat\tau_i)$ of ${\mathcal P}_{\tau}$, 
we let $c_i(t)$ be unit speed geodesics in $P(\tau,\hat\tau_i)$
forward asymptotic to $\xi$.
They may be chosen strongly asymptotic by Lemma~\ref{lem:exstrasygeo}. 
The parallel set of $c_i$ equals $P(\tau,\hat\tau_i)$ 
because $\xi$ is an interior point of $\tau$. 
Hence an isometry in $N_{\tau}$ carrying $c_1$ to $c_2$
carries $P(\tau,\hat\tau_1)$ to $P(\tau,\hat\tau_2)$. 
\qed

\begin{rem}
One also obtains 
that every geodesic asymptotic to an ideal point $\xi\in\D\tau$
can be carried by an isometry in $N_{\tau}$ to any other strongly asymptotic geodesic.
However, 
$N_{\tau}$ does not preserve strong asymptote classes at $\xi$ in that case. 
\end{rem}

\begin{lem}
If $n\in N_{\tau}$ 
preserves a leaf $P(\tau,\hat\tau)$, 
$n\hat\tau=\hat\tau$, 
then it acts trivially on it. 
\end{lem}
\proof
Then $n$ fixes $\st(\tau)$ and $\hat\tau$ pointwise, 
and hence their Tits convex hull $\geo P(\tau,\hat\tau)$. 
Thus $n$ preserves every maximal flat $F$ in $P(\tau,\hat\tau)$. 
Moreover it preserves the Busemann functions $b_{\xi}$ 
at all $\xi\in\geo F\cap\st(\tau)$, 
and hence must fix $F$ pointwise,
compare Lemma~\ref{lem:strasycrit}.
\qed

\medskip
Thus, $N_{\tau}$ provides consistent preferred identifications between 
the parallel sets $P(\tau,\hat\tau)$. 
The $N_{\tau}$-orbits are submanifolds orthogonal to the parallel sets. 
They form a smooth foliation 
\begin{equation}
\label{eq:horfol}
{\mathcal H}_{\tau}={\mathcal P}_{\tau}^{\perp}
\end{equation}
refining ${\mathcal F}_{\tau}^{\perp}$, 
which we call the {\em horocyclic foliation} 
and its leaves the {\em horocycles} at $\tau$. 
We denote the horocycle at $\tau$ through the point $x$ by $Hc^{\tau}_x$, i.e.\ $Hc^{\tau}_x=N_{\tau}x$.

For incident faces $\ups\subset\tau$, 
the associated subgroups and foliations 
are contained in each other. 
For instance, we have 
$\st(\ups)\supset\st(\tau)$ and 
$N_{\ups}\subset N_{\tau}$. 

\medskip
We next relate horocycles and strong asymptote classes. 

\begin{prop}[Strong asymptote classes are horocycles]\label{prop:Strong asymptote classes are horocycles}
The sectors $V(x_1,\tau)$ and $V(x_2,\tau)$ are strongly asymptotic 
if and only if $x_1$ and $x_2$ lie in the same horocycle at $\tau$. 
\end{prop}
\proof
Let $\xi\in\interior(\tau)$. 
By Corollary~\ref{cor:horotranspar}(ii), 
$N_{\tau}$ acts transitively on every strong asymptote class 
of geodesics at $\xi$, and hence of sectors at $\tau$. 
Thus, strong asymptote classes are contained in horocycles. 

On the other hand, they cannot be strictly larger 
because every strong asymptote class intersects every parallel set, 
cf.\ the surjectivity of (\ref{eq:parasy}), 
and every horocycle intersects every parallel set 
exactly once. 
\qed

\medskip
Our discussion shows that there is the exact sequence 
\begin{equation*}
N_{\tau}\to P_{\tau}\to\Isom(X_{\tau}^{par}) .  
\end{equation*}

\begin{rem}
Note that the homomorphism $P_{\tau}\to\Isom(X_{\tau}^{par})$ is (in general) not surjective.
Namely, let $f\times CS(s)$ denote the de Rham decomposition with the maximal flat factor $f$. 
Then the image of the above homomorphism (if $\dim(f)>1$) does not contain 
the full group of rotations of $f$. This can be corrected as follows. Let 
$A_\tau$ denote the group of translations of $f$, and let $M_\tau$ be the isometry group of 
$CS(s)$. Then the above exact sequence is a part of the  {\em Langlands' decomposition} of the group $P_\tau$, 
$$
1\to N_\tau \to P_\tau \to A_\tau\times M_\tau \to 1, 
$$
which, on the level of Lie algebras, is a split exact sequence.  
\end{rem}

We return now to Lemma~\ref{lem:constrhiso}.
For later use, 
we elaborate on the special case 
when the geodesics $c_i$ are contained in the parallel set of a 
singular flat of dimension rank minus one.

For an half-apartment $h\subset\geo X$, 
we define its {\em star} $\st(h)$ as the union of the stars $\st(\tau)$
where $\tau$ runs through all simplices with $\interior(\tau)\subset\interior(h)$,
equivalently, which are spanned by an interior point of $h$.
Similarly, we define the {\em open star} $\ost(h)$ as the union of the corresponding open stars $\ost(\tau)$.
Note that $\interior(h)\subset\ost(h)$.
Furthermore, 
we define the closed subgroup $N_h\subset G$ as the intersection of the horocyclic subgroups $N_{\tau}$
at these simplices $\tau$.

We observe that $N_h$ 
preserves the strong asymptote classes of geodesics at all ideal points $\xi\in\ost(h)$,
and it permutes the maximal flats $F$ asymptotic to $h$, $\geo F\supset h$.
The next result shows that it acts transitively on them:

\begin{lem}
\label{lem:constrhisospec}
Let $h\subset\geo X$ be a half-apartment,
and let $F_1,F_2\subset P(\D h)$ be maximal flats asymptotic to $h$, $\geo F_i\supset h$.
Then there exists an isometry $n\in N_h$ with the properties:

(i) $nF_1=F_2$.

(ii) $n$ fixes $\st(h)$ pointwise.

(iii) $b_{\eta}\circ n\equiv b_{\eta}$ for all $\eta\in\st(h)$.
\end{lem}
\proof
In the metric decomposition $P(\D h)\cong\R^{rank(X)-1}\times CS(\D h)$ (see \eqref{eq:parsplit}), 
the maximal flats $F_i$ correspond to a pair of asymptotic, and hence strongly asymptotic geodesics $\bar c_i$ 
in the rank one symmetric space $CS(\D h)$.

Let $\xi,\xi'\in\interior(h)$. 
Let $c_1(t)$ and $c_2(t)$ be strongly asymptotic geodesics at $\xi$ so that $c_i\subset F_i$.
Then they project (up to reparameterization) to the strongly asymptotic geodesics $\bar c_i$ in $CS(\D h)$,
and their projections to the Euclidean de Rham factor of $P(\D h)$ coincide.
Analogously, let $c'_i\subset F_i$ be strongly asymptotic geodesics at $\xi'$.
Their parametrizations can be chosen so that 
their projections to $CS(\D h)$ coincide with the projections of the $c_i$. 
Then 
\begin{equation*}
T_{-t}^{c_2}\circ T_t^{c_1} = T_{-t}^{c'_2}\circ T_t^{c'_1}
\end{equation*}
and the isometry $n$ produced in the proof of Lemma~\ref{lem:constrhiso} belongs to both $N_{\tau}$ and $N_{\tau'}$
where $\tau.\tau'$ denote the simplices spanned by $\xi,\xi'$.
Varying $\xi$ or $\xi'$ yields the assertion. 
\qed

\medskip
We obtain an analogue of Corollary~\ref{cor:horotranspar}:
\begin{cor}
\label{cor:horotransparspec}
$N_h$ acts transitively on 

(i) every strong asymptote class of geodesics at every ideal point $\xi\in\interior(h)$;

(ii) the set of maximal flats $F$ asymptotic to $h$, $\geo F\supset h$.
\end{cor}

We describe a consequence of our discussion for the horocyclic foliations. 

The maximal flats asymptotic to $h$ are contained in the parallel set $P(\D h)\cong\R^{rank(X)-1}\times CS(\D h)$
and form the leaves of a smooth foliation ${\mathcal P}_h$ of $P(\D h)$.
This foliation is the pullback of the one-dimensional foliation of the rank one symmetric space $CS(\D h)$ by the geodesics 
asymptotic to the ideal point $\zeta\in\geo CS(\D h)$ corresponding to the center (pole) of $h$.
We call the foliation of $P(\D h)$ normal to ${\mathcal P}_h$
the {\em horocyclic} foliation ${\mathcal H}_h$.
Its leaves are of the form point times a horosphere in $CS(\D h)$ centered at $\zeta$.
We call them {\em horocycles at $h$}.
Corollary~\ref{cor:horotransparspec} implies that they are the $N_h$-orbits,
and we denote by $Hc^h_x=N_hx$ the horocycle through the point $x\in P(\D h)$.

Let $\tau$ be a simplex so that $\interior(\tau)\subset\interior(h)$.
Then the foliation ${\mathcal P}_{\tau}$ of $X$ by parallel sets
restricts on $P(\D h)$ to the foliation ${\mathcal P}_h$ by maximal flats,
and the horocyclic foliation ${\mathcal H}_{\tau}$
restricts to the horocyclic foliation ${\mathcal H}_h$.
(This follows from the fact that the foliations ${\mathcal P}_{\tau}$ and ${\mathcal H}_{\tau}$
are normal to each other, cf.\ (\ref{eq:horfol}).)
In other words, the foliations ${\mathcal H}_{\tau}$ for the various simplices $\tau$ {\em coincide} 
on the parallel set $P(\D h)$.

\subsubsection{Distances to parallel sets versus angles}
\label{sec:angledist}

In this section we collect further geometric facts regarding 
parallel sets in symmetric spaces, primarily dealing with estimation of distances from points in $X$ to parallel sets. 

 We first strengthen Proposition \ref{prop:Strong asymptote classes are horocycles}:

\begin{lemma}\label{lem:strongly asymptotic geodesic to a star} 
Suppose that $\tau_\pm$ are antipodal simplices in $\tits X$. Then 
every geodesic ray $\ga$ asymptotic to a point $\xi\in \ost(\tau_+)$, 
is strongly asymptotic to a geodesic ray in $P(\tau_-, \tau_+)$. 
\end{lemma} 
\proof If $\xi$ belongs to the interior of the simplex $\tau_+$, then the assertion follows from 
Proposition \ref{prop:Strong asymptote classes are horocycles}. We now consider the general case. Suppose, that $\xi$ belongs to an open simplex 
$\interior(\tau')$, such that $\tau$ is a face of $\tau'$. Then there exists an apartment $a\subset \tits X$ containing both $\xi$ (and, hence, $\tau'$ as well as $\tau$) and the simplex $\tau_-$. Let $F\subset X$ be the maximal flat with $\geo F=a$. Then $F$ contains a geodesic asymptotic to points in 
$\tau_-$ and $\tau_+$. Therefore, $F$ is contained in $P(\tau_-, \tau_+)$. On the other hand, by the same Proposition \ref{prop:Strong asymptote classes are horocycles} applied to the simplex $\tau'$, we conclude that  $\ga$ is strongly asymptotic to a geodesic ray in $F$. \qed 

\medskip 
The following lemma provides a quantitative strengthening of the conclusion of Lemma \ref{lem:strongly asymptotic geodesic to a star}: 

\begin{lem}
\label{lem:decay} 
Let $\Theta$ be a compact subset of $\ost(\tau_+)$. Then 
those rays $x\xi$ with $\theta(\xi)\in\Theta$ 
are uniformly strongly asymptotic to $P(\tau_-,\tau_+)$, 
i.e.\ $d(\cdot,P(\tau_-,\tau_+))$ decays to zero along them 
uniformly in terms of $d(x, P(\tau_-, \tau_+))$ and $\Theta$. 
\end{lem}
\proof Suppose that the assertion of lemma is false, i.e., there exists $\eps>0$, a sequence $T_i\in \R_+$ diverging to infinity, and 
a sequence of  rays  $\rho_i=x_i \xi_i$ with  $\xi_i\in \Theta$ and $d(x_i, P(\tau_-,\tau_+))\le d$, so that
\begin{equation}\label{eq:outside}
d(y, P(\tau_-,\tau_+))\ge \eps, \forall y\in \rho([0, T_i]).
\end{equation}
Using the action of the stabilizer of 
$P(\tau_-,\tau_+)$, we can assume that the points $x_i$ belong to a certain compact subset of $X$. Therefore, the sequence of rays 
$x_i \xi_i$ subconverges to a ray $x\xi$ with $d(x, P(\tau_-,\tau_+))\le d$ and $\xi\in \Theta$. The inequality \eqref{eq:outside} then implies that 
the entire limit ray $x\xi$ is contained outside of the open $\eps$-neighborhood of the parallel set $P(\tau_-,\tau_+)$. However, in view of Lemma 
\ref{lem:strongly asymptotic geodesic to a star}, the ray $x\xi$ is strongly asymptotic to a geodesic in $P(\tau_-,\tau_+)$. Contradiction. \qed

\medskip 
We next relate distance from points $x\in X$ to parallel sets and certain angles at $x$. 
Fix a generic point $\zeta=\zeta_{mod}$ in $\tau_{mod}$ and 
consider a pair of simplices $\tau_{\pm}$ of type $\tau_{mod}$;  
set $\zeta_{\pm}=\zeta(\tau_{\pm})$. We begin with the following elementary observation:

\begin{rem}
\label{rem:antip}
We observe that the ideal points $\zeta_{\pm}$ are opposite,
$\tangle(\zeta_-,\zeta_+)=\pi$, 
if and only if they can be seen under angle $\simeq\pi$ (i.e., close to $\pi$) 
from some point in $X$. More precisely, 
there exists $\eps(\zeta_{mod})$ such that:
 
{\em If $\angle_x(\zeta_-,\zeta_+)>\pi-\eps(\zeta_{mod})$ 
for some point $x$ then $\zeta_{\pm}$ are opposite. }

\noindent This follows from the angle comparison 
$\angle_x(\zeta_-,\zeta_+)\leq\tangle(\zeta_-,\zeta_+)$
and the fact that the Tits distance between ideal points 
of the fixed type $\zeta_{mod}$ 
takes only finitely many values. 
\end{rem}

Suppose now that the simplices $\tau_{\pm}$, equivalently, 
the ideal points $\zeta_{\pm}$, are opposite. 
Then $\angle_x(\zeta_-,\zeta_+)=\pi$ 
if and only if $x$ lies in the parallel set $P(\tau_-,\tau_+)$. 
Furthermore, 
$\angle_x(\zeta_-,\zeta_+)\simeq\pi$
if and only if $x$ is close to $P(\tau_-,\tau_+)$, 
and both quantities control each other near the parallel set. 
More precisely: 

\begin{lem}
\label{lem:distangcontr}
(i)
If $d(x,P(\tau_-,\tau_+))\leq d$, 
then $\angle_x(\zeta_-,\zeta_+)\geq\pi-\eps(d)$ 
with $\eps(d)\to0$ as $d\to0$. 

(ii)
For sufficiently small $\eps$, $\eps\leq\eps'(\zeta_{mod})$, we have:
The inequality $\angle_x(\zeta_-,\zeta_+)\geq\pi-\eps$ implies that  
$d(x,P(\tau_-,\tau_+))\leq d(\eps)$ for some function $d(\eps)$ which converges to $0$ as $\eps\to0$.
\end{lem}
\proof
The intersection of parabolic subgroups 
$P_{\tau_-}\cap P_{\tau_+}$
preserves the parallel set $P(\tau_-,\tau_+)$ and acts transitively on it.
Compactness and the continuity of $\angle_{\cdot}(\zeta_-,\zeta_+)$
therefore imply that $\pi-\angle_{\cdot}(\zeta_-,\zeta_+)$ 
attains on the boundary of the tubular $r$-neighborhood of $P(\tau_-,\tau_+)$ 
a strictly positive maximum and minimum, 
which we denote by $\phi_1(r)$ and $\phi_2(r)$. 
Furthermore, $\phi_i(r)\to0$ as $r\to0$. 
We have the estimate:
\begin{equation*}
\pi-\phi_1(d(x,P(\tau_-,\tau_+)))
\leq\angle_x(\zeta_-,\zeta_+)\leq
\pi-\phi_2(d(x,P(\tau_-,\tau_+)))
\end{equation*}
The functions $\phi_i(r)$ are (weakly) monotonically increasing. 
This follows from the fact that, 
along rays asymptotic to $\zeta_-$ or $\zeta_+$, 
the angle $\angle_{\cdot}(\zeta_-,\zeta_+)$ is monotonically increasing 
and the distance $d(\cdot,P(\tau_-,\tau_+))$ is monotonically decreasing. 
The estimate implies the assertions. 
\qed

\medskip
The control of 
$d(\cdot,P(\tau_-,\tau_+))$  and $\angle_{\cdot}(\zeta_-,\zeta_+)$ 
``spreads'' along the Weyl cone $V(x,\st(\tau_+))$, 
since the latter is asymptotic to the parallel set $P(\tau_-,\tau_+)$. 
Moreover, the control improves,  
if one enters the cone far into a $\tau_{mod}$-regular direction. 
More precisely: 
\begin{lem}
\label{lem:distangdec}
Let $y\in V(x,\st_{\Theta}(\tau_+))$
be a point with $d(x,y)\geq l$. 

(i) 
If $d(x,P(\tau_-,\tau_+))\leq d$,  
then 
\begin{equation*}
d(y,P(\tau_-,\tau_+))\leq D'(d,\Theta,l)\leq d
\end{equation*}
with $D'(d,\Theta,l)\to0$ as $l\to+\infty$. 

(ii)
For sufficiently small $\eps$, $\eps\leq\eps'(\zeta_{mod})$, we have:
If $\angle_x(\zeta_-,\zeta_+)\geq\pi-\eps$, 
then 
\begin{equation*}
\angle_y(\zeta_-,\zeta_+)\geq\pi-\eps'(\eps,\Theta,l)\geq\pi-\eps(d(\eps))
\end{equation*}
with $\eps'(\eps,\Theta,l)\to0$ as $l\to+\infty$. 
\end{lem}
\proof
The distance from $P(\tau_-,\tau_+)$
takes its maximum at the tip $x$ of the cone $V(x,\st(\tau_+))$, 
because it is monotonically decreasing along the rays $x\xi$ 
for $\xi\in\st(\tau_+)$. 
This yields the right-hand bounds $d$ and, 
applying Lemma~\ref{lem:distangcontr} twice, $\eps(d(\eps))$. 

Those rays $x\xi$ 
with uniformly $\tau_{mod}$-regular type $\theta(\xi)\in\Theta$ 
are uniformly strongly asymptotic to $P(\tau_-,\tau_+)$, 
i.e.\ $d(\cdot,P(\tau_-,\tau_+))$ decays to zero along them 
uniformly in terms of $d$ and $\Theta$, see Lemma \ref{lem:decay}.  
This yields the decay $D'(d,\Theta,l)\to0$ as $l\to+\infty$. 
The decay of $\eps'$ follows by applying Lemma~\ref{lem:distangcontr} again. 
\qed

\subsection{Dynamics of transvections at infinity}
\subsubsection{Identifications of horocycles and contraction}

We continue with the notation of section \ref{sec:horocycles}. 

Fix a simplex $\tau\subset\geo X$.
Since every leaf of the foliation ${\mathcal H}_{\tau}$ intersects every leaf of ${\mathcal P}_{\tau}$ exactly once,
and the leaves of the latter foliation correspond to the flags in $C(\tau)$,
we have consistent simultaneous $N_{\tau}$-equivariant smooth identifications  
\begin{equation*}
Hc^{\tau}\buildrel\cong\over\to {Hc'}^{\tau}
\end{equation*}
between the horocycles at $\tau$,
and 
\begin{equation*}
Hc^{\tau}\buildrel\cong\over\to C(\tau)
\end{equation*}
of the horocycles  with the open Schubert stratum $C(\tau)$. 
A point $x\in Hc^{\tau}$ corresponds to a point $x'\in {Hc'}^{\tau}$ 
iff they lie in the same leaf of ${\mathcal P}_{\tau}$, i.e.\ parallel set $P(\tau,\hat\tau)$ 
for a (unique) simplex $\hat\tau$ opposite to $\tau$. 

Let $h\subset\geo X$ be a half-apartment such that $\interior(\tau)\subset\interior(h)$. 
Then the horocycles at $\tau$ intersect the parallel set $P(\D h)$ in the horocycles at $h$,
and the identifications restrict to $N_h$-equivariant smooth identifications 
\begin{equation*}
Hc^h\buildrel\cong\over\to {Hc'}^h
\end{equation*}
between the horocycles at $h$,
and 
\begin{equation*}
Hc^h \buildrel\cong\over\to C(h)
\end{equation*}
of the horocycles with a submanifold $C(h)\subset C(\tau)$.

We discuss now contraction-expansion properties of these identifications.

The horocycles $Hc^h_x$ in $P(\D h)$ isometrically project onto the horospheres $Hc^{\zeta}_{\bar x}$ in $CS(\D h)$. 
Under these projections, 
the natural identifications $\pi^h_{x'x}: Hc^h_x\to Hc^h_{x'}$ correspond to the natural identifications 
$\pi^{\zeta}_{\bar x'\bar x}: Hc^{\zeta}_{\bar x}\to Hc^{\zeta}_{\bar x'}$ 
of horospheres.
The cross sections $CS(\D h)$ are rank one symmetric spaces.
There are only finitely many isometry types of them occurring in $X$.
In a rank one symmetric space, the natural identifications between horospheres contract and expand 
exponentially.
(This follows from the fact that the exponential decay rate of decaying Jacobi fields along geodesics 
is bounded below and above.)
We therefore obtain the following estimate for the contraction-expansion of the identifications 
$\pi^h_{x'x}:Hc^h_x\to Hc^h_{x'}$:
If $b_{\zeta}(x)-b_{\zeta}(x')=b_{\zeta}(\bar x)-b_{\zeta}(\bar x')\geq0$
and $x'=\pi^h_{x'x}(x)$, 
then for corresponding points $y\in Hc^h_x$ and $y'=\pi^h_{x'x}(y)\in Hc^h_{x'}$ we have 
\begin{equation}
\label{ineq:contrestrk1}
e^{-c_1d(\bar x,\bar x')} \leq\frac{d(x',y')}{d(x,y)}\leq e^{-c_2d(\bar x,\bar x')}
\end{equation}
with constants $c_1\geq c_2>0$ depending only on $X$.

\subsubsection{Infinitesimal contraction of transvections}

We describe now the action of transvections at infinity
using the natural identifications of horocycles.

Suppose that $x,x'\in P(\tau,\hat\tau)$ are distinct points. 
Let $\vartheta_{xx'}$ denote the transvection with axis $l_{xx'}$ through $x$ and $x'$ 
mapping $x'\mapsto x$. 
It preserves $\geo  P(\tau,\hat\tau)\subset\geo X$
and fixes the singular sphere $s(\tau,\hat\tau)\subset\geo  P(\tau,\hat\tau)$ pointwise.
In particular, it fixes the simplices $\tau$ and $\hat\tau$.
We consider the action of $\vartheta_{xx'}$ on $C(\tau)$ 
and its differential at the fixed point $\hat\tau$. 
Modulo the identifications $C(\tau)\buildrel\cong\over\to Hc^{\tau}_x$ and 
$\pi^{\tau}_{x'x}:Hc^{\tau}_x\to Hc^{\tau}_{x'}$,
the action of $\vartheta_{xx'}$ on $C(\tau)$ 
corresponds to the action of $\vartheta_{xx'}\pi^{\tau}_{x'x}$ on $Hc^{\tau}_x$, 
and the differential of $\vartheta_{xx'}$ at $\hat\tau$ 
to the differential of $\vartheta_{xx'}\pi^{\tau}_{x'x}$ at $x$. 

We first consider the case that $\vartheta_{xx'}$ translates towards $\hat\tau$,
i.e.\ when $\xi:=l_{xx'}(-\infty)\in\ost(\tau)$,
equivalently, when $x'$ lies in the interior of the Weyl cone $V(x,\st(\tau))$. 
\begin{lem}
\label{lem:diffposev}
If $\xi\in\ost(\tau)$,
then $(d\vartheta_{xx'})_{\hat\tau}$ 
is diagonalizable with strictly positive eigenvalues.
\end{lem}
\proof
If $\xi\in\ost(\tau)$,
then there is a natural $N_{\tau}$-equi\-va\-ri\-ant identification of $C(\tau)$ 
with the strong asymptote class of the geodesic $l_{xx'}$.
Namely, 
the simplex $\hat\tau'\in C(\tau)$ corresponds to the unique geodesic $l'\subset P(\tau,\hat\tau')$ 
strongly asymptotic to $l_{xx'}$,
and in particular $\hat\tau$ corresponds to $l_{xx'}$. 
Accordingly,
tangent vectors to $C(\tau)$ at $\hat\tau$ one-to-one correspond to Jacobi fields along $l_{xx'}$
which are orthogonal to $P(\tau,\hat\tau)$ and decay to zero at $\xi$.
The effect of the differential $(d\vartheta_{xx'})_{\hat\tau}$ on $C(\tau)$ is given in terms of these Jacobi fields $J$ by
\begin{equation*}
J\mapsto d\vartheta_{xx'}\circ J\circ \vartheta_{x'x}
\end{equation*}
The Jacobi fields, which are of the form exponential function times a parallel vector field along $l_{xx'}$,
correspond to eigenvectors of the differential of $\vartheta_{x'x}$ with strictly positive eigenvalues.
It is a standard fact from the Riemannian geometry of symmetric spaces
that every decaying Jacobi field 
(orthogonally) decomposes as a sum of such special Jacobi fields.
(One trivializes the normal bundle along the geodesic using the one parameter group of transvections along it.
Since the curvature tensor of a Riemannian  symmetric space is parallel, 
the Jacobi equation becomes an ODE with constant coefficients.) 
Thus the eigenvectors of $(d\vartheta_{xx'})_{\hat\tau}$ for positive eigenvalues span $T_{\hat\tau}C(\tau)$.
\qed

\begin{lem}
\label{lem:evestinst}
If $\xi\in\ost(\tau)$,
then the eigenvalues $\la$ of $(d\vartheta_{xx'})_{\hat\tau}$ satisfy an estimate
\begin{equation*}
-\log\la\geq c\cdot d(x',\D V(x,\st(\tau)))
\end{equation*}
with a constant $c>0$ depending only on $X$.
\end{lem}
\proof
We continue the argument in the previous proof.

Let $F\supset l_{xx'}$ be a maximal flat. 
Then $F\subset P(\tau,\hat\tau)$. 
A smooth variation of the geodesic $l_{xx'}$ by strongly asymptotic geodesics
extends to a smooth variation of $F$ by maximal flats asymptotic to $\st(\tau)\cap\geo F$. 
Indeed, 
it can be induced by a smooth curve in $N_{\tau}$ through $1$,
and this curve can be used to vary $F$.

The Jacobi field $J$ along $l_{xx'}$ corresponding to a tangent vector $v\in T_{\hat\tau}C(\tau)$
therefore extends to a Jacobi field $\hat J$ along $F$ which decays to zero at all ideal points $\eta\in\ost(\tau)\cap\geo F$.
The decomposition of Jacobi fields on symmetric spaces mentioned in the previous proof
works in the same way along flats.
(One trivializes the normal bundle using the abelian transvection subgroup.)
Hence $\hat J$ decomposes as an orthogonal sum of Jacobi fields along $F$ of the form 
$e^{-\al}V$
with an affine linear form $\al$ on $F$ 
and a parallel orthogonal vector field $V$ along $F$.
Up to additive constants, only finitely many affine linear forms $\al$ occur, 
since $G$ acts transitively on maximal flats.
(The possible forms are determined by the root system of $G$, but we do not need this fact here.)

We may normalize the forms $\al$ occurring in our decomposition so that $\al(x)=0$.
Since $\hat J$ decays to zero at $\ost(\tau)\cap\geo F$, 
they have the property that 
$\al\geq0$ on $V(x,\st(\tau)\cap\geo F)$,
equivalently, 
that $\st(\tau)\cap\geo F\subset\geo\{\al\geq0\}$.
Moreover, $\al>0$ on the interior of the cone $V(x,\st(\tau)\cap\geo F)$, because $\al\not\equiv0$.
One can estimate 
\begin{equation*}
\al(x')\geq c(\al)\cdot d(x',\D V(x,\st(\tau)))
\end{equation*}
with a constant $c(\al)>0$.
Taking $c$ to be the minimum over the finitely many constants $c(\al)$,
we obtain the assertion, 
because the eigenvalues of $(d\vartheta_{xx'})_{\hat\tau}$ are bounded above by the maximal possible value of $e^{-\al(x')}$.
\qed

\medskip
The previous two lemmas yield:
\begin{cor}
\label{cor:wcontrincone}
If $x'\in V(x,\st(\tau))$, 
then $(d\vartheta_{xx'})_{\hat\tau}$ is weakly contracting on $T_{\hat\tau}C(\tau)$ with norm 
\begin{equation*}
\|(d\vartheta_{xx'})_{\hat\tau}\|\leq e^{-c\cdot d(x',\D V(x,\st(\tau)))}
\end{equation*}
where the constant $c>0$ depends only on $X$.
In particular, it is strongly contracting if $x'$ lies in the interior of $V(x,\st(\tau))$.
\end{cor}
\proof
For $x'$ in the interior of the Weyl cone, and hence $\xi\in\ost(\tau)$, this is a direct consequence. 
For $x'$ on the boundary of the cone, it follows by continuity.
\qed

\medskip
In order to show that $\vartheta_{xx'}$ has expanding directions at $\hat\tau$ if $x'$ lies outside the Weyl cone, 
we consider its action on certain invariant submanifolds of $C(\tau)$
corresponding to parallel sets of singular hyperplanes.

Let again $F$ be a maximal flat with 
$l_{xx'}\subset F\subset  P(\tau,\hat\tau)$.
Moreover, let $h\subset\geo F$ be a half-apartment such that $\interior(\tau)\subset\interior(h)$. 
Then $\vartheta_{xx'}$ fixes $\geo F$ pointwise.
Hence it preserves $h$, the parallel set $P(\D h)$ and the submanifold 
$C(h)=N_h\hat\tau\subset C(\tau)$. 

If $l_{xx'}$ is parallel to the euclidean factor of $P(\D h)$,
equivalently, if $\xi\in\geo l_{xx'}\subset\D h$,
then $\vartheta_{xx'}$ acts trivially on the cross section $CS(\D h)$,
and hence also trivially on $C(h)\cong\geo CS(\D h)-\{\zeta\}$.

In the general case,
the action of $\vartheta_{xx'}$ on $C(h)$ corresponds to the restriction of the action of 
$\vartheta_{xx'}\pi^{\tau}_{x'x}$ to $Hc^h_x=Hc^{\tau}_x\cap P(\D h)$.
When projecting to $CS(\D h)$, 
the latter action in turn corresponds to  the action of 
$\vartheta_{\bar x\bar x'}\pi^{\zeta}_{\bar x'\bar x}$ on $Hc^{\zeta}_{\bar x}$.
Here, $\vartheta_{\bar x\bar x'}$ denotes the transvection on $CS(\D h)$ 
with axis $l_{\bar x\bar x'}$ through $\bar x$ and $\bar x'$ 
mapping $\bar x'\mapsto\bar x$, 
and $\pi^{\zeta}_{\bar x'\bar x}:Hc^{\zeta}_{\bar x}\to Hc^{\zeta}_{\bar x'}$ 
denotes the natural identification of horospheres at $\zeta$.
The projection of $F$ to $CS(\D h)$ is a $\vartheta_{\bar x\bar x'}$-invariant geodesic line 
asymptotic to $\zeta$ and another ideal point $\hat\zeta$.
The fixed point $\hat\tau$ of $\vartheta_{xx'}$ on $C(\tau)$ 
corresponds to the fixed point $\hat\zeta$ of $\vartheta_{\bar x\bar x'}$ on 
$\geo CS(\D h)-\{\zeta\}$.

We prove analogues of Lemmata~\ref{lem:diffposev} and~\ref{lem:evestinst}.
\begin{lem}
If $\xi\in\interior(h)$, 
then the restriction of the differential of $\vartheta_{xx'}$ at $\hat\tau$ 
to the invariant subspace $T_{\hat\tau}C(h)\subset T_{\hat\tau}C(\tau)$
is diagonalizable with strictly positive eigenvalues.
\end{lem}
\proof
Modulo the canonical identification $C(h)\cong\geo CS(\D h)-\{\zeta\}$,
$\vartheta_{xx'}$ restricts to $\vartheta_{\bar x\bar x'}$.
The argument for $\vartheta_{\bar x\bar x'}$ on $\geo CS(\D h)$ is then the same as for Lemma~\ref{lem:diffposev}.
We have that $\zeta=l_{\bar x\bar x'}(-\infty)$ because $\xi\in\interior(h)$.
Therefore, 
the tangent vectors to $C(h)$ at $\hat\zeta$ one-to-one correspond to the orthogonal Jacobi fields along $l_{\bar x\bar x'}$
which decay to zero at $\zeta$.
We conclude as before 
that $(d\vartheta_{\bar x\bar x'})_{\hat\zeta}$ 
is diagonalizable with strictly positive eigenvalues.
\qed

\begin{lem}
\label{lem:evestinstrk1}
If $\xi\in\interior(h)$,
then the eigenvalues $\la$ of $(d\vartheta_{xx'})_{\hat\tau}|_{T_{\hat\tau}C(h)}$ satisfy an estimate
\begin{equation}
\label{ineq:eigvsubsp}
c_1 \leq\frac{-\log\la}{b_{\zeta}(x)-b_{\zeta}(x')}\leq c_2
\end{equation}
with constants $c_1,c_2>0$ depending only on $X$.
\end{lem}
\proof
Since $\xi\in\interior(h)$, we have that 
$b_{\zeta}(x)-b_{\zeta}(x')=d(\bar x,\bar x')>0$.
The assertion follows from the contraction estimate (\ref{ineq:contrestrk1})
and the diagonalizability of the differential.
\qed

\begin{cor}
\label{cor:diffnotcontr}
If $x'\in P(\tau,\hat\tau)-V(x,\st(\tau))$, 
then $(d\vartheta_{xx'})_{\hat\tau}$ is not weakly contracting on $T_{\hat\tau}C(\tau)$.
\end{cor}
\proof
By our assumption, we have that $\xi\not\in\st(\tau)$.
Therefore, the half-apartment $h\subset\geo F$ can be chosen so that its interior contains, besides $\interior(\tau)$, 
also $l_{xx'}(+\infty)$.
(Recall that the convex subcomplex $\st(\tau)\cap\geo F$ is an intersection of half-apartments in $\geo F$.)
Then the estimate (\ref{ineq:eigvsubsp}) applied to $\vartheta_{x'x}=\vartheta_{xx'}^{-1}$ 
yields that $(d\vartheta_{xx'})_{\hat\tau}^{-1}$ has eigenvalues in $(0,1)$.
\qed

\medskip
We also can deduce an upper estimate on the strength of the contraction 
if $x'\in V(x,\st(\tau))$, 
complementing Lemma~\ref{lem:evestinst}:

\begin{lem}
\label{lem:contrestinst}
If $\xi\in\st(\tau)$,
then $(d\vartheta_{xx'})_{\hat\tau}$ has an eigenvalue 
satisfying an estimate
\begin{equation*}
-\log\la\leq C\cdot d(x',\D V(x,\st(\tau)))
\end{equation*}
with a constant $C>0$ depending only on $X$.
\end{lem}
\proof
A nearest point $y'$ to $x'$ on $\D V(x,\st(\tau))$ lies 
on $\D V(x,\st(\tau))\cap F=\D V(x,\st(\tau)\cap\geo F)$.
Hence we can choose the half-apartment $h$ so that $b_{\zeta}(y')=b_{\zeta}(x)$ and 
\begin{equation*}
d(x',\D V(x,\st(\tau))) =
b_{\zeta}(x)-b_{\zeta}(x').
\end{equation*}
Now let $\la$ be an eigenvalue of $(d\vartheta_{xx'})_{\hat\tau}|_{T_{\hat\tau}C(h)}$
and apply the upper estimate (\ref{ineq:eigvsubsp}).
\qed

\medskip
Putting the information 
(Corollaries~\ref{cor:wcontrincone}, \ref{cor:diffnotcontr} and Lemma~\ref{lem:contrestinst})
together, we obtain:
\begin{thm}[Infinitesimal contraction of transvections at infinity]
\label{thm:infcontrtrans}
Let $\tau,\hat\tau\subset\geo X$ be a pair of opposite simplices,
and let $\vartheta$ be a nontrivial transvection which has an axis $c\subset P(\tau,\hat\tau)$
through the point $x=c(0)$.
Then the following hold for the differential $d\vartheta_{\hat\tau}$ of $\vartheta$ on $C(\tau)$
at the fixed point $\hat\tau$:

(i) 
$d\vartheta_{\hat\tau}$ is weakly contracting on $T_{\hat\tau}C(\tau)$,
if and only if $\vartheta^{-1}x\in V(x,\st(\tau))$,
and strongly contracting if and only if $\vartheta^{-1}x\in\interior(V(x,\st(\tau)))$.

(ii)
Suppose that $\vartheta^{-1}x\in V(x,\st(\tau))$.
Then the exponential contraction rate of $d\vartheta_{\hat\tau}$
is comparable to $d(\vartheta^{-1}x,\D V(x,\st(\tau))$,
i.e.\ there is an estimate 
\begin{equation*}
c_1\cdot d(\vartheta^{-1}x,\D V(x,\st(\tau))) \leq
-\log\| d\vartheta_{\hat\tau}\|
\leq c_2\cdot d(\vartheta^{-1}x,\D V(x,\st(\tau)))
\end{equation*}
with constants $c_1,c_2>0$ depending only on $X$.
\end{thm}

We will later use the following consequence of the theorem for general isometries in $G$.
\begin{cor}[Infinitesimal expansion of isometries at infinity]
\label{cor:expand}
Let $\tau\subset\geo X$ be a simplex of type $\tau_{mod}$,
$x\in X$ a point and $g\in G$ an isometry such that 
$d(gx,V(x,\st(\tau))<r$.
Then the exponential expansion rate $\log\eps(g^{-1},\tau)$ of $g^{-1}$ on $\Flagt$ at $\tau$ 
is comparable to $d(gx,\D V(x,\st(\tau))$ in the sense that
\begin{equation*}
C^{-1}\cdot d(gx,\D V(x,\st(\tau))) - A\leq
\log \eps(g^{-1},\tau)
\leq C\cdot d(gx,\D V(x,\st(\tau))) + A
\end{equation*}
with constants $C,A>0$ depending only on $X$, $r$ and the chosen background Riemannian metric on $\Flagt$.
\end{cor}
\proof 
We can write the isometry $g$ as a product $g=tb$ 
of a transvection $t$ along a geodesic $l$ through $x$ asymptotic to $\st(\tau)$, $l(+\infty)\in\st(\tau)$,
and an isometry $b\in G$ which is bounded in terms of the radius $r$.
Then $t$ fixes $\tau$ on $\Flagt$, 
and the expansion factor $\eps(g^{-1},\tau)$ equals 
$\eps(t^{-1},\tau)$
up to a multiplicative constant depending on $X$, $r$ and the background Riemannian metric on $\Flagt$.
Furthermore, $\eps(t^{-1},\tau)=\|dt_{\tau}\|^{-1}$.

Let $\hat\tau$ denote the simplex opposite to $\tau$ with respect to $x$.
Applying Theorem~\ref{thm:infcontrtrans}(ii) to $t=\vartheta^{-1}$ while exchanging the roles of $\tau$ and $\hat\tau$,
we obtain that the exponential contraction rate of $dt_{\tau}$
is comparable to 
$d(t^{-1}x,\D V(x,\st(\hat\tau))=d(tx,\D V(x,\st(\tau))$, i.e.:
\begin{equation*}
c_1\cdot d(tx,\D V(x,\st(\tau))) \leq
-\log\| dt_{\tau}\|
\leq c_2\cdot d(tx,\D V(x,\st(\tau)))
\end{equation*}
Since $d(gx,tx)=d(bx,x)$ is bounded in terms of $r$ and $X$,
the assertion follows.
\qed

\section{Topological dynamics preliminaries} \label{sec:dyn prelim}

In this section we collect various definitions and results from topological dynamics; most of them are rather standard but some are new. 

Throughout this section, we let $\Ga$ be a discrete group, i.e.\ a group equipped with the discrete topology.  
We say that a sequence $(\ga_n)$ of elements of the discrete group $\Ga$ {\em diverges to infinity}, $\ga_n\to\infty$, if the map $\N \to \Ga, n\mapsto \ga_n$ is proper.  
We consider continuous actions $\Ga \acts Z$
on compact metric spaces $(Z,d)$.

\subsection{Expanding actions}
\label{sec:trexpand}

The following notion 
due to Sullivan \cite[\S 9]{Sullivan}
will be of basic importance to us: 

\begin{dfn}[Expanding action]\label{def:expanding}
We say that the action $\Ga\acts Z$ is {\em expanding} 
at the {\em point} $z\in Z$ 
if there exists an element $\ga\in\Ga$ 
which is {\em uniformly expanding} on a neighborhood $U$ of $z$, 
i.e.\ for some constant $c>1$ and all points $z_1,z_2\in U$ we have  
\begin{equation*}
d(\ga z_1,\ga z_2)\geq c\cdot d(z_1,z_2) .
\end{equation*}
We say that the action of $\Ga$ is {\em expanding} 
at a compact $\Ga$-invariant {\em subset} $E\subset Z$ 
if it is expanding at all points $z\in E$. 
\end{dfn}

\subsection{Convergence actions} \label{sec:conv actions}

Let $Z$ be a compact metric space. We define the space $\Trip(Z)$ to be the subset of $Z^3$ consisting of triples of pairwise distinct points in $Z$. Every topological action $\Ga\acts Z$ induces a topological action 
$\Ga \acts\Trip(Z)$.  

\begin{definition}[Convergence action]
The action $\Ga\acts Z$ is called a {\em convergence action} and the image of $\Ga$ in $\Homeo(Z)$ is said to be a {\em convergence group} if one of the following equivalent conditions holds:

(i) The action $\Ga \acts\Trip(Z)$ is properly discontinuous. 

(ii) For every  sequence $\ga_n\to\infty$ in $\Ga$ 
there exist points $z_{\pm}\in Z$ and a subsequence of 
$(\ga_n)$  which converges  to the constant map $\equiv z_+$ 
uniformly on compacts in $Z-\{z_-\}$. The points $z_+$ and $z_-$ are called the {\em limit point} 
and the {\em exceptional point} of this subsequence. 

A convergence action $\Ga\acts Z$ is said to be {\em uniform} 
if the action $\Ga\acts\Trip(Z)$ is cocompact. 
\end{definition}

A proof for the equivalence of both definitions can be found in 
\cite{Bowditch_config}.

The main example of convergence actions comes from the following fact: 
Every discrete group $\Ga$ of isometries of a proper Gromov hyperbolic geodesic metric space $X$ acts as a convergence group on the Gromov boundary $\geo X$ of $X$. 
Furthermore, every word hyperbolic group $\Ga$ acts on its Gromov boundary $\geo \Ga$ as a uniform convergence group. 

Bowditch proved that, vice versa, 
this dynamical behavior characterizes the natural actions of 
word hyperbolic groups 
on their boundaries: 
\begin{thm}[{\cite[Thm.\ 0.1]{Bowditch1998}}]
\label{thm:charhypbow}
Let $\Ga \acts Z$ be a uniform convergence action 
on a non-empty perfect compact metric space. 
Then $\Ga$ is word hyperbolic and $Z$ is equivariantly homeomorphic 
to $\geo\Ga$. 
\end{thm}
The uniformity of a convergence action 
is in turn equivalent to all points being conical. 
The notion of conical limit point 
(e.g.\ of a Kleinian group)
can be expressed purely in terms of the dynamics at infinity 
and, therefore, extends to the more general context considered here:
\begin{dfn}[Intrinsically conical {\cite[\S 8]{Bowditch1998}}]
\label{dfn:intrcon}
Let $\Ga\acts Z$ be a convergence action. 
A point $z\in Z$ is called {\em intrinsically conical} 
if there exists a sequence $\ga_n\to\infty$ in $\Ga$ 
such that the sequence of points $\gamma_n^{-1}z$ converges 
and the sequence of maps $\gamma_n^{-1}|_{Z-\{z\}}$ converges 
(uniformly on compacta)
to a constant map with value $\neq\lim_{n\to\infty}\gamma_n^{-1}z$. 
\end{dfn}

We note that the locally uniform convergence of $\gamma_n^{-1}|_{Z-\{z\}}$ to a constant map 
implies that $z$ is a limit point of $\Ga$. We, thus, will refer to such point $z$ as an {\em intrinsically conical limit point of $\Ga$.} 

\begin{thm}[{\cite[Thm.\ 8.1]{Bowditch1998}}, \cite{Tukia}]
\label{thm:bowditch-conical}
A convergence action $\Ga \acts Z$ on a perfect compact metric space $Z$ 
is uniform  if and only if every point in $Z$ is intrinsically conical. 
\end{thm}
\begin{rem}
The easy direction is that uniformity implies conicality. 
This can be seen as follows: 
Let $z'\neq z$ and $z''_n\to z$ be a sequence of points different from $z$. 
By uniformity, 
there exist elements $\ga_n\to\infty$ in $\Ga$ 
such that we have convergence 
$\ga_n^{-1}z\to z_{\infty}$, 
$\ga_n^{-1}z'\to z'_{\infty}$ and 
$\ga_n^{-1}z''_n\to z''_{\infty}$ 
with pairwise different limits. 
Since $\Ga\acts Z$ is a convergence action, 
we have convergence of $\ga_n^{-1}$ to a constant map 
on $Z-\{z\}$ or on $Z-\{z'\}$. 
The latter is impossible 
because the convergence is locally uniform 
and $\ga_n^{-1}z''_n\to z''_{\infty}\neq z_{\infty}$. 
Thus, the point $z$ is intrinsically conical. 
\end{rem}
The following result connects expanding actions with Bowditch's theorem.
Note that 
if we equip the boundary of a word hyperbolic group $\Ga$ 
with a visual metric $d$, 
then the natural action $\Ga\acts(\geo\Ga,d)$ is expanding, 
see e.g.\ \cite{CP}. 

\begin{lem}\label{lem:conical}
If $\Gamma\acts Z$ is an expanding convergence action 
on a perfect compact metric space $Z$, 
then all points in $Z$ are intrinsically conical. 
\end{lem}
\proof 
We start with a general remark concerning expanding actions. 
For every point $z\in Z$ there exist an element $\ga\in\Ga$
and constants $r>0$ and $c>1$ such that 
$\ga$ is a $c$-expansion on the ball $B(z,r)$ and 
$\ga(B(z,r'))\supset B(\ga z,cr')$ for all radii $r'\leq r$. 
To see this, 
suppose that $c$ is a local expansion factor for $\ga$ at $z$ and, 
by contradiction, 
that there exist sequences of radii $r_n\to0$ 
and points $z_n\not\in B(z,r_n)$ 
such that $\ga z_n\in B(\ga z,cr_n)$. 
Then $z_n\to z$ due to the continuity of $\ga^{-1}$ and, 
for large $n$, 
we obtain a contradiction to the local $c$-expansion of $\ga$. 
Since $Z$ is compact, 
the constants $r$ and $c$ can be chosen uniformly. 
It follows by iterating expanding maps 
that for every point $z$ and every neighborhood $V$ of $z$ 
there exists $\ga\in\Ga$ such that 
$\ga(V)\supset B(\ga z,r)$, 
equivalently, 
$\ga(Z-V)\subset Z-B(\ga z,r)$. 

To verify that a point $z$ is intrinsically conical, 
let $V_n$ be a shrinking sequence of neighborhoods of $z$, 
$$
\bigcap_{n} V_n= \{z\},
$$ 
and let $\ga_n\in\Ga$ be elements such that 
$\ga_n^{-1}(Z-V_n)\subset Z-B(\ga_n^{-1}z,r)$. 
Since $V_n$ is shrinking 
and $\ga_n^{-1}(V_n)\supset B(\ga_n^{-1}z,r)$ 
contains balls of uniform radius $r$, 
it follows that 
the $\ga_n^{-1}$ do not subconverge uniformly 
on any neighborhood of $z$; 
here we use that $Z$ is perfect. 
In particular, $\ga_n\to\infty$. 
The convergence action property implies that, 
after passing to a subsequence, 
the $\ga_n^{-1}$ must converge locally uniformly on $Z-\{z\}$. 
Moreover, we can assume that the sequence of points $\ga_n^{-1}z$ converges. 
By construction, 
its limit will be different 
(by distance $\geq r$) from the limit of the sequence of maps 
$\ga_n^{-1}|_{Z-\{z\}}$. 
Hence the point $z$ is intrinsically conical. 
\qed

\medskip
Suppose that $\Ga\acts Z$ is a convergence action. The set of limit points of sequences in $\Ga$ 
is called the {\em limit set} $\La(\Ga)$ of $\Ga$; the limit set can be also described as the set of 
exceptional points of sequences in $\Ga$. The group $\Ga$ is called {\em elementary} if 
$\La(\Ga)$ contains at most 2 points and {\em nonelementary} otherwise. 

We will  need the following theorem, proven in the case of groups acting 
on spheres by Gehring and Martin  \cite[Theorem 6.3]{GehringMartin} and by Tukia \cite[Theorem 2S]{Tukia1994} in general:  

\begin{thm}\label{thm:minimal action}
If $\Ga$ is nonelementary then the action of $\Ga$ on its limit set $\La(\Ga)$ is minimal and 
$\La(\Ga)$ is perfect. 
\end{thm}

\section{Convex cocompact groups of isometries 
of rank one symmetric spaces}
\label{sec:rank1}

In this section we review equivalent definitions and properties of convex cocompact groups of isometries of negatively curved symmetric spaces. Most of this discussion remains valid in the case of isometry groups of proper $CAT(-1)$ spaces. The main reference for this material is the paper of Bowditch \cite{Bowditch1995}. We also refer the reader to \cite{Kapovich2007} for a survey of discrete isometry groups of rank one symmetric spaces (primarily focused on higher-dimensional real-hyperbolic spaces). 

Let $\Gamma\subset G=\Isom(X)$ be a discrete subgroup of the group of isometries of a negatively curved symmetric space $X$. We let $\Lambda=\Lambda(\Gamma)\subset \geo X$  denote the {\em limit set} of $\Gamma$, i.e.\ the accumulation set of a $\Gamma$-orbit in $X$. Note that $\La$ is necessarily closed in $\geo X$. Then $\Omega=\Omega(\Gamma)= \geo X-\Lambda$ is the {\em domain of discontinuity} of $\Gamma$, which is 
also the wandering set for the action $\Ga\acts \geo X$, and hence is 
the largest open subset of $\geo X$ where $\Gamma$ acts properly discontinuously. 
The {\em Nielsen hull} $N(\La)$ of $\La$ is defined as the smallest closed convex subset in $X$ whose ideal boundary contains $\La$. The set $N(\La)$ exists provided that $\La$ contains at least two points; 
in this case, $\geo N(\La)=\La$. In what follows, we will consider only {\em nonelementary} discrete subgroups $\Ga$, i.e., subgroups for which $\La(\Ga)$ contains more than 2 points. 

The following definition explains the terminology {\em convex cocompact}. 
\begin{dfn}
[C1] A discrete subgroup $\Gamma\subset G$ is called {\em convex cocompact} if $N(\La(\Ga))/\Ga$ is compact. 
\end{dfn}

In particular, such group $\Ga$ is finitely presented and, moreover, word hyperbolic. 

A limit point $\la\in \La$ is called {\em a conical limit point} of $\Gamma$ if for some (every) geodesic ray $\rho$ in $X$ asymptotic to $\la$ there exists a sequence $\gamma_i x\subset X$ converging to $\la$ in an $R$-neighborhood of $\rho$ for some $R<\infty$.   

\begin{dfn}
[C2] A discrete subgroup $\Gamma\subset G$ is called {\em convex cocompact} if every limit point of $\Ga$ is conical. 
\end{dfn} 

In fact, one can get $R$ to be uniform for $x\in N(\La)$ and $\la\in \La$.

Recall (see section \ref{sec:conv actions}) that for a set $Z$, $\Trip(Z)$ denotes the set of triples of pairwise distinct points in $Z$. 

\begin{dfn}
[C3] A discrete subgroup $\Gamma\subset G$ is called {\em convex cocompact} if the action $\Gamma\acts \Trip(\La)$ is cocompact.  
\end{dfn} 
Every discrete group $\Gamma\subset G$ acts properly discontinuously on $X\cup \Om$, which we equip with the subset topology induced from $\bar{X}=X\cup \geo X$. 
\begin{dfn}
[C4] A discrete subgroup $\Gamma\subset G$ is called {\em convex cocompact} if the action $\Gamma\acts 
X\cup \Om$ is cocompact. 
\end{dfn} 

Definition C4 implies that for every convex cocompact subgroup $\Gamma\subset G$, the quotient $\Om/\Ga$ is compact. 
The converse is false, as the following examples show. 
\begin{example}
1. Consider a cyclic group $\Gamma$ of parabolic isometries of the hyperbolic plane $\H^2=X$. 
Then $\La(\Ga)$ is a singleton, $\Om/\Ga$ is homeomorphic to $S^1$, while  
$$
(X\cup \Omega)/\Ga \cong [0,\infty)\times S^1
$$
is noncompact. Thus, $\Gamma$ is not convex cocompact. 
In this case, of course, $\Gamma$ contains unipotent (parabolic) elements. The next three examples contain only loxodromic elements. 

2. Let $S$ denote a closed hyperbolic surface, $\pi:=\pi_1(S)$. Then $\pi$ admits a discrete and faithful representation $\rho: \pi\to G=\Isom(\H^3)$, so that its image $\Ga=\rho(\pi)$ is a 
{\em totally-degenerate purely loxodromic subgroup} of $G$: $\Om(\Ga)$ is simply connected and nonempty, $\Gamma$ contains no parabolic elements and
$$
(X\cup \Om)/\Ga
$$
is homeomorphic to $S\times [0,\infty)$, where $S\times \{0\}$ corresponds to $\Om/\Gamma$, 
see \cite{Bers}. 
Thus, $\Gamma$ is not convex cocompact. 

3. Let $M$ be a closed oriented hyperbolic $m$-manifold with a nonseparating oriented closed totally-geodesic hypersurface $N$. Such manifolds exist for all $m$ (see \cite{Millson(1976)}). 
Let $\tilde{M}\to M$ denote the infinite cyclic cover determined by the homomorphism $\pi_1(M)\to \Z$ corresponding to the element of $H^1(M,\Z)$ Poincar\'e dual to the homology class $[N]$. Then $N$ lifts to a submanifold $N_0\subset \tilde{M}$ which is isometric to $N$ and which separates $\tilde{M}$ in two components $\tilde{M}_-, \tilde{M}_+$. Let $M'$ denote the metric completion of the Riemannian manifold $\tilde{M}_+$. Then $M'$ is a complete hyperbolic manifold with 
single geodesic boundary component isometric to $N_0$ and injectivity radius bounded below. The fundamental group $\Gamma$ of $M'$ is not finitely generated. The hyperbolic structure on $M'$ determines a discrete isometric action $\Gamma \acts \H^m$, so that $\Gamma$ contains no parabolic elements. Then 
$$
(\H^m\cup\Om)/\Gamma 
$$
is homeomorphic to $M'$. In particular, $\Om/\Gamma$ is compact and nonempty, while $\Gamma$ is not even 
finitely generated; in particular, $\Gamma$ is not convex cocompact.  

4. Similarly, for every simple Lie group $G$ of rank 1, there are discrete subgroups $\Ga\subset G$ whose limit set is the entire sphere $\geo X$, but $\Ga$ is not  finitely generated. For instance, one can 
start with a uniform lattice $\widehat{\Ga}\subset G$; being a non-elementary word hyperbolic group, $\widehat{\Ga}$ admits a normal subgroup $\Ga$ which is 
isomorphic to the free groups of countably infinite rank. The limit set of $\Ga$ is necessarily the entire sphere $\geo X$. 
Furthermore, when $X$ is a real-hyperbolic 3-space or complex-hyperbolic plane or complex-hyperbolic 3-space, there are examples of finitely generated subgroups $\Ga\subset G$ whose limit set is $\geo X$, but $\Ga$ is not a lattice in $G$. In the case $X=\H^3$, such examples can be constructed, for instance, using normal surface subgroups in 
fundamental groups of closed hyperbolic  3-manifolds fibering over the circle. For examples in $\C\H^2$ and $\C\H^3$ see e.g. \cite{Kapovich2013}: These are normal subgroups 
in complex-hyperbolic manifolds which admit (singular)  holomorphic fibrations over hyperbolic Riemann surfaces. 
\end{example}

On the other hand, the phenomenon in the Examples 2 and 3 can essentially only occur in the real-hyperbolic case: 

\begin{thm}
[See \cite{Ramachandran}] Let $X$ be a negatively curved rank one symmetric space which is not real-hyperbolic. Suppose that $\Gamma\subset G=\Isom(X)$ is a discrete torsion-free subgroup without unipotent elements so that $\Om(\Ga)/\Ga$ is compact and nonempty. Then $\Gamma$ is convex cocompact provided that  $X$ is not isometric to $\C\H^2$. In the case $X=\C \H^2$, the same result holds provided that the Riemannian manifold $X/\Ga$  has injectivity radius bounded below. 
\end{thm}
Let $\Ga\subset G$ be a discrete subgroup. 
Pick a point $x\in X$ which is not fixed by any nontrivial element of $\Ga$, 
and define the {\em Dirichlet fundamental domain} $D_x$ of $\Ga$ as
$$
D_x=\{y\in X: d(x,y)\le d(\ga x, y), \forall \ga\in \Ga\}. 
$$
Note that $D_x$ is convex if $X$ is real-hyperbolic, but is not convex otherwise. In general, $D_x$ is starlike with the center $x$; since $X$ is Gromov hyperbolic, this implies that $D_x$ is {\em quasiconvex} in $X$. Subsets of $D_x$ of the form
$$
D_x\cap \ga D_x, \quad\ga\in \Ga,
$$
are called {\em faces} of $D_x$. 

Let $\bar{D}_x$ denote the union
$$
D_x\cup (\geo D_x \cap \Om(\Ga)),
$$
which is a certain partial compactification of $D_x$. It follows (almost) immediately from C4 that $\Gamma$ is convex cocompact if and only if $\bar{D}_x$  is compact. The following definition is a more elaborate version of this observation:

\begin{dfn}
[C5] A discrete subgroup $\Gamma\subset G$ containing no parabolic elements is called {\em convex cocompact} if one (every) Dirichlet fundamental domain $D_x$ has finitely many faces. 
\end{dfn}

Note that a cyclic unipotent subgroup of $\Isom(\H^4)$ can have  a
Dirichlet domain with infinitely many faces.

\begin{dfn}
[C6] A discrete subgroup $\Gamma\subset G$ is convex cocompact whenever $\Gamma$ is word hyperbolic as an abstract group and there exists an equivariant homeomorphism 
$$
\beta: \geo \Gamma \to \La(\Ga),
$$
where $\geo \Ga$ is the Gromov boundary of $\Gamma$. 
\end{dfn}

Note that the injectivity of $\beta$ is critical here: 
\begin{thm}
[Mj, \cite{Mj1, Mj2}] 
Suppose that $\Ga\subset\Isom(\H^3)$ is a word hyperbolic subgroup 
(not necessarily convex cocompact). Then 
there always exists a continuous equivariant map $\beta: \geo \Gamma\to \La(\Ga)$; the map 
$\beta$ is called the {\em Cannon-Thurston map} 
for the subgroup $\Ga\subset G$.  
\end{thm}

Let $\Ga\subset G$ be a finitely generated subgroup of a Lie group $G$ with finitely many connected components; 
we will equip $\Ga$ with a word metric. 
A point $x\in X$ defines the  {\em orbit map} 
$\Ga \to \Ga x\subset X$. 
The subgroup $\Ga$ is called {\em undistorted} in $G$, 
if some (any) orbit map $\Ga \to X$ is a quasi-isometric embedding, 
equivalently, 
if the inclusion $\Ga\to G$ is a quasi-isometric embedding, 
where $G$ is equipped with a left invariant Riemannian metric. 
\begin{dfn}
[C7] 
A discrete subgroup $\Gamma\subset G$ is convex cocompact  
if it is undistorted. 
\end{dfn}
Note that, in view of the hyperbolicity of $X$, 
undistortion of $\Ga$ implies that 
the quasi-isometrically embedded $\Ga$-orbits 
are quasi-convex subsets of $X$. In particular:

1. $\Ga$ is word hyperbolic, and hence the orbit maps $\Ga\to\Ga x$ 
continuously extend at infinity to 
an equivariant homeomorphism $\beta: \geo \Ga \to \La(\Ga)$. 

2. The $\Ga$-equivariant relation in $X\times \Ga$ 
given by the nearest-point projection to an orbit $\Ga x$ 
is a coarse Lipschitz retraction $X\to\Ga x$. 

The converse to this is also easy: 
\begin{dfn}
[C8] A discrete subgroup $\Gamma\subset G$ is convex cocompact if for 
some (every) $\Ga$-orbit $\Ga x\subset X$ 
there exists a $\Ga$-equivariant coarse Lipschitz retraction $X\to \Ga x$. 
\end{dfn}

The equivariance condition for the retraction can be omitted: 
\begin{dfn}
[C9] A discrete subgroup $\Gamma\subset G$ is convex cocompact  if for 
some (every) $\Ga$-orbit $ \Ga x\subset X$ 
there exists a coarse Lipschitz retraction $X\to  \Ga x$. 
\end{dfn}
Our last characterization of convex cocompactness is 
in terms of expanding actions. 
We fix a visual metric $d$ on $S=\geo X$. 
\begin{dfn}
[C10] A discrete subgroup $\Gamma\subset G$ is convex cocompact if 
its action on $\geo X$ is expanding at every point 
of $\La(\Ga)$, see Definition \ref{def:expanding}. 
\end{dfn}

This interpretation of convex cocompactness appears 
in Sullivan's paper \cite{Sullivan}. 
\begin{thm}\label{thm:rk1}
The definitions C1--C10 are equivalent. 
\end{thm}
The equivalence of Definitions C1--C9 can be found for instance in \cite{Bowditch1995}. 
The implication C10 $\Rightarrow$ C2  
is a corollary of Lemma \ref{lem:conical}. 
In the case of real-hyperbolic space, the implication 
C5 $\Rightarrow$ C10 is immediate by taking a Ford fundamental domain (bounded by finitely many isometric spheres  $I(\gamma_i), I(\gamma_i^{-1})$, $i=1,\ldots,k$) and observing that $\gamma_i$ is a strict expansion on every compact contained in the open round ball bounded by $I(\gamma_i)$. For the remaining rank one symmetric spaces the implication C2 $\Rightarrow$ C10 is a corollary of our Proposition \ref{prop:conicimplexpatau} in section \ref{sec: Expansion at conical limit flags}.

\section{Weakly regular subgroups and their limit sets}\label{sec:regularity and limit sets}
\subsection{Weak regularity}

In this section we introduce and discuss an asymptotic regularity property 
for discrete subgroups $\Ga$ of a semisimple group $G$. 
It is a condition on the asymptotics of divergent sequences in the subgroup 
and is defined with respect to a $\iota$-invariant face type 
$\tau_{mod}\subset\si_{mod}$.

We first define a stronger {\em uniform} version of regularity 
which can be stated more directly in terms of the limit set $\La(\Ga)=\ol{\Ga x}\cap\geo X$ in the visual boundary.
 
We recall that the subset of {\em $\tau_{mod}$-regular types} in $\si_{mod}$
is the {\em open star} $\ost(\tau_{mod})$, that is,
the union of all open faces of $\si_{mod}$ which contain $\tau_{mod}$ in their closure,
see section~\ref{sec:star}. 
The {\em $\tau_{mod}$-regular part} of the ideal boundary is then defined as the subset
\begin{equation*}
\geo^{\tau_{mod}-reg} X=\theta^{-1}(\ost(\tau_{mod}))\subset\geo X 
\end{equation*}
of all $\tau_{mod}$-regular ideal points. 
It consists of the ideal points for which there is a unique closest (with respect to the Tits metric) 
simplex $\tau\subset\geo X$ of type $\tau_{mod}$. 
It contains all open chambers and is in particular 
dense in $\geo X$ (also in the Tits topology). 

\begin{dfn}[Uniformly weakly regular]
We call the subgroup $\Ga$ {\em uniformly $\tau_{mod}$-regular} 
if its visual limit set consists only of $\tau_{mod}$-regular ideal points, 
$\La(\Ga)\subset\geo^{\tau_{mod}-reg}X$. 
\end{dfn}
Note that $\La(\Ga)$ is compact, as is its type projection to $\si_{mod}$. 
A quantitative version of uniform regularity is given by:
\begin{dfn}[$\Theta$-regular]
Let $\Theta\subset\ost(\tau_{mod})$ be compact. 
The subgroup $\Ga$ is called {\em $\Theta$-regular} 
if its visual limit set consists only of ideal points of type $\Theta$, 
$\La(\Ga)\subset\theta^{-1}(\Theta)\subset\geo^{\tau_{mod}-reg}X$. 
\end{dfn}
These notions apply in the same way to divergent sequences $(x_n)$ in $X$ 
and $(g_n)$ in $G$, 
the latter by looking at the accumulation set in $\geo X$ of an associated orbit sequence $(g_nx)$ for a(ny) base point $x$.

Now we define (non-uniform) weak regularity itself.
We denote by $\D\ost(\tau_{mod})=\si_{mod}-\ost(\tau_{mod})$ 
the set of {\em $\tau_{mod}$-singular} types. 

\begin{dfn}
[Weakly regular sequence] 
\label{def:tauregseq}
(i) A sequence $\de_n\to\infty$ in $\De=V_{mod}$ 
is {\em $\tau_{mod}$-regular} if 
\begin{equation}
\label{eq:deldistnonsing}
d(\de_n,V(0,\si_{mod}-\ost(\tau_{mod})))\to+\infty.
\end{equation}

(ii) A sequence $x_n\to\infty$ in $X$ 
is {\em $\tau_{mod}$-regular}
if for some (any) base point $x$ the sequence
of $\De$-valued lengths $d_{\De}(x,x_n)$ in $\De$ 
has this property. 

(iii)
A sequence $g_n\to\infty$ in $G$ 
is {\em $\tau_{mod}$-regular} 
if some (any) orbit $(g_nx)$ in $X$ has this property. 
\end{dfn}
\begin{rem}
(i) 
The independence of the base point and the orbit 
in parts (ii) and (iii) of the definition 
is due to the triangle inequality
$|d_{\De}(x,y)-d_{\De}(x',y')|\leq d(x,x')+d(y,y')$.

(ii)
{\em Uniform} $\tau_{mod}$-regularity is equivalent to {\em linear} (with respect to $|\delta_n|$) divergence in (\ref{eq:deldistnonsing}). 

(iii)
If the sequence $(x_n)$ is $\tau_{mod}$-regular, then every sequence $(x'_n)$ uniformly close to it, 
$\sup_n d(x_n,x'_n)<+\infty$, is also $\tau_{mod}$-regular, again by the triangle inequality.
Similarly, if the sequence $(g_n)$ in $G$ is $\tau_{mod}$-regular, 
then for all bounded sequences $(b_n)$ and $(b'_n)$ in $G$, 
the sequence $(b_ng_nb'_n)$ is also $\tau_{mod}$-regular. 

(v)
(Uniform) $\tau_{mod}$-regularity implies 
(uniform) $\tau'_{mod}$-regularity for all face types $\tau'_{mod}\subset\tau_{mod}$, 
because $\ost(\tau_{mod}) \subset \ost(\tau'_{mod})$.

(vi)
Every diverging sequence has a weakly regular subsequence,
i.e.\ a subsequence which is $\tau_{mod}$-regular for some face type $\tau_{mod}$. 
\end{rem}

\begin{dfn}
[Weakly regular subgroup]
\label{def:tauregsubs}
A subgroup $\Ga\subset G$ 
is {\em $\tau_{mod}$-regular} 
if all sequences $\ga_n\to\infty$ in $\Ga$ have this property. 
When we do not want to specify $\tau_{mod}$, we refer to $\Ga$ simply as {\em weakly regular}.
\end{dfn}

\begin{rem}
(i)
The definition of $\tau_{mod}$-regularity for {\em sequences} makes sense also if $\tau_{mod}$ is not $\iota$-invariant. 
Then a sequence $(g_n)$ in $G$ 
is (uniformly) $\tau_{mod}$-regular 
if and only if the sequence $(g_n^{-1})$ of inverses
is (uniformly) $\iota\tau_{mod}$-regular,
cf.\ the symmetry property (\ref{eq:symprop}) for $\De$-lengths. 
When defining $\tau_{mod}$-regularity for {\em subgroups},
it is therefore natural to require the $\iota$-invariance of the simplex $\tau_{mod}$,
and this is why we impose this condition in the entire paper.

(ii)
A discrete subgroup $\Ga\subset G$ needs not be $\tau_{mod}$-regular (for any $\tau_{mod}$)  
even if all its nontrivial elements 
are regular (transvections). 
This can happen e.g.\ for free abelian subgroups of rank $\geq 2$. 
\end{rem}

\begin{rem}
[Relation to visual compactifications for Finsler metrics]
\label{rem:regfinscomp}
We recall that a sequence $x_n\to\infty$ converges to an ideal point in the visual compactification 
if and only if the normalized distance functions $d(\cdot,x_n)-d(p,x_n)$ converge (locally uniformly), 
where $p$ is some base point. 
Ideal boundary points can thus be identified with normalized Busemann functions. 
The same construction can be carried out for $G$-invariant {\em Finsler} metrics on $X$ 
and one obtains modified visual compactifications. 
To a face type $\tau_{mod}$ we can associate a Finsler metric as follows.
Fix a $\iota$-invariant unit vector $v\in\De$ pointing to an interior point of $\tau_{mod}$ 
(e.g.\ its center), 
and define $d_v(x,x'):=\langle d_{\De}(x,x'),v\rangle$. 
In the visual boundary with respect to this Finsler metric 
there is a unique $G$-orbit which is a copy of $\Flag(\tau_{mod})$;
its points correspond to collapsed open stars around simplices of type $\tau_{mod}$ in $\geo X$. 
(The union $X\cup\Flag(\tau_{mod})$ is the natural {\em $\tau_{mod}$-bordification} of $X$.)
The relation with $\tau_{mod}$-regularity is as follows: 
A subgroup $\Ga$ is $\tau_{mod}$-regular if and only if its visual limit set with respect to the modified compactification of $X$
is contained in the $G$-orbit $\Flag(\tau_{mod})$.
\end{rem}

\subsection{Contraction-expansion dynamics at infinity}
\label{sec:contrexp}

In this section, 
we will describe 
how the weak regularity of sequences in $G$ can be read off their dynamics at infinity. 
Roughly speaking, 
on appropriate flag manifolds, 
certain almost full subsets are contracted, asymptotically, almost to points.
Dually, small balls are expanded to almost full subsets.
This ``contraction-expansion dynamics" is a generalization of convergence dynamics in rank one to arbitrary rank.

\subsubsection{Contraction and expansion}

We now formulate the contraction and expansion properties for sequences in $G$ 
and show that they are satisfied by weakly regular sequences. 
We first state and discuss contraction.
\begin{dfn}[Contraction]
Let $(g_n)$ be a sequence in $G$.
We call a sequence of subsets $U_n\subset\Flag(\tau_{mod})$ 
a {\em contraction sequence} for $(g_n)$ if it satisfies the following properties:

(i) 
It is {\em exhaustive} in the sense that 
there exist a bounded sequence $(b_n)$ in $G$ and a simplex $\tau$ of type $\tau_{mod}$ 
such that $b_nU_n\to C(\tau)$, 
meaning that every compact subset of the open Schubert stratum $C(\tau)$ 
is contained in $b_nU_n$ for $n$ sufficiently large.

(ii) 
The image subsets $g_nU_n$ {\em shrink} in the sense that 
there exists a bounded sequence $(b_n)$ in $G$ 
such that $b_ng_nU_n$ Hausdorff converges to a point, 
equivalently, 
if $\diam(g_nU_n)\to0$ with respect to a Riemannian background metric on $\Flag(\tau_{mod})$. 

We call the sequence $(g_n)$ {\em contracting} on $\Flag(\tau_{mod})$ if it admits a contraction sequence.
\end{dfn}

Note that open Schubert strata are dense open subsets of full volume 
with respect to any auxiliary smooth probability measure on $\Flagt$.
Hence the subsets in a (measurable) exhausting sequence 
have asymptotically full volume. 

Property (i) in the definition means that the subsets $U_n$ asymptotically fill out the 
sequence of moving open Schubert strata $b_n^{-1}C(\tau)=C(b_n^{-1}\tau)$,
and property (ii) means that the images $g_nU_n$ asymptotically concentrate to points
which are also allowed to move.

Note that if the sequence $(g_n)$ in $G$ is contracting,
then for all bounded sequences $(b_n)$ and $(b'_n)$ in $G$, 
the sequence $(b_ng_nb'_n)$ is also contracting.

We will say that a contraction sequence $(U_n)$ is {\em opposite} to a sequence $(\tau_n)$ in $\Flagt$
if $U_n\subset C(\tau_n)$ for large $n$.
Note that in the definition the sequence $(U_n\cap C(b_n^{-1}\tau))$ 
is still exhaustive, and hence a contraction sequence for $(g_n)$ opposite to the sequence $(b_n^{-1}\tau)$.

The contraction property can be reformulated in terms of concentration of measures,
which shows that it essentially agrees with other notions of contraction (on flag manifolds) used in the literature,
cf.\ \cite[\S 3]{GoldsheidMargulis}:
If a sequence $(g_n)$ is contracting in the above sense, 
then there exists a bounded sequence $(b_n)$ such that the sequence of measures $(b_ng_n)_*\mu$
converges to a Dirac mass for any smooth probability measure $\mu$ on $\Flagt$. 
Moreover,
one can also prove the converse. 
Note that we do not require the measures ${g_n}_*\mu$ themselves to converge, but allow them to concentrate around any sequence 
of moving points. 
Since flag manifolds are compact, one can always achieve by passing to a subsequence that 
also the measures ${g_n}_*\mu$ converge. 

\begin{rem}[Proximal]
\label{rem:prox}
Related to the notion of contraction is the notion of {\em proximality}, see e.g.\ \cite{Abels}.
Say that an element $g\in G$ acts as a {\em proximal transformation} on $\Flag(\tau_{mod})$ 
if  the sequence $(g^n)_{n\in\N}$ converges to a point $\tau_+\in\Flagt$ 
locally uniformly on an open Schubert stratum $C(\tau_-)\subset\Flagt$, 
where $\tau_-$ is a simplex antipodal to $\tau_+$.  
We will refer to $\tau_+$ as the {\em attractive} fixed point of $g$. 
(Note that it is necessarily fixed by $g$).
The sequence $(g^n)$ is then contracting on $\Flagt$. 
Note however that vice versa sequences, which are contracting on $\Flagt$, need not contain proximal elements.
For instance, in the rank one case all divergent sequences in $G$ are contracting
due to convergence dynamics.
\end{rem}

The expansion property is dual to contraction:
\begin{dfn}[Expansion]
Let $(g_n)$ be a sequence in $G$.
We call a sequence of subsets $V_n\subset\Flag(\tau_{mod})$ 
an {\em expansion sequence} for $(g_n)$ if it satisfies the following properties:

(i) 
It shrinks, $\diam(V_n)\to0$.

(ii) 
The sequence of image subsets $g_nV_n$ is exhaustive.

We call the sequence $(g_n)$ {\em expanding} on $\Flag(\tau_{mod})$ if it admits an expansion sequence.
\end{dfn}
The duality means that
$(U_n)$ is a contraction sequence for $(g_n)$ if and only if $(g_nU_n)$ is an expansion sequence for $(g^{-1}_n)$. 

\subsubsection{Strong asymptoticity of Weyl cones}

Before proceeding, we need to prove a fact from the geometry of symmetric spaces.

Let $\tau\subset\geo X$ be a simplex of type $\tau_{mod}$. 
For a point $x\in X$ 
we have on the open Schubert stratum $C(\tau)\subset\Flag(\tau_{mod})$ 
the function
\begin{equation*}
\tau'\mapsto d(x,P(\tau,\tau')).
\end{equation*}
It is continuous and proper.
(This follows from the fact that $C(\tau)$ and $X$ 
are homogeneous spaces for the parabolic subgroup $P_{\tau}$.)
It has a unique minimum zero in the chamber $\hat\tau_x$ 
which is opposite to $\tau$ with respect to $x$. 

We define the following open subsets of $C(\tau)$ 
which can be regarded as ``shadows'' of balls with respect to $\tau$. 
For $x\in X$ and $r>0$, we put
\begin{equation*}
U_{\tau,x,r}:=\{\tau'\in C(\tau) | d(x,P(\tau,\tau'))<r\} .
\end{equation*}
Then the subsets $U_{\tau,x,r}$ for fixed $\tau$ and $x$ 
form a neighborhood basis of $\hat\tau_x$. 

The next fact expresses the uniform strong asymptoticity
of asymptotic Weyl cones. 
\begin{lem}
\label{lem:expconvsect}
For $r,R>0$ exists $d=d(r,R)>0$ such that:

If $y\in V(x,\st(\tau))$ with $d(y,\D V(x,\st(\tau)))\geq d(R,r)$, 
then 
$U_{\tau,x,R}\subset U_{\tau,y,r}$.
\end{lem}
\proof
If $U_{\tau,x,R}\not\subset U_{\tau,y,r}$ 
then there exists $x'\in B(x,r)$ such that $d(y,V(x',\st(\tau)))\geq r$. 
Thus, if the assertion is wrong, 
there exist a sequence $x_n\to x_{\infty}$ in $B(x,r)$ 
and a sequence $y_n\to\infty$ in $V(x,\st(\tau))$ 
such that $d(y_n,\D V(x,\st(\tau)))\to+\infty$
and $d(y_n,V(x_n,\st(\tau)))\geq r$.  

Let $\rho:[0,+\infty)\to V(x,\tau))$ be a geodesic ray with initial point $x$ and asymptotic to an interior point of $\tau$.
Then the sequence $(y_n)$ eventually enters every Weyl sector $V(\rho(t),\st(\tau))$. 
Since the distance function $d(\cdot,V(x_n,\st(\tau)))$ is convex and bounded, and hence non-increasing 
along rays asymptotic to $\st(\tau)$, 
we have that 
\begin{equation*}
R\geq d(x,V(x_n,\st(\tau))) 
\geq d(\rho(t),V(x_n,\st(\tau)))\geq d(y_n,V(x_n,\st(\tau)))\geq r
\end{equation*}
for $n$ large. 
It follows that 
\begin{equation*}
R\geq d(\rho(t),V(x_{\infty},\st(\tau)))\geq r
\end{equation*}
for all $t$. 
However, if $\rho$ is asymptotic to $V(x_{\infty},\st(\tau))$, then it is strongly asymptotic, a contradiction. 
\qed

\subsubsection{Regularity implies contraction-expansion}

We are now ready to show 
that weakly regular sequences in $G$ are contracting-expanding on suitable flag manifolds.

The following flexibile and base point independent notion of radial projection to infinity will be useful
for describing contraction and expansion sequences.
\begin{dfn}[Shadows at infinity]
A {\em shadow sequence} of a sequence $(x_n)$ in $X$ 
is a sequence  $(\tau_n)$ of simplices in $\Flag(\tau_{mod})$ such that 
\begin{equation*}
\sup_n d(x_n,V(x,\st(\tau_n)))<+\infty 
\end{equation*}
for some (any) base point $x$. 
A shadow sequence of a sequence $(g_n)$ in $G$ is a sequence of shadows 
for an orbit sequence $(g_nx)$ for some (any) base point $x$. 
\end{dfn}

\begin{prop}
\label{prop:regcompr}
$\tau_{mod}$-regular sequences in $G$ are contracting and expanding on $\Flag(\tau_{mod})$.
More precisely,
let $(g_n)$ be a $\tau_{mod}$-regular sequence in $G$.
Then for a shadow sequence $(\tau_n^-)$ of $(g_n^{-1})$ there exists 

(i) a contraction sequence $(U_n)$ for $(g_n)$ opposite to $(\tau_n^-)$,
i.e.\ $U_n\subset C(\tau_n^-)$.

(ii) an expansion sequence $(V_n)$ for $(g_n)$ containing $(\tau_n^-)$, i.e.\ $\tau_n^-\in V_n$. 
\end{prop} 
\proof
(i)
We fix a base point $x$ 
and denote by $g_n^{-1}x_n$ the nearest point projection of $g_n^{-1}x$ to $V(x,\st(\tau_n^-))$. 
Then the sequence $(x_n)$ is bounded. 
Due to $\tau_{mod}$-regularity, we have that 
\begin{equation*}
d(g_n^{-1}x_n,\D V(x,\st(\tau_n^-)))\to+\infty .
\end{equation*}
Lemma~\ref{lem:expconvsect} yields that 
for any $r,R>0$ the inclusion 
\begin{equation*}
U_{\tau_n^-,x,R} \subset U_{\tau_n^-,g_n^{-1}x_n,r}
\end{equation*}
holds for large $n$. 
Therefore there exist sequences of positive numbers $R_n\to+\infty$ and $r_n\to0$ such that 
\begin{equation*}
g_nU_{\tau_n^-,x,R_n} \subset U_{g_n\tau_n^-,x_n,r_n}
\end{equation*}
for large $n$. 
The sequence $(U_{\tau_n^-,x,R_n})$ is exhaustive, 
and the sequence $(U_{g_n\tau_n^-,x_n,r_n})$ shrinks because $(x_n)$ is bounded.
Hence $(U_{\tau_n^-,x,R_n})$ is a contraction sequence for $(g_n)$,
and $U_{\tau_n^-,x,R_n}\subset C(\tau_n^-)$.

(ii)
Let $g_n^{-1}\tau_n^+$ be the simplices opposite to $\tau_n^-$ with respect to $x$.
Similarly,
Lemma~\ref{lem:expconvsect} implies that 
\begin{equation*}
U_{g_n^{-1}\tau_n^+,g_n^{-1}x_n,R_n} \subset U_{g_n^{-1}\tau_n^+,x,r_n} \ni \tau_n^-
\end{equation*}
for suitable sequences $R_n\to+\infty$ and $r_n\to0$.
The sequence $(U_{\tau_n^+,x_n,R_n})$ is exhaustive because $(x_n)$ is bounded, 
and it follows that $(U_{g_n^{-1}\tau_n^+,x,r_n})$ 
is an expansion sequence for $(g_n)$. 
\qed

\medskip
The proposition 
has several useful consequences. 

Firstly, the various shadow sequences of a weakly regular sequence approach each other. 
We will use this below to define flag convergence. 
\begin{lem}[Asymptotic uniqueness of shadows]
\label{lem:shadasuniq}
For a $\tau_{mod}$-regular sequence $(g_n)$ in $G$ 
any two shadow sequences $(\tau_n)$ and $(\tau'_n)$ asymptotically coincide,
i.e.\ $d(\tau_n,\tau'_n)\to0$.
\end{lem}
\proof
By Proposition~\ref{prop:regcompr} (Part (ii)), 
there exist expansion sequences $(V_n)$ and $(V'_n)$ for $(g_n^{-1})$ such that 
$\tau_n\in V_n$ and $\tau'_n\in V'_n$.
The assertion follows from the fact that any two expansion sequences asymptotically coincide,
i.e.\ $\diam(V_n\cup V'_n)\to0$. 

To see the latter, 
note that $(g_n^{-1}V_n)$ and $(g_n^{-1}V'_n)$ are contraction sequences for $(g_n)$.
In particular, they are exhaustive 
and therefore asymptotically intersect by volume reasons,
that is $g_n^{-1}V_n\cap g_n^{-1}V'_n\neq\emptyset$ for large $n$.
So, $V_n\cap V'_n\neq\emptyset$ for large $n$, 
and hence $\diam(V_n\cup V'_n)\to0$, 
as claimed.
\qed

\medskip
Secondly, 
shadow sequences asymptotically agree with the values on contraction subsets:
\begin{lem}
\label{lem:shadagrcontr}
Let $(g_n)$ be a $\tau_{mod}$-regular sequence in $G$. 
Then for any contraction sequence $(U_n)$ and any shadow sequence $(\tau_n)$, 
the subsets $g_nU_n\cup\{\tau_n\}$ shrink, $\diam(g_nU_n\cup\{\tau_n\})\to0$.
\end{lem}
\proof
Note that $(g_nU_n)$ is an expansion sequence for $(g_n^{-1})$.
According to Proposition~\ref{prop:regcompr} (Part (ii)) applied to $(g_n^{-1})$
there exists another expansion sequence $(V_n^-)$ for $(g_n^{-1})$
such that $\tau_n\in V_n^-$. 
The assertion follows from the asymptotic uniqueness of expansion sequences, 
compare the proof of Lemma~\ref{lem:shadasuniq}, 
which yields that $\diam(g_nU_n\cup V'_n)\to0$. 
\qed

\begin{rem}
The last result relates 
the asymptotics of orbits in $X$
and the dynamics at infinity. 
One can promote it to showing that 
$\tau_{mod}$-regular sequences in $G$ 
have contraction-expansion dynamics on the bordification $X\cup\Flag(\tau_{mod})$.
\end{rem}

\subsubsection{Contraction implies regularity}

We consider now sequences $(g_n)$ in $G$ which are contracting on $\Flagt$, 
and show that they are $\tau_{mod}$-regular.

The key step is a converse to Proposition~\ref{prop:regcompr} (i),
essentially saying that 
sequences opposite to contraction sequences for $(g_n)$ asymptotically coincide with shadow sequences for $(g_n^{-1})$:
\begin{lem}
\label{lem:shadalmopp}
Let $(U_n)$ be a contraction sequence for $(g_n)$, 
and let $(\tau_n^-)$ be a shadow sequence of $(g_n^{-1})$.
Suppose that $(\tau_n)$ is a sequence in $\Flagt$ such that $(U_n\cap C(\tau_n))$ is still exhaustive.
Then  $d(\tau_n,\tau_n^-)\to0$.
\end{lem}
\proof
By passing to a subsequence, 
we may assume convergence $\tau_n^-\to\tau_-$ and $\tau_n\to\tau$.
It suffices to show that then $\tau_-=\tau$. 

We look at the dynamics of $(g_n)$ on the space of maximal flats. 
Recall that a sequence of maximal flats $F_n\subset X$ is bounded
if $d(x,F_n)$ is bounded for a base point $x\in X$. 

Suppose that $(F_n)$ is a bounded sequence of maximal flats, 
such that the sequence $(g_nF_n)$ of image flats is also bounded. 
We will see that its position relative to $\tau$ is restricted.
Consider a maximal flat $F$ which is the limit of a subsequence of $(F_n)$.
\begin{sublem}
The apartment $\geo F\subset\geo X$ contains exactly one simplex opposite to $\tau$.
\end{sublem}
\proof
In a spherical building, every point has an antipode in every apartment.
Hence, $\geo F$ contains at least one simplex $\hat\tau$ opposite to $\tau$. 
Suppose that it contains another simplex $\hat\tau'$ opposite to $\tau$. 
After passing to a subsequence, we may assume that $F_n\to F$.
Let $\hat\tau_n,\hat\tau'_n\subset\geo F_n$ be approximating simplices, 
$\hat\tau_n\to\hat\tau$ and $\hat\tau'_n\to\hat\tau'$. 
Since $(g_nF_n)$ is bounded, after passing to a subsequence, the sequences 
$(g_n\hat\tau_n)$ and $(g_n\hat\tau'_n)$ converge to {\em distinct} limit simplices.
On the other hand, 
in view of $\hat\tau,\hat\tau'\in C(\tau)$, we have that $\hat\tau_n,\hat\tau'_n\in U_n$ for large $n$. 
Since $(g_nU_n)$ shrinks, it follows that $d(g_n\hat\tau_n,g_n\hat\tau'_n)\to0$, 
a contradiction. 
\qed

\medskip
We need the following general fact from spherical building geometry.

\begin{sublem}
Let $\xi$ be a point in a spherical building $\B$ and let $a\subset \B$ be an apartment. 
If $\xi$ has only one antipode in $a$, 
then $\xi\in a$.
\end{sublem}
\proof
Suppose that $\xi\not\in a$ and let $\hat\xi\in a$ be an antipode of $\xi$. (It always exists.)
We choose a "generic" segment $\xi\hat\xi$ of length $\pi$ tangent to $a$ at $\hat\xi$ as follows.
The suspension $\B(\xi,\hat\xi)\subset \B$ contains an apartment $a'$ with the same unit tangent sphere at $\hat\xi$, $\Si_{\hat\xi}a'=\Si_{\hat\xi}a$. 
Inside $a'$ there exists a segment $\xi\hat\xi$ whose interior does not meet simplices of codimension $\geq2$. 
Hence $\hat\xi\xi$ leaves $a$ at an interior point $\eta\neq\xi,\hat\xi$ of a panel $\pi\subset a$,
i.e.\ $a\cap\xi\hat\xi=\eta\hat\xi$ and $\pi\cap\xi\hat\xi=\eta$,
and $\eta\xi$ initially lies in a chamber adjacent to $\pi$ but not contained in $a$. 
Let $s\subset a$ be the wall (codimension one singular sphere) containing $\pi$. 
By reflecting $\hat\xi$ at $s$, one obtains a second antipode for $\xi$ in $a$, 
contradiction. 
\qed

\medskip
Returning to the proof of the lemma, it follows that $\tau\subset\geo F$.

Sequences $(F_n)$ of maximal flats satisfying our assumptions are obtained as follows.
Since $(\tau_n^-)$ is a shadow of $(g_n^{-1})$ in $\Flagt$,
there exist chambers $\si_n^-\supset\tau_n^-$ 
such that $(\si_n^-)$ is a shadow of $(g_n^{-1})$ in $\D_FX$,
that is 
\begin{equation*}
\label{eq:bdddistwsect}
\sup_n d(g_n^{-1}x,V(x,\si_n^-)) <+\infty
\end{equation*}
for a(ny) point $x\in X$. 
Then the flats $F_n$ extending the euclidean Weyl chambers $V(x,\si_n^-)$ have the property that 
both sequences $(F_n)$ and $(g_nF_n)$ are bounded. 

A subsequence of $(F_n)$ converges iff the corresponding subsequence of $(\si_n^-)$ of chambers converges,
and the limit flat contains the Weyl sector $V(x,\tau_-)$.
In particular, it is itself contained in the parallel set $P(\tau_-,\hat\tau_-)$ 
for the simplex $\hat\tau_-$ opposite to $\tau_-$ with respect to $x$. 
It follows from the above that $\tau\subset\geo P(\tau_-,\hat\tau_-)$.
Since the point $x$ is arbitrary, $\hat\tau_-$ can be any simplex opposite to $\tau_-$,
and we obtain that
\begin{equation*}
\tau\subset\bigcap_{\hat\tau_-\in C(\tau_-)} \geo P(\tau_-,\hat\tau_-) =\st(\tau_-) ,
\end{equation*}
which implies that $\tau=\tau_-$.
\qed

\medskip
Lemma~\ref{lem:shadalmopp} has various implications of asymptotic uniqueness.
Firstly, it yields the asymptotic uniqueness of shadows of the inverse of a contracting sequence. 
From this, in turn, follows regularity 
and we obtain the converse of Proposition~\ref{prop:regcompr}:
\begin{prop}[Contraction implies regularity]
\label{prop:comprreg}
Sequences in $G$, which are contracting on $\Flagt$, are $\tau_{mod}$-regular. 
\end{prop}
\proof
Suppose that the sequence $(g_n)$ in $G$ 
is contracting on $\Flagt$ but not $\tau_{mod}$-regular. 
Then the sequence $(g_n^{-1})$ is not $\tau_{mod}$-regular either and, after passing to a subsequence, 
we may assume that there exists 
a converging sequence $\pi_n^-\to\pi_-$ of panels (codimension one simplices) 
of fixed face type $\theta(\pi_n^-)\not\supset\tau_{mod}$
such that 
\begin{equation}
\label{eq:bdddistsect}
\sup_n d(g_n^{-1}x,V(x,\pi_n^-)) <+\infty
\end{equation}
for a(ny) point $x\in X$. 

Let $\tau_n^-$ be simplices of type $\tau_{mod}$ such that 
$\tau_n^-$ and $\pi_n^-$ are faces of the same chamber.
(This is equivalent to $\tau_n^-\subset\D\st(\pi_n^-)$, respectively, to $\pi_n^-\subset\D\st(\tau_n^-)$.)
Note that $\tau_n^-\not\subset\pi_n^-$ because $\theta(\pi_n^-)\not\supset\tau_{mod}$,
and hence the $\tau_n^-$ are non-unique.
Any such sequence $(\tau_n^-)$ is a shadow of $(g_n^{-1})$ in $\Flagt$.
These sequences $(\tau_n^-)$ can accumulate at any simplex $\tau_-$ 
of type $\tau_{mod}$ contained in $\D\st(\pi_-)$.
Again, since $\theta(\pi_-)\not\supset\tau_{mod}$, there are several such simplices $\tau_-$.
In particular, shadow sequences of $(g_n^{-1})$ are not asymptotically unique. 
This contradicts Lemma~\ref{lem:shadalmopp}. 
\qed

\medskip
Combining Propositions~\ref{prop:regcompr} and~\ref{prop:comprreg}, 
we obtain a characterization
of weak regularity in terms of dynamics at infinity:
\begin{thm}[Contraction characterizes regularity]
\label{thm:regeqcompr}
A sequence in $G$ is $\tau_{mod}$-regular if and only if it is contracting on $\Flagt$. 
\end{thm}

\medskip
Secondly, 
an immediate consequence of Lemma~\ref{lem:shadalmopp} 
is the following asymptotic uniqueness statement for contraction sequences 
complementing the asymptotic uniqueness of expansion and shadow sequences,
cf.\ Lemma~\ref{lem:shadasuniq}:
\begin{lem}[Asymptotic uniqueness of contraction sequences]
\label{lem:contrasuniq}
Suppose that $(U_n)$ and $(U'_n)$ are contraction sequences for $(g_n)$, 
and that 
$(\tau_n)$ and $(\tau'_n)$ are sequences in $\Flagt$ such that 
$(U_n\cap C(\tau_n))$ and $(U'_n\cap C(\tau'_n))$ are still exhaustive.
Then  $d(\tau_n,\tau'_n)\to0$.
\end{lem}

Another consequence of Lemma~\ref{lem:shadalmopp}  
and our earlier discussion of the contraction-expansion dynamics of weakly regular sequences is:
\begin{lem}[Contraction and expansion sequences are asymptotically opposite]
Let $(g_n)$ be a $\tau_{mod}$-regular sequence in $G$. 
Suppose that $(U_n)$ is a contraction and $(V_n)$ an expansion sequence for $(g_n)$. 
Furthermore, let $(\tau_n)$ be a sequence on $\Flagt$ such that $(U_n\cap C(\tau_n))$ 
is still exhaustive. 
Then $V_n\cup\{\tau_n\}$ shrinks, 
$\diam(V_n\cup\{\tau_n\})\to0$.
\end{lem}
\proof
Combine Lemma~\ref{lem:shadalmopp} with 
Proposition~\ref{prop:regcompr} and the proof of Lemma~\ref{lem:shadasuniq}.
\qed

\subsection{Flag convergence}
\label{sec:tauconv}

The asymptotic uniqueness of shadow sequences, 
see Lemma~\ref{lem:shadasuniq},
leads to a notion of convergence at infinity for weakly regular sequences
with limits in the appropriate flag manifolds.
Namely, we can define convergence as the convergence of their shadows:
\begin{dfn}[Flag convergence]
A $\tau_{mod}$-regular sequence in $X$ or $G$ {\em flag converges} 
to a simplex $\tau\in\Flag(\tau_{mod})$, 
if one (any) of its shadow sequences in $\Flag(\tau_{mod})$ converges to $\tau$.
\end{dfn}
We denote the flag convergence of a sequence $(x_n)$ or $(g_n)$ by $x_n\tauto\tau$, respectively, $g_n\tauto\tau$.
If we want to refer to the face type $\tau_{mod}$,
we will sometimes also say that the sequence {\em $\tau_{mod}$-converges}
or speak of flag convergence {\em of type $\tau_{mod}$}.

Note that a $\tau_{mod}$-regular sequence always has a $\tau_{mod}$-converging subsequence
due to the compactness of $\Flag(\tau_{mod})$.

\begin{rem}
(i) 
Flag convergence of type $\tau_{mod}$ captures convergence ``transversely to the stars of the simplices of type $\tau_{mod}$". 
It is related to the ``usual" convergence at infinity 
with respect to the visual compactification $\bar X=X\cup\geo X$ as follows. 
If a $\tau_{mod}$-regular sequence flag converges to $\tau\in\Flagt$, 
then it accumulates in the visual compactification at $\st(\tau)\subset\geo X$;
however, the converse is in general not true.
For {\em uniformly} $\tau_{mod}$-regular sequences, one has equivalence:
They flag converge to $\tau$ if and only if they accumulate in $\bar X$ at $\ost(\tau)$. 

(ii)
Flag convergence of type $\tau_{mod}$ can be understood as convergence at infinity 
with respect to a {\em modified} visual compactification, 
namely as convergence in the $\tau_{mod}$-bordification $X\cup\Flag(\tau_{mod})$ of $X$, 
compare Remark~\ref{rem:regfinscomp}. 

(iii)
$\tau_{mod}$-Convergence implies $\tau'_{mod}$-convergence for smaller face types $\tau'_{mod}\subset\tau_{mod}$,
and the limits correspond under the natural forgetful projection 
$\Flag(\tau_{mod})\to\Flag(\tau'_{mod})$. 
\end{rem}

Flag convergence of weakly regular sequences in $G$ can be characterized in different ways 
in terms of the dynamics at infinity.

Firstly,
one can use the close relation between the asymptotics of orbits in $X$ and the dynamics on flag manifolds,
as expressed by Lemma~\ref{lem:shadagrcontr}.
Since shadow sequences asymptotically agree with the values on contraction subsets,
one has the equivalence:
A $\tau_{mod}$-regular sequence $(g_n)$ in $G$ flag converges
to $\tau\in\Flag(\tau_{mod})$
if and only if for the image sequences $(g_nU_n)$ of its contraction sequences $(U_n)$ in $\Flag(\tau_{mod})$ 
shrink to $\tau$, $g_nU_n\to\tau$.

Secondly,
one can read off flag convergence from the dynamics of the sequence of inverses:

\begin{lem}
\label{lem:charflagconvinv}
For a $\tau_{mod}$-regular sequence $(g_n)$ in $G$ the following two conditions 
are equivalent:

(i)
$g_n\tauto\tau\in\Flagt$. 

(ii)
$(g_n^{-1})$ admits a contraction sequence $(U_n^-)$ in $\Flagt$ opposite to $\tau$,
i.e.\ $U_n^-\subset C(\tau)$.
\end{lem}
\proof
The direction (i)$\Ra$(ii) follows from our construction of contraction sequences:
Let $(\tau_n)$ be a shadow of $(g_n)$ in $\Flagt$. 
Then $\tau_n\to\tau$, and we apply Proposition~\ref{prop:regcompr}(i) to $(g_n^{-1})$.

The reverse direction follows from the asymptotic uniqueness of contraction sequences:
Namely, 
invoking Proposition~\ref{prop:regcompr}(i) again,
there exists another contraction sequence $({U'}_n^-)$ for $(g_n^{-1})$ 
which is opposite to a shadow $(\tau'_n)$ of $(g_n)$.
Lemma~\ref{lem:contrasuniq} then implies that $\tau'_n\to\tau$.
\qed

\medskip
Our next observation 
concerns the relation between flag convergence of sequences in $G$
and their convergence as sequences of maps (homeomorphisms) of flag manifolds.

For a $\tau_{mod}$-converging sequence $g_n\tauto\tau$ in $G$,
one can in general not conclude {\em pointwise} convergence of $(g_n)$ to $\tau$ anywhere on $\Flag(\tau_{mod})$. 
The reason being that in general 
no nested (monotonic) contraction sequence $(U_n)$ exists 
because there is no control on the sequence of open Schubert strata which it approximates;
if $(U_n)$ is opposite to a sequence $(\tau_n^-)$ in $\Flagt$, $U_n\subset C(\tau_n^-)$,
then $(\tau_n^-)$ can be arbitrary.

However, after passing to a subsequence so that 
also $(g_n^{-1})$ flag converges, 
one obtains for $(g_n)$ and $(g_n^{-1})$ 
{\em locally uniform} convergence on open Schubert strata:

\begin{prop}[Attraction-repulsion]
\label{prop:atrrrep}
For a $\tau_{mod}$-regular sequence $(g_n)$ in $G$, 
we have flag convergence $g_n^{\pm1}\tauto\tau_{\pm}\in\Flagt$ 
if and only if 
\begin{equation*}
g_n^{\pm1}\to\tau_{\pm}
\end{equation*} 
locally uniformly on $C(\tau_{\mp})$ as homeomorphisms of $\Flag(\tau_{mod})$. 
\end{prop}
\proof
Suppose that $g_n^{\pm1}\tauto\tau_{\pm}$. 
Then there exist contraction sequences $(U_n^{\pm})$ for $(g_n^{\pm1})$ in $\Flagt$ opposite to $\tau_{\mp}$, 
cf.\ Lemma~\ref{lem:charflagconvinv}.
Moreover, 
$g_n^{\pm1}U_n^{\pm}\to\tau_{\pm}$,
cf.\ Lemma~\ref{lem:shadagrcontr} or our remark above.
Together, this means that 
$g_n^{\pm1}\to\tau_{\pm}$ locally uniformly on $C(\tau_{\mp})$.

Conversely,
if $g_n^{\pm1}\to\tau_{\pm}$ locally uniformly on $C(\tau_{\mp})$,
then there exist contraction sequences $(U_n^{\pm})$ for $(g_n^{\pm1})$ 
such that $g_n^{\pm1}U_n^{\pm}\to\tau_{\pm}$,
and Lemma~\ref{lem:shadagrcontr} implies that $g_n^{\pm1}\tauto\tau_{\pm}$. 
\qed

\begin{rem}
There is no restriction on the relative position of the pair of simplices $\tau_{\pm}$.
They may even agree. 
\end{rem}

\begin{example}\label{ex:convergence}
Fix a maximal flat 
$F\subset X$ and a simplex $\tau\subset\geo F$. Let $(\vartheta_n)$ be a sequence of transvections 
along $F$ such that the expansion factors $\eps(\vartheta_n^{-1},\tau)$ for the action of $\vartheta_n$ on the 
flag  manifold $\Flag(\theta(\tau))=G\tau$ satisfy
$$
\lim_{n\to\infty} \eps(\vartheta_n^{-1},\tau)=+\infty. 
$$
In view of Theorem \ref{thm:infcontrtrans}, 
this is equivalent to $\vartheta_nx\in V(x, \st(\tau)\cap\geo F)$ for large $n$,
with $x\in F$ fixed,
and to the sequence $(\vartheta_n)$ being $\tau_{mod}$-regular, 
where $\tau_{mod}=\theta(\tau)$.
Then $\vartheta_n\tauto\tau$,
because the constant sequence $(\tau)$ is a shadow of $(\vartheta_n)$ in $\Flagt$. 
\end{example}

\subsection{Flag limit sets}
\label{sec:limset}

We now consider discrete subgroups $\Ga\subset G$. 

We recall that the {\em visual limit set} $\La(\Ga)\subset\geo X$ 
is defined as the set of accumulation points of an(y) orbit $\Ga x\subset X$
in the visual compactification $\bar X=X\cup\geo X$,
i.e.\ $\La(\Ga)=\ol{\Ga x}\cap\geo X$.
The notion of flag convergence 
allows to associate to $\Ga$ in an analogous way 
visual limit sets in the flag manifolds associated to $G$. 

\begin{dfn}[Flag limit set]
Let $\Ga\subset G$ be a discrete subgroup. 
We define its {\em flag limit set} of type $\tau_{mod}$ or {\em $\tau_{mod}$-limit set}
$\La_{\tau_{mod}}(\Ga)\subset\Flag(\tau_{mod})$ 
as the set of all limit simplices of $\tau_{mod}$-converging sequences in $\Ga$. 
We call $\La_{\si_{mod}}(\Ga)\subset\D_FX$ 
the {\em chamber limit set} of $\Ga$. 
\end{dfn}
In other words, 
$\La_{\tau_{mod}}(\Ga)$ is the set of accumulation points 
$\ol{\Ga x}\cap\Flag(\tau_{mod})$ 
of an(y) orbit $\Ga x\subset X$
in the $\tau_{mod}$-bordification $X\cup\Flag(\tau_{mod})$. 

\begin{rem}
(i)
$\La_{\tau_{mod}}(\Ga)$ is compact and $\Ga$-invariant. 

(ii)
$\La_{\tau_{mod}}(\Ga)$ is nonempty if and only if $\Ga$ contains a $\tau_{mod}$-regular sequence $\ga_n\to\infty$, 
i.e.\ if part of the $\Ga$-action on $\Flag(\tau_{mod})$ is contracting.
In particular, $\La_{\tau_{mod}}(\Ga)$ is nonempty if $\Ga$ is an infinite $\tau_{mod}$-regular subgroup. 

(iii)
If $\Ga$ is uniformly $\tau_{mod}$-regular, 
then $\La_{\tau_{mod}}(\Ga)$ is the image of $\La(\Ga)$ 
under the natural projection $\geo^{\tau_{mod}-reg}X\to\Flag(\tau_{mod})$. 

(iv)
If $\tau'_{mod}\subset\tau_{mod}$,
then $\La_{\tau_{mod}}(\Ga)$ maps onto $\La_{\tau'_{mod}}(\Ga)$
via the natural forgetful projection 
$\Flag(\tau_{mod})\to\Flag(\tau'_{mod})$. 
\end{rem}

For sufficiently generic subgroups
the limit sets are perfect and the groups act on them minimally:
\begin{prop}[Minimality and perfectness]
\label{prop:miniperf}
Suppose that $\La_{\tau_{mod}}(\Ga)\neq\emptyset$
and that for all $\la\in\La_{\tau_{mod}}(\Ga)$ and $\tau\in\Flag(\tau_{mod})$ we have 
$\Ga\tau\cap C(\la)\neq\emptyset$. 
Then $\La_{\tau_{mod}}(\Ga)$ 
is the unique minimal nonempty $\Ga$-invariant compact subset of $\Flag(\tau_{mod})$. 

If in addition $\La_{\tau_{mod}}(\Ga)$ is infinite, then it is perfect.
\end{prop}
\proof
Let $\la_+\in\La_{\tau_{mod}}(\Ga)$. 
For minimality, 
we must show that $\la_+$ is contained in the closure of every $\Ga$-orbit $\Ga\tau$ in $\Flag(\tau_{mod})$.
By definition of the limit set, 
there exists a $\tau_{mod}$-regular sequence $(\ga_n)$ in $\Ga$ 
such that $\ga_n^{\pm1}\tauto\la_{\pm}$ with some $\la_-\in\La_{\tau_{mod}}(\Ga)$, 
and hence $\ga_n^{\pm1}\to\la_{\pm}$ locally uniformly on $C(\la_{\mp})$. 
By assumption, we have $\Ga\tau\cap C(\la_-)\neq\emptyset$. 
It follows that $\la_+\in\ol{\Ga\tau}$. 

It remains to show perfectness. 
Suppose that $\La_{\tau_{mod}}(\Ga)$ contains isolated points.
The nonisolated limit points (simplices) form a closed subset which, by minimality, must be empty.
Hence, all limit points are isolated and, by compactness, the limit set is finite. 
\qed

\begin{rem}
(i)
The condition that every $\Ga$-orbit on $\Flag(\tau_{mod})$ intersects every open Schubert stratum
is satisfied by {\em Zariski dense} subgroups $\Ga$.
Indeed, 
the action $G\acts\Flag(\tau_{mod})$ is an algebraic action of a semisimple algebraic group on a projective variety. 
In this setting, 
Zariski dense subgroups of $G$ have Zariski dense orbits 
(because the Zariski closure of a subgroup preserves the Zariski closure of any of its orbits). 
So, no orbit can avoid an open Schubert stratum 
because its complement, a Schubert cycle, is a subvariety.
If $\La_{\tau_{mod}}(\Ga)$ is nonempty and finite, 
then $\Ga$ is virtually contained in a parabolic subgroup and can therefore not be Zariski dense. 
Hence, the proposition applies in particular to Zariski dense $\tau_{mod}$-regular subgroups. 

(ii)
Essentially the same notion of limit set 
had been introduced and studied by Benoist in \cite{Benoist}.
He showed that for Zariski dense subgroups $\Ga$ the chamber limit set $\La_{\si_{mod}}(\Ga)$ 
(and hence every $\tau_{mod}$-limit set) is nonempty, 
perfect and the $\Ga$-action on it is minimal,
see \cite[3.6]{Benoist}.
Furthermore, 
the attractive fixed points of the proximal elements in $\Ga$,
cf.\ Remark~\ref{rem:prox},
lie dense in the limit set.
The minimality of the $\Ga$-action on $\La_{\si_{mod}}(\Ga)$ implies 
that the intersection of the visual limit set $\La(\Ga)$ 
with every limit chamber $\si\in\La_{\si_{mod}}(\Ga)$ 
is independent of $\si$ 
modulo the canonical mutual identifications of Weyl chambers,
in other words, the set of types $l(\Ga):=\theta(\si\cap\La(\Ga))\subset\si_{mod}$ is independent. 
One of the main results of \cite{Benoist} regarding the structure of limit sets 
of Zariski dense subgroups $\Ga$ asserts that 
$l(\Ga)$ is convex with nonempty interior, 
see \cite[1.2-1.3]{Benoist}.
(We will not use this result in our paper.)
In particular, if $\Ga$ is uniformly regular and Zariski dense in $G$, then $\La(\Ga)$ is $\Ga$-equivariantly homeomorphic to the product  $l(\Ga)\times \La_{\si_{mod}}(\Ga)$. This product decomposition comes from the fact that 
$\geo^{reg} X$ splits naturally as the product of $\D_F X$ and 
the open simplex $\interior (\si_{mod})$. 
We note that the Zariski density assumption not used in our paper, is essential to Benoist's work. 
On the other hand, $\tau_{mod}$-regularity assumptions  which are key for us, play no role in \cite{Benoist}. 

(iii)
There are other notions of limit sets for actions of discrete subgroups of $G$ on partial flag manifolds, 
see \cite{CNS} for details.
\end{rem}

\begin{rem}
We recall that the limit set of a discrete isometry group of a rank 1 symmetric space either consists of $\le 2$ points or has cardinality of continuum. Situation in the case of symmetric spaces of rank $\ge 2$ is different: It may happen that flag limit set i finite with more than two points. For a specific example, 
let $\Ga\subset G$ be a virtually abelian discrete subgroup which preserves a maximal flat 
$F\subset X$, and acts on $F$ cocompactly and such that the action $\Ga\acts \geo F$ is via the 
Weyl  group of $G$. Then for every face type $\tau_{mod}$ the flag limit set $\LatGa$ is finite and consists of the type $\tau_{mod}$ flags in $\geo F$. Clearly, the group $\Ga$ therefore satisfies the assumptions of the first part of Proposition \ref{prop:miniperf}. 
\end{rem}

\subsection{Antipodal subgroups}

In the remainder of this paper, 
we will only consider weakly regular subgroups 
whose limit set satisfies the following natural additional property.

\begin{dfn}[Antipodal]
\label{dfn:antip}
A subset of $\Flagt$ is {\em antipodal}
if it consists of pairwise opposite simplices.
A discrete subgroup $\Ga\subset G$ 
is {\em $\tau_{mod}$-antipodal} 
if $\La_{\tau_{mod}}(\Ga)$ is antipodal.
\end{dfn}

It would be nice to know whether or not weak regularity implies antipodality. 
We currently have no examples of $\tau_{mod}$-regular discrete subgroup $\Ga\subset G$ 
which are not $\tau_{mod}$-antipodal. 

For antipodal subgroups,
the attraction-repulsion dynamics on flag manifolds as given by Proposition~\ref{prop:atrrrep} 
implies convergence dynamics on the limit set
(also without the assumption of weak regularity).
We obtain the following version of Proposition~\ref{prop:miniperf}:
\begin{prop}[Dynamical properties of antipodal subgroups]
\label{prop:dynantip}
Let $\Ga\subset G$ be a $\tau_{mod}$-antipodal discrete subgroup. 
Then the action $\Ga\acts\La_{\tau_{mod}}(\Ga)$ 
is a convergence action.

If in addition $|\La_{\tau_{mod}}(\Ga)|\geq3$,
then the convergence action $\Ga\acts\La_{\tau_{mod}}(\Ga)$ 
is minimal 
and $\La_{\tau_{mod}}(\Ga)$ is perfect. 
\end{prop}
\proof
The argument is similar to the proof of Proposition~\ref{prop:miniperf}.
Let $\la_+\in\La_{\tau_{mod}}(\Ga)$, 
and let $(\ga_n)$ be a $\tau_{mod}$-regular sequence in $\Ga$ such that  
$\ga_n^{\pm1}\to\la_{\pm}$ locally uniformly on $C(\la_{\mp})$
with $\la_-\in\La_{\tau_{mod}}(\Ga)$. 
Due to antipodality, 
$\La_{\tau_{mod}}(\Ga)\cap C(\la_{\mp})=\La_{\tau_{mod}}(\Ga)-\{\la_{\mp}\}$. 
This establishes the convergence property.

If $\La_{\tau_{mod}}(\Ga)$ has at least three points, 
the minimality of the $\Ga$-action on $\La_{\tau_{mod}}(\Ga)$ and the perfectness of $\La_{\tau_{mod}}(\Ga)$ 
follow from results by Gehring-Martin and Tukia, cf.\ Proposition~\ref{thm:minimal action}. 
\qed

\begin{rem}
Under the assumptions of the proposition,
let $A\subset\Flagt$ be a nonempty $\Ga$-invariant compact antipodal subset. 
Then $\La_{\tau_{mod}}(\Ga)\subset A$, 
compare Proposition~\ref{prop:miniperf}, 
and $\Ga\acts A$ is a convergence action. 
The argument is the same as in the proof of the last proposition.
\end{rem}

\section{Asymptotic conditions for discrete subgroups}\label{sec:Asymptotic properties}

\subsection{Conicality}
\subsubsection{Conical convergence}

Following the notion of conical convergence at infinity for sequences in rank one symmetric spaces,
we will say that a flag converging sequence converges {\em conically} 
if it goes straight towards its limit flag in a suitable sense.
In the context of flag convergence it is natural to require the sequence to stay within bounded distance 
not from a geodesic ray but from the Weyl cone over the star of its limit flag.

\begin{dfn}[Conical convergence]
\label{dfn:conicaltau}
We say that a $\tau_{mod}$-converging sequence in $X$ 
converges {\em conically} to its limit flag $\tau\in\Flagt$
if it is contained in a certain tubular neighborhood of the Weyl cone $V(x,\st(\tau))$. 
Similarly, we say that a $\tau_{mod}$-converging sequence $(g_n)$ in $G$
converges {\em conically} if an(y) orbit $(g_nx)$ in $X$ converges conically.
\end{dfn}

In other words, the convergence of a sequence in $X$ or $G$ is conical with limit flag $\tau$ if and only if 
the constant sequence $(\tau)$ is a shadow sequence.

We will need later the following observation. 
\begin{lem}
\label{lem:clparimplcon}
A flag converging sequence $x_n\tauto\tau$ converges conically 
if and only if it is contained in a tubular neighborhood 
of a parallel set $P(\tau,\hat\tau)$ for some simplex $\hat\tau$ opposite to $\tau$. 
\end{lem}
\proof
Suppose that $(x_n)$ is contained in a tubular neighborhood of $P(\tau,\hat\tau)$.
Then there exist ideal points $\xi_n\in\geo P(\tau,\hat\tau)$ such that $x_n$ has uniformly bounded distance from the rays $x\xi_n$,
where $x\in P(\tau,\hat\tau)$ is some base point. 
Since $\geo P(\tau,\hat\tau)$ is a top-dimensional subbuilding of $\geo X$,
the points $\xi_n$ are contained in chambers $\si_n\subset\geo P(\tau,\hat\tau)$. 
The type $\tau_{mod}$ faces of these chambers then form a shadow sequence of $(x_n)$.
Hence there exists a shadow sequence $(\tau_n)$ of $(x_n)$ 
consisting of simplices $\tau_n\subset\geo P(\tau,\hat\tau)$. 
Since $\tau$ is isolated among the type $\tau_{mod}$ simplices in $\geo P(\tau,\hat\tau)$,
compare Lemma~\ref{lem:poleisol},
the convergence $\tau_n\to\tau$ implies that $\tau_n=\tau$ for large $n$.
\qed

\subsubsection{Recognizing conical convergence at infinity}

As it is the case for flag convergence,
also the conicality of the convergence can be read off the dynamics at infinity,
namely it can be recognized from the dynamics on pairs of flags. 

Suppose that the sequence $(g_n)$ in $G$ flag converges, $g_n\tauto\tau$.
Then there exists a contraction sequence for $(g_n^{-1})$ which exhausts $C(\tau)$, 
cf.\ Lemma~\ref{lem:charflagconvinv}.
Hence  the orbits $(g_n^{-1}\hat\tau)$ for all simplices $\hat\tau\in C(\tau)$ 
are asymptotic to each other.
(More generally, 
this remains true for all sequences $(g_n^{-1}\hat\tau_n)$ where $(\hat\tau_n)$ is a bounded sequence in $C(\tau)$.)
However,
there is in general no control on the asymptotics of the orbits of the other simplices $\tau'\in\Flagt-C(\tau)$, 
in particular of the orbit $(g_n^{-1}\tau)$.

We will now see that $g_n\tauto\tau$ conically
if and only if the orbit $(g_n^{-1}\tau)$ stays away from the orbits $(g_n^{-1}\hat\tau)$ for $\hat\tau\in C(\tau)$ 
in the sense that the sequence of pairs $g_n^{-1}(\tau,\hat\tau)$ is bounded in the space 
\begin{equation*}
(\Flag(\tau_{mod})\times\Flag(\tau_{mod}))^{opp} \subset \Flag(\tau_{mod})\times\Flag(\tau_{mod})
\end{equation*}
of pairs of opposite flags of type $\tau_{mod}=\theta(\tau)$.

\begin{lem}
\label{lem:opppairsparsets}
A set of opposite pairs $A\subset (\Flag(\tau_{mod})\times\Flag(\tau_{mod}))^{opp}$ is bounded 
if and only if the corresponding family of parallel sets $P(\tau_-,\tau_+)$ in $X$ for $(\tau_-,\tau_+)\in A$ is bounded,
i.e.\ 
\begin{equation*}
\sup_{(\tau_-,\tau_+)\in A} d(x,P(\tau_-,\tau_+)) < +\infty
\end{equation*}
for a base point $x$.
\end{lem}
\proof
The forward direction can be deduced from the fact that 
$(\Flag(\tau_{mod})\times\Flag(\tau_{mod}))^{opp}$ is a homogeneous $G$-space.
Indeed, if $A$ is bounded, 
then there exists a compact subset $C\subset G$ and a reference pair $a_0=(\tau_0^-,\tau_0^+)$ such that 
$A\subset Ca_0$. 
It follows that the parallel sets $P(\tau_-,\tau_+)$ for $(\tau_-,\tau_+)\in A$ intersect the compact set $Cx_0$,
where $x_0$ is a point in $P(\tau_0^-,\tau_0^+)$.

For the converse direction we use that the set of triples $(\tau_-,x',\tau_+)$, such that 
$\tau_{\pm}$ are simplices of type $\tau_{mod}$ opposite to each other with respect to the point $x'\in X$,
is still a homogeneous $G$-space. 
As a consequence, every parallel set $P(\tau_-,\tau_+)$ intersecting the ball $B(x_0,R)$ 
is of the form $gP(\tau_0^-,\tau_0^+)$ with $g$ such that $d(x_0,gx_0)\leq R$, 
i.e.\ $g$ belongs to a compact subset of $G$.
It follows that the set of these pairs $(\tau_-,\tau_+)=g(\tau_0^-,\tau_0^+)$ 
is bounded in $(\Flag(\tau_{mod})\times\Flag(\tau_{mod}))^{opp}$. 
\qed

\begin{prop}[Recognizing conical convergence at infinity]
\label{prop:recogconconvinf}
Suppose that the sequence $(g_n)$ in $G$ flag converges, $g_n\tauto\tau$.
Then the following are equivalent:

(i) $g_n\tauto\tau$ conically.

(ii)
For some flag $\hat\tau\in C(\tau)$,
the sequence of pairs $g_n^{-1}(\tau,\hat\tau)$ is bounded in $(\Flag(\tau_{mod})\times\Flag(\tau_{mod}))^{opp}$. 

(ii')
For any bounded sequence $(\hat\tau_n)$ in $C(\tau)$,
the sequence of pairs $g_n^{-1}(\tau,\hat\tau_n)$ is bounded in $(\Flag(\tau_{mod})\times\Flag(\tau_{mod}))^{opp}$. 
\end{prop}
\proof
Let $(\hat\tau_n)$ be a bounded sequence in $C(\tau)$,
i.e.\ as a subset it is relatively compact. 
Then the set of pairs 
$(\tau,\hat\tau_n)$ is bounded 
in $(\Flag(\tau_{mod})\times\Flag(\tau_{mod}))^{opp}$,
and hence the family of parallel sets $P(\tau,\hat\tau_n)$ is bounded, i.e.\
\begin{equation*}
\sup_n d(x,P(\tau,\hat\tau_n)) < +\infty 
\end{equation*}
for a base point $x$,
cf.\ Lemma~\ref{lem:opppairsparsets}. 
We have the estimate 
\begin{equation*}
d(x,g_n^{-1}P(\tau,\hat\tau_n)) = d(g_nx,P(\tau,\hat\tau_n)) \leq d(g_nx,V(x,\st(\tau)) + d(x,P(\tau,\hat\tau_n)) ,
\end{equation*}
using that $\st(\tau)\subset\geo P(\tau,\hat\tau_n)$ and hence 
$d(\cdot,P(\tau,\hat\tau_n))|_{V(x,\st(\tau))}$ is maximal at $x$.

The right hand side is bounded iff $g_n\to\tau$ conically.
The left hand side is bounded iff
the sequence of pairs $g_n^{-1}(\tau,\hat\tau_n)$ is bounded in $(\Flag(\tau_{mod})\times\Flag(\tau_{mod}))^{opp}$, 
again by Lemma~\ref{lem:opppairsparsets}. 
This shows the implication (i)$\Ra$(ii').

Conversely, assume the weaker condition (ii). 
Then the sequence $(g_nx)$ is contained in a tubular neighborhood of $P(\tau,\hat\tau)$ 
and Lemma~\ref{lem:clparimplcon} implies (i).
\qed

\subsubsection{Conical limit set}

Conicality is a condition on the asymptotic geometry of the orbits of a discrete subgroup $\Ga\subset G$,
namely on how limit flags can be approached by sequences in orbits $\Ga x\subset X$.

Following the definition of conicality in the rank one case, 
cf.\ section~\ref{sec:rank1}, 
a limit point $\xi\in\La(\Ga)\subset\geo X$ may be called {\em ray conical} 
if it is the limit of a sequence of orbit points $\ga_nx$ 
which are contained in a tubular neighborhood 
of some geodesic ray asymptotic to $\xi$. 
However, ray conicality too restrictive in higher rank.
It is satisfied by con\-vex-co\-comp\-act subgroups, 
but these are rare, cf.\ \cite{convcoco},
and e.g.\ one can show that RCA Schottky subgroups 
are, in general, not ray conical, see section~\ref{sec:schottky actions} for the construction. 

The following notion of conicality 
considered by Albuquerque 
\cite[Def.\ 5.2]{Albuquerque}
is more flexible
and useful in higher rank. 

\begin{dfn}[Conical limit set]
\label{dfn:conlim}
Let $\Ga\subset G$ be a discrete subgroup. 
We call a limit flag $\la\in\La_{\tau_{mod}}(\Ga)$ {\em conical} 
if there exists a $\tau_{mod}$-regular sequence $(\ga_n)$ in $\Ga$ 
such that $\ga_n\tauto\la$ conically.
The {\em conical $\tau_{mod}$-limit set} $\La_{\tau_{mod}}^{con}(\Ga)\subset\La_{\tau_{mod}}(\Ga)$
is the subset of conical limit flags.
We say that $\Ga$ has {\em conical $\tau_{mod}$-limit set} 
or is {\em $\tau_{mod}$-conical} 
if $\La_{\tau_{mod}}(\Ga)=\La_{\tau_{mod}}^{con}(\Ga)$.
\end{dfn}

We deduce from Proposition~\ref{prop:recogconconvinf}
how one can recognize conical limit flags 
from the dynamics at infinity,
compare the notion of intrinsically conical point in the case of convergence actions
cf.\ Definition~\ref{dfn:intrcon}.

\begin{lem}[Recognizing conical limit flags]
\label{lem:charconinf}
A flag $\tau\in\Flagt$ is a conical limit flag if and only if there exists a $\tau_{mod}$-regular sequence 
$(\ga_n)$ in $\Ga$ such that 

(i) the maps $\ga_n^{-1}|_{C(\tau)}$ converge locally uniformly to a constant map, 
$\ga_n^{-1}|_{C(\tau)}\to\hat\tau_-$, 
and 

(ii)
the points $\ga_n^{-1}\tau$ converge, $\ga_n^{-1}\tau\to\tau_-$, 
with limit $\tau_-$ opposite to $\hat\tau_-$.
\end{lem}
\proof
Let $\tau$ be a conical limit flag. 
Then there exists a $\tau_{mod}$-regular sequence $(\ga_n)$ in $\Ga$ 
such that $\ga_n\tauto\tau$ conically. 
By passing to a subsequence we can obtain further convergence properties,
namely that also the sequence $(\ga_n^{-1})$ flag converges 
(not necessarily conically),
$\ga_n^{-1}\tauto\hat\tau_-$,
and that $\ga_n^{-1}\tau\to\tau_-$. 
Then Proposition~\ref{prop:atrrrep} yields that 
$\ga_n^{-1}|_{C(\tau)}\to\hat\tau_-$ locally uniformly.

To see that the flags $\tau_-$ and $\hat\tau_-$ are opposite,
we use Proposition~\ref{prop:recogconconvinf}.
For any flag $\hat\tau\in C(\tau)$ we have the convergence of pairs
$\ga_n^{-1}(\tau,\hat\tau)\to(\tau_-,\hat\tau_-)$.
On the other hand, the proposition implies that the
sequence $\ga_n^{-1}(\tau,\hat\tau)$ is relatively compact
in the space of opposite pairs
$(\Flag(\tau_{mod})\times\Flag(\tau_{mod}))^{opp}$.
Hence the limit pair $(\tau_-,\hat\tau_-)$ must also lie in this space, 
i.e.\ $\tau_-$ and $\hat\tau_-$ are opposite.

Conversely, suppose that $(\ga_n)$ is a $\tau_{mod}$-regular sequence in $\Ga$ satisfying (i) and (ii).
Then property (i) implies flag convergence $\ga_n\tauto\tau$, compare Lemma~\ref{lem:charflagconvinv},
and we can apply Proposition~\ref{prop:recogconconvinf}.
It follows that $\ga_n\tauto\tau$ conically. 
\qed

\begin{rem}
As mentioned in the proof, 
property (i) implies flag convergence $\ga_n\tauto\tau$, compare Lemma~\ref{lem:charflagconvinv},
and in particular that $\tau$ is a limit flag, $\tau\in\LatGa$.
\end{rem}

\subsubsection{Comparing extrinsic and intrinsic conicality}

If $\Ga\subset G$ is a $\tau_{mod}$-antipodal discrete subgroup,
then $\Ga\acts\LatGa$ is a convergence action, cf.\ Proposition~\ref{prop:dynantip},
and hence there is a notion of intrinsic conical point for this action, 
cf.\ Definition~\ref{dfn:intrcon}.
We show now that the intrinsic and extrinsic notions of conicality coincide
for non-elementary weakly regular antipodal subgroups:

\begin{prop}
[Conical equivalent to intrinsically conical]
\label{prop:concon}
Let $\Ga\subset G$ be a $\tau_{mod}$-antipodal $\tau_{mod}$-regular
discrete subgroup
and suppose that $|\LatGa|\geq3$.

Then a limit flag in $\LatGa$ is conical 
if and only if it is intrinsically conical 
for the convergence action $\Ga\acts\LatGa$.
\end{prop}
\proof
Suppose that the limit flag $\la\in\LatGa$ is conical.
Then it follows by restricting the $\Ga$-action on $\Flagt$ to the convergence action on $\LatGa$
and applying Lemma~\ref{lem:charconinf},
that $\la$ is intrinsically conical. 

Conversely, 
suppose that $\la$ is intrinsically conical for the convergence action $\Ga\acts\LatGa$.
Then there exist a sequence $\ga_n\to\infty$ in $\Ga$ 
and a pair 
of opposite limit flags in $\la_-,\hat\la_-\in\LatGa$ such that
$\ga_n^{-1}\la\to\la_-$ and 
$\ga_n^{-1}|_{\LatGa-\{\la\}}\to\hat\la_-$ locally uniformly.
Note that $\LatGa-\{\la\}\neq\emptyset$ because $|\LatGa|\geq3$.

Now we consider the $\Ga$-action on the entire flag manifold $\Flagt$.
By assumption,
$\Ga$ is $\tau_{mod}$-regular,
and hence in particular the sequence $(\ga_n)$. 
After passing to a subsequence, we may assume that $\ga_n^{\pm1}\tauto\tau_{\pm}\in\Flagt$,
and therefore
\begin{equation*}
\ga_n^{\pm1}|_{C(\tau_{\mp})}\to\tau_{\pm}
\end{equation*} 
locally uniformly, 
as a consequence of Proposition~\ref{prop:atrrrep}.
Necessarily, $\tau_{\pm}\in\LatGa$.

We first observe that $\tau_+=\la$.
Indeed, assume that $\tau_+\neq\la$.
Then $\LatGa-\{\la\}$ and $C(\tau_+)$ cover $\LatGa$.
Since $|\LatGa|\geq3$,
they also intersect 
and it follows that $\ga_n^{-1}$ converges uniformly on the entire limit set $\LatGa$ to a constant map,
which is absurd.
So, $\ga_n\tauto\la$.

Now we can apply Proposition~\ref{prop:recogconconvinf}.
For every limit flag $\hat\la\in\LatGa-\{\la\}\subset C(\la)$ we have 
$\ga_n^{-1}(\la,\hat\la)\to(\la_-,\hat\la_-)$.
Such limit flags $\hat\la\neq\la$ exist,
and implication (ii)$\Ra$(i) of the proposition yields that $\ga_n\tauto\la$ conically,
i.e.\ $\la$ is conical.
\qed

\subsubsection{Expansion at conical limit flags}\label{sec: Expansion at conical limit flags}

If a weakly regular sequence $(g_n)$ in $G$ flag converges, 
$g_n\tauto\tau\in\Flagt$, 
then its sequence $(g_n^{-1})$ of inverses admits 
an expansion sequence $(V_n^-)$ in $\Flagt$ 
with $V_n^-\to\tau$,
see Proposition~\ref{prop:regcompr}(ii).
We will show now that, 
if the convergence is conical, 
then there is a stronger form of expansion at $\tau$ 
for the dynamics of $(g_n^{-1})$ on $\Flagt$.

Generalizing the definition of expansion point for a group action,
cf.\ Definition~\ref{def:expanding}, 
from groups to sequences, 
we say that the sequence $(g_n^{-1})$ in $G$
is {\em expanding} on $\Flagt$ at $\tau$ if
$$
\lim_{n\to+\infty}\eps(g_n^{-1},\tau)=+\infty
$$
with respect to an auxiliary Riemannian metric on $\Flagt$.
This implies in particular 
that for large $n$
the map $g_n^{-1}$ is uniformly expanding on some neighborhood $V_n$ of $\tau$ in $\Flagt$ 
with expansion factor $c_n\to+\infty$.

\begin{lem}
\label{lem:expconlim}
If a $\tau_{mod}$-regular sequence $(g_n)$ in $G$ flag converges conically, $g_n\tauto\tau\in\Flagt$,
then the sequence $(g_n^{-1})$ is expanding at $\tau$. 
\end{lem}
\proof
By assumption,
the orbit sequence $(g_nx)$ for a point $x\in X$ is contained in a tubular neighborhood of the Weyl cone $V(x,\st(\tau))$
and, due to $\tau_{mod}$-regularity, we have that 
$$
\lim_{n\to+\infty}d(g_nx,\D V(x,\st(\tau)))=+\infty.
$$
Corollary~\ref{cor:expand} implies that the infinitesimal expansion of $g_n^{-1}$ at $\tau$ 
becomes arbitrarily strong,
$$
\lim_{n\to+\infty}\eps(g_n^{-1},\tau)=+\infty
$$
(with respect to the auxiliary Riemannian metric on $\Flag(\tau_{mod})$). 
Thus for sufficiently large $n$, 
the map $g_n^{-1}$ 
is uniformly expanding with arbitrarily large expansion factor on some neighborhood of $\tau$. 
\qed

\medskip
Applied to group actions, the lemma yields:
\begin{prop}[Expansion at conical limit flags]
\label{prop:conicimplexpatau}
Let $\Ga\subset G$ be a discrete subgroup. 
If the limit flag $\la\in\La_{\tau_{mod}}(\Ga)$ is conical, 
then the action $\Ga\acts\Flag(\tau_{mod})$ is expanding at $\la$. 
In particular, 
if $\Ga$ has conical flag limit set $\La_{\tau_{mod}}(\Ga)$, 
then the action $\Ga\acts\Flag(\tau_{mod})$ is expanding at 
$\La_{\tau_{mod}}(\Ga)$. 
\end{prop}

\subsection{Equivalence of certain asymptotic conditions}

In section~\ref{sec:rank1}, we discussed discrete groups of isometries on symmetric spaces of rank one 
and formulated a number of conditions (C2-C10) which are equivalent to convex cocompactness. 
We will now generalize some of these conditions to weakly regular discrete groups 
in arbitrary rank 
and show that they remain equivalent to each other. 

The first condition generalizes the conical limit set property (C2) in rank one:
\begin{dfn}[RCA]
\label{dfn:rcatau}
We call a discrete subgroup of $G$ {\em $\tau_{mod}$-RCA} 
if it is $\tau_{mod}$-regular, $\tau_{mod}$-conical and $\tau_{mod}$-antipodal. 
We call it {\em weakly RCA} if it is $\tau_{mod}$-RCA for some $\tau_{mod}$, and 
{\em RCA} if it is $\si_{mod}$-RCA. 
\end{dfn}

In this section we prove that weakly RCA groups are word hyperbolic 
and that their $\tau_{mod}$-limit set is equivariantly homeomorphic to their Gromov boundary.
We will also prove a converse of this result and establish a similar equivalence of between weak RCA  and an expansion property (subject to $\tau_{mod}$-regularity and  antipodality conditions). 

The second condition generalizes condition C6 in rank one, 
requesting that the subgroup is intrinsically word hyperbolic 
and its limit set is an equivariantly embedded copy of its Gromov boundary as a word hyperbolic group: 
\begin{dfn}[Asymptotically embedded]
\label{defn:hyperbolictau}
We call an $\tau_{mod}$-antipodal $\tau_{mod}$-regular discrete subgroup $\Ga\subset G$ 
{\em $\tau_{mod}$-asymptotically embedded} 
if $\Gamma$ is word hyperbolic 
and there exists a $\Ga$-equivariant homeomorphism 
\begin{equation}
\label{eq:mapalphatau}
\alpha: \geo \Gamma \buildrel\cong\over\to 
\La_{\tau_{mod}} (\Ga)\subset \Flag(\tau_{mod})
\end{equation}
of its Gromov boundary onto its $\tau_{mod}$-limit set. 
\end{dfn}
Note that in view of the minimality of the action on the flag limit set, cf.\ Proposition~\ref{prop:dynantip}, 
it suffices to assume that $\al$ is an equivariant embedding into $\La_{\tau_{mod}}(\Ga)$.

The terminology ``asymptotically embedded'' is justified by the next observation 
which can be understood as saying that for asymptotically embedded subgroups the orbit maps $\Ga\to\Ga x\subset X$ 
continuously extend to maps $\bar\Ga\to X\cup\Flagt$ 
from the visual compactification of the group to the bordification of $X$:
\begin{lem}[Continuity at infinity]
\label{lem:continuity_at_infinity}
Suppose that $\Ga\subset G$ is $\tau_{mod}$-asymptotically embedded 
and $|\LatGa|\geq3$. 
Then a sequence $\ga_n\to\infty$ in $\Ga$ 
converges to $\zeta\in\geo\Ga$ 
if and only if, as a sequence in $G$, it flag converges to 
$\al(\zeta)\in\La_{\tau_{mod}}(\Ga)$. 
\end{lem}
Note that the assertion may fail when $\Ga$ is elementary (as a hyperbolic group),
for instance when $\Ga$ is cyclic and hence acts trivially on its ideal boundary $\geo\Ga$
which consists of two points. 

\proof
We use the characterization of flag convergence in terms of the dynamics at infinity 
and the analogous fact for hyperbolic groups. 

Suppose that the assertion is wrong. 
Then there exists a sequence $\ga_n\to\infty$ in $\Ga$ 
such that $\ga_n\to\zeta$ in $\bar\Ga$ 
and $\ga_n\tauto\la\in\LatGa$ in $G$,
but $\la\neq\al(\zeta)$.
According to Lemma~\ref{lem:charflagconvinv},
this means that $(\ga_n^{-1})$
admits contraction sequences $(U_n)$ on $\Flagt$ 
opposite to $\la$, $U_n\subset C(\la)$,
and $(U'_n)$ on $\geo\Ga$ opposite to $\zeta$, 
$U'_n\subset\geo\Ga-\{\zeta\}$.
Due to antipodality, 
$(\al^{-1}(U_n))$ is another contraction sequence on $\geo\Ga$,
but opposite to $\al^{-1}(\la)\neq\zeta$.
It follows that for large $n$ the subsets $U'_n$ and $\al^{-1}(U_n)$ cover $\geo\Ga$ and,
since $|\geo\Ga|=|\LatGa|\geq3$, intersect.
This implies that $\ga_n^{-1}(\geo\Ga)$ shrinks to a point,
a contradiction. 
\qed

\medskip
Our third condition extends the 
expansion at the limit set property (C10) to higher rank:
\begin{dfn}[Expanding at infinity]\label{def:expanding at infinity}
We call a $\tau_{mod}$-antipodal $\tau_{mod}$-regular discrete subgroup $\Ga\subset G$ 
{\em $\tau_{mod}$-expanding at infinity}
if the action $\Ga\acts\Flag(\tau_{mod})$ is expanding at $\La_{\tau_{mod}}(\Ga)$ 
(with respect to a Riemannian background metric). 
\end{dfn}

We note that it is plausible that in the three definitions above the antipodality assumption 
is implied by regularity and hence redundant.

\begin{thm}[Equivalence of asymptotic conditions]
\label{thm:regopphyp}
Suppose that $\Ga\subset G$ is a $\tau_{mod}$-antipodal 
$\tau_{mod}$-regular discrete subgroup such that $|\LatGa|\geq3$.
Then the following are equivalent: 

(i) $\Ga$ is $\tau_{mod}$-RCA. 

(ii) $\Ga$ is $\tau_{mod}$-asymptotically embedded. 

(iii) $\Ga$ is $\tau_{mod}$-expanding at infinity.
\end{thm}
\proof
Due to the equivalence of the intrinsic and extrinsic notions of conicality, 
see Proposition~\ref{prop:concon}, 
property (i) is equivalent to $\La_{\tau_{mod}}(\Ga)$ 
being intrinsically conical for the convergence action $\Ga\acts\LatGa$.

Suppose that $\La_{\tau_{mod}}(\Ga)$ is intrinsically conical.
According to Proposition~\ref{prop:dynantip},
$\La_{\tau_{mod}}(\Ga)$ is perfect. 
Therefore, 
we can apply 
Bowditch's Theorems~\ref{thm:bowditch-conical} and~\ref{thm:charhypbow}. 
($\La_{\tau_{mod}}(\Ga)$ is metrizable as a subset of the manifold 
$\Flag(\tau_{mod})$.) 
Theorem~\ref{thm:bowditch-conical} implies that 
the convergence action $\Ga\acts\La_{\tau_{mod}}(\Ga)$ is uniform, 
and Theorem~\ref{thm:charhypbow} implies property (ii). 
The converse implication holds because 
the action of a word hyperbolic group on its Gromov boundary 
is intrinsically conical. 
Hence (i) and (ii) are equivalent. 

The equivalence with (iii) can be seen as follows.
The (extrinsic) conicality of $\LatGa$ implies that 
the action $\Ga\acts\Flag(\tau_{mod})$ is expanding at 
$\La_{\tau_{mod}}(\Ga)$,
see Lemma~\ref{lem:expconlim}.
Vice versa,
if the action $\Ga\acts\LatGa$ is expanding, 
then $\La_{\tau_{mod}}(\Ga)$ is intrinsically conical by Lemma \ref{lem:conical}. 
\qed

\begin{cor}
\label{cor:fingen}
If $\Ga\subset G$ is $\tau_{mod}$-RCA or $\tau_{mod}$-expanding at infinity, 
then it is word hyperbolic.
\end{cor}

\subsection{Boundary embeddings}

The asymptotic embedding property requires an equivariant embedding from the Gromov 
boundary  to the $\tau_{mod}$-limit set of an intrinsically word hyperbolic $\tau_{mod}$-regular discrete subgroup. 
In this section we consider a weakening of this property. 

\begin{dfn}[Boundary embedded]
\label{dfn:nicastau}
Let $\Ga$ be a non-elementary (i.e.\ not virtually cyclic) word hyperbolic group. 
We say that an isometric action $\rho:\Ga\acts X$ 
is {\em $\tau_{mod}$-boundary embedded} if there exists a $\Ga$-equivariant continuous embedding  
\begin{equation}
\label{eq:mapbeta}
\beta: \geo \Gamma \to \Flag(\tau_{mod}).
\end{equation}
which maps different boundary points to opposite flags. 
\end{dfn}
Note that such actions are necessarily properly discontinuous, 
because $\Ga$ acts on $\beta(\geo\Ga)$ as a discrete convergence group. 
Henceforth, we identify $\Ga$ with its image in $G$. 

Note moreover, that we do not a priori assume that the subgroup $\Ga\subset G$ 
is $\tau_{mod}$-regular. 
Even if it is $\tau_{mod}$-regular, 
boundary embeddedness is a priori weaker 
than asymptotic embeddedness, because  
it does not assume that 
$\beta(\geo\Ga)\subseteq\La_{\tau_{mod}}(\Ga)$
or that $\Ga$ is $\tau_{mod}$-antipodal. Nevertheless, 
in the case $\tau_{mod}=\sigma_ {mod}$ we will show that being boundary embedded and regular implies the
stronger asymptotic embeddedness property,
see Proposition~\ref{prop:homeochlim} below. Furthermore, one can show   
the same implication for arbitrary $\tau_{mod}$ and Zariski dense subgroups
(cf. \cite{GW}), but we will not prove this in our paper. 

\medskip 
Note also that for  $\tau_{mod}$-asymptotically embedded groups there are in general other equivariant embeddings 
$\geo\Gamma \to \Flag(\tau_{mod})$
besides the one onto the $\tau_{mod}$-limit set, even if $\tau_{mod}=\sigma_{mod}$:  

\begin{example}[Nonuniqueness of boundary maps]
\label{ex:alnebe}
One can construct totally geodesic embeddings $Y\embed X$ of symmetric spaces, 
e.g.\ of equal rank $\geq2$, 
such that for the induced boundary map at infinity 
$\geo Y\embed\geo X$
Weyl chambers of $Y$ break up into several Weyl chambers of $X$. 
Then there are several induced embeddings 
$\D_FY\embed\D_FX$ of F\"urstenberg boundaries.
As a consequence, 
for suitable hyperbolic groups $\Ga$ acting on $Y$ 
one obtains several equivariant embeddings 
$\geo\Ga\embed\D_FX$. 

For instance, 
consider the embedding of Weyl groups 
$W_{A_1\circ A_1}\subset W_{B_2}$ 
and the corresponding refinement of Coxeter complexes 
where an $A_1\circ A_1$-Weyl arc of length $\pihalf$ 
breaks up into two $B_2$-Weyl arcs of length $\piforth$. 
Let $Y\embed X$ be an isometric embedding of symmetric spaces
inducing this embedding of Weyl groups, 
and $H\subset G$ a corresponding embedding of semisimple Lie groups, 
for instance, 
$SO(2,1)\times SO(2,1)\subset SO(4,2)$. 
The symmetric space $Y$ is reducible 
and decomposes as a product of rank one spaces, 
$Y\cong Y_1\times Y_2$. 
Accordingly, its Tits boundary decomposes as a spherical join, 
$\tits Y\cong\tits Y_1\circ\tits Y_2$. 
The Weyl chambers of $Y$ are arcs $\xi_1\circ\xi_2$ of length $\pihalf$ 
with endpoints $\xi_i\in\geo Y_i$. 
The embedding $\geo Y\embed\geo X$ of visual boundaries 
sends the Weyl chamber $\xi_1\circ\xi_2$ to an arc 
denoted by $\xi_1\xi_2$. 
With respect to the spherical building structure on $\geo X$ 
it decomposes as the union of two Weyl chambers 
$\xi_i\mu$ of length $\piforth$, 
where $\mu$ is the midpoint of $\xi_1\xi_2$. 
We see that there are {\em two} $H$-equivariant continuous embeddings 
of F\"urstenberg boundaries 
$\iota_i:\D_FY\embed\D_FX$ 
obtained by assigning to each Weyl chamber $\xi_1\circ\xi_2$ of $Y$ 
its half $\xi_i\mu$ of type $i$. 
The embeddings send opposite chambers to opposite chambers. 
It is easy to construct regular Schottky subgroups $\Ga\subset H$ 
which remain regular in $G$, 
and by composing the embeddings $\geo\Ga\to\D_FY$ with $\iota_i$, 
one obtains two $\Ga$-equivariant embeddings $\geo\Ga\to\D_FX$ 
mapping distinct boundary points to opposite chambers. 
\end{example}

\subsection{Coarse extrinsic geometry}

In this section we study the coarse geometry of discrete subgroups $\Ga\subset G$
satisfying one of the asymptotic conditions introduced above,
i.e.\ one of the three equivalent conditions 
``RCA", ``asymptotically embedded" and ``expanding at infinity"
or the weaker condition ``boundary embedded".
Note that such subgroups are intrinsically word hyperbolic and hence finitely generated, 
see Corollary~\ref{cor:fingen} and Definition~\ref{dfn:nicastau}.
We will show that RCA subgroups are undistorted, i.e. that the orbit maps $\Ga\to\Ga x\subset X$ are quasi-isometric embeddings. 
Equivalently, 
they send uniform quasigeodesics in $\Ga$ to uniform quasigeodesics in $X$. 
We will in fact prove a stronger form of undistortion,
namely that 
the images of quasigeodesics in $\Ga$ under the orbit maps 
satisfy a generalized version of the Morse Lemma (for quasigeodesics in negatively curved spaces):
the images of quasirays stay close to Weyl cones. 
This is indeed a strong further restriction 
because in higher rank quasigeodesics are quite flexible. 

\subsubsection{Boundary embedded groups}
\label{sec:bdembgp}

In this section we consider boundary embedded groups; this is a weakening of the asymptotic embeddedness which will be, however, sufficient for establishing some some preliminary control 
on the images of quasigeodesics in $\Ga$ 
under the orbit maps, namely that they are uniformly close to parallel sets.

Let $\Ga\subset G$  be a $\tau_{mod}$-boundary embedded discrete subgroup.
The boundary map $\beta:\geo\Gamma \to \Flag(\tau_{mod})$ induces a map $(\beta,\beta)$ 
of pairs of boundary points. 
By assumption, $(\beta,\beta)$ maps pairs of distinct boundary points into 
the  open dense $G$-orbit 
\begin{equation*}
(\Flag(\tau_{mod})\times\Flag(\tau_{mod}))^{opp}
\subset\Flag(\tau_{mod})\times\Flag(\tau_{mod})
\end{equation*}
consisting of the pairs of opposite flags,
and we obtain a continuous $\Ga$-equivariant embedding
\begin{equation*}
\label{eq:leafmapmov}
(\geo\Ga\times\geo\Ga)-\mathrm{Diag} \xrightarrow{(\beta,\beta)}
(\Flag(\tau_{mod})\times\Flag(\tau_{mod}))^{opp} .
\end{equation*}
Here and in what follows, $\mathrm{Diag}$ denotes the diagonal in the product. 

\begin{lem}
\label{lem:qgeoflat}
An $(L,A)$-quasigeodesic $q:\Z\to\Ga$ with ideal endpoints 
$\zeta_{\pm}\in\geo\Ga$ 
is mapped by the orbit map $\Ga\to\Ga x\subset X$ 
into a tubular neighborhood 
of uniform radius $r=r(\Ga,L,A,x)$
of the parallel set $P(\beta(\zeta_-),\beta(\zeta_+))$. 
\end{lem}
\proof
We consider the map from the space of  $(L,A)$-quasigeodesics
$q:\Z\to\Ga$ to the space
$$
(\Flag(\tau_{mod})\times\Flag(\tau_{mod}))^{opp}\times X,
$$
assigning to $q$ the pair consisting of the (directed) parallel set
$P(\beta(q(-\infty) ),\beta(  q(+\infty) ))$ and the orbit point
$q(0) x$.
This map is $\Gamma$-equivariant and continuous, where the space of $(L,A)$-quasigeodesics
is equipped with the topology of pointwise convergence.
(The continuity of the assignment $q\mapsto (q(-\infty), q(+\infty))$ uses 
the Morse Lemma in the hyperbolic group $\Ga$.)
By composing this map with the distance  between  $q(0) x$ and the parallel set,
we obtain a $\Gamma$-periodic continuous function on the space of $(L,A)$-quasigeodesics.
Since $\Gamma$ acts cocompactly on this space (by Arzela-Ascoli), this function is bounded.
The assertion follows by shifting the parametrization of $q$.
\qed

\subsubsection{The regular case}

We restrict now to the case $\tau_{mod}=\sigma_{mod}$,
i.e.\ we assume that the subgroup $\Ga$ is $\sigma_{mod}$-boundary embedded. In particular, it is then regular. Our proofs, while less general, will be more straightforward and motivate the more difficult arguments for general $\tau_{mod}$ 
and asymptotically embedded subgroups in section \ref{sec:coarse_asymptotically_embedded} below.

According to Lemma~\ref{lem:qgeoflat}, the images of quasigeodesics in $\Ga$ under the orbit maps 
are now uniformly close to maximal flats; 
a quasigeodesics asymptotic to a pair of ideal points $\zeta_{\pm}\in\geo\Ga$
is mapped into a tubular neighborhood of the maximal flat $F(\beta(\zeta_-),\beta(\zeta_+))$ 
asymptotic to the pair of opposite chambers $ \beta(\zeta_-)$ and $\beta(\zeta_+)$.
The next result restricts the position of the image along the maximal flat.
Namely, the images of quasirays are uniformly close to euclidean Weyl chambers
and move towards limit chambers at infinity:

\begin{lem}
\label{lem:qgeoconvinf}
There exists a $\Ga$-equivariant embedding $\beta':\geo\Gamma \to\LasGa\subset\D_FX$
sending distinct ideal points to antipodal chambers, such that
for every $(L,A)$-quasigeodesic $q:\Z\to\Ga$ with ideal endpoints $\zeta_{\pm}$
we have

(i)
$\beta'(\zeta_{\pm})\subset\geo F(\beta(\zeta_-),\beta(\zeta_+))$.

(ii)
$q(m\pm n)x$ is contained in a tubular neighborhood 
of uniform radius $r'=r'(\Ga,L,A,x)$ 
of the euclidean Weyl chamber $V(q(m)x,\beta'(\zeta_{\pm}))$ 
for $m,n\in\N$. 
\end{lem}
\proof
Let $q:\Z\to\Ga$ be an $(L,A)$-quasigeodesic.
By Lemma~\ref{lem:qgeoflat}, 
$q(m\pm n)x$ is contained in the $r(\Ga,L,A,x)$-neighborhood of 
the maximal flat $F(\beta(\zeta_-),\beta(\zeta_+))$
and hence (by the triangle inequality) in the $2r(\Ga,L,A,x)$-neighborhood 
of the euclidean Weyl chamber $V(q(m)x,\si(m,\pm n))$ 
for some chamber
$\si(m,\pm n)\subset\geo F(\beta(\zeta_-),\beta(\zeta_+))$. 

The regularity of $\Ga$
implies that for every  $D>0$ we have
\begin{equation}
\label{ineq:asyreg}
d(d_{\De}(x,\ga x),\D\De)\geq D
\end{equation}
for all $\ga\in\Ga$ with $d_{\Ga}(\ga,e)\geq R=R(\Ga,d_{\Ga},x,D)$. Here
$d_\Ga$ denotes a word metric on $\Ga$.

It follows 
that $\si(m,\pm n)$ stabilizes as $n\to+\infty$ 
independently of $m$,
i.e.\ $\si(m,\pm n)=\si(m,\pm\infty)$
for $n\geq n(\Ga,L,A,x)$.  
In particular, we have chamber convergence $q(m\pm n)x\tauto\si(m,\pm\infty)$.
Since any two asymptotic quasirays in $\Ga$ have finite Hausdorff distance from each other,
the chamber limit $\si(m,\pm\infty)$ depends only on $\zeta_{\pm}$.
Putting $\beta'(\zeta_{\pm})=\si(m,\pm\infty)$,
we thus obtain a well-defined map $\beta':\geo\Ga\to\LasGa\subset\D_FX$
satisfying properties (i) and (ii).

The equivariance of $\beta'$ is clear from the construction. 
To verify its continuity, we argue by contradiction. 
Suppose that $\zeta_k\to\zeta$ in $\geo\Ga$, but $\beta'(\zeta_k)\to\si\neq\beta'(\zeta)$ in $\D_FX$. 
Since $\Ga$ is a word hyperbolic group, there exist
uniform quasigeodesics $q_k:\Z\to\Ga$ with $q_k(0)=1_{\Ga}$ and $q_k(+\infty)=\zeta_k$. 
After passing to a subsequence, we may assume that they converge (pointwise) 
to a quasigeodesic $q:\Z\to\Ga$ with $q(0)=1_{\Ga}$ and $q(+\infty)=\zeta$. 
Then there exists a sequence of natural numbers $n_k\to+\infty$ 
such that the points $q_k(n_k)$ are contained in a tubular neighborhood of the quasiray $q(\N)$,
i.e.\ $q_k(n_k)\to\zeta$ conically. 
Using property (ii), 
it follows that both the sequence of chambers $(\beta'(\zeta_k))$ and the constant sequence $(\beta'(\zeta))$
are shadow sequences in $\D_FX$ for the sequence $(q_k(n_k)x)$ of points in $X$.
The asymptotic uniqueness of shadows (Lemma~\ref{lem:shadasuniq})
implies that $\beta'(\zeta_k)\to\beta'(\zeta)$, contradicting our assumption. 
We conclude that $\beta'$ is continuous. 

It remains to verify that $\beta'$ is antipodal.
If $\zeta_{\pm}\in\geo\Ga$ are distinct ideal points,
then there exists a quasigeodesic $q$ asymptotic to them, $q(\pm\infty)=\zeta_{\pm}$.
For the nearest point projection $\bar q$ of $qx$ to 
$F(\beta(\zeta_-),\beta(\zeta_+))$ 
we have for large $n$ that 
$\bar q(\pm n)\in V(\bar q(\mp n),\beta'(\zeta_{\pm}))$
and the segment $\bar q(-n)\bar q(n)$ is regular.
This implies that the chambers $\beta'(\zeta_{\pm})$ are opposite to each other.
\qed

\medskip
Using the information on quasirays,
we can now show that the image of the modified boundary map $\beta'$ 
fills out the chamber limit set.
We conclude that the weaker asymptotic condition of boundary embeddedness already implies the stronger ones in the regular case:
\begin{prop}[Boundary embedded regular implies asymptotically embedded]
\label{prop:homeochlim}
Every $\sigma_{mod}$-boundary embedded discrete subgroup $\Ga\subset G$
is $\sigma_{mod}$-asymptotically embedded. 
\end{prop}
\proof
Lemma~\ref{lem:qgeoconvinf} yields that 
$\beta'(\geo\Ga)\subseteq\La_{\si_{mod}}(\Ga)$.
It suffices to prove that $\La_{\si_{mod}}(\Ga)=\beta'(\geo\Ga)$. 

The argument is similar to the proof of the continuity of $\beta'$ in Lemma~\ref{lem:qgeoconvinf}.
Let $\si\in\La_{\si_{mod}}(\Ga)$ 
and let $(\ga_n)$ be a sequence in $\Ga$ 
chamber-converging to $\si$, 
$\ga_n\tauto\si$. 
Since $\Ga$ is word hyperbolic, 
there exists a sequence of uniform quasigeodesics 
$q_n:\Z\to\Ga$ such that $q_n(0)=1_{\Ga}$ 
and $\ga_n\in q_n(\N)$. 
Let $\zeta_n\in\geo\Ga$ denote their forward ideal endpoints. 
According to Lemma~\ref{lem:qgeoconvinf}, 
the distance from $\ga_nx$ 
to the euclidean Weyl chamber $V(x,\beta'(\zeta_n))$ is uniformly bounded. 
Hence $(\beta'(\zeta_n))$ is a shadow of $(\ga_nx)$ in $\D_FX$,
and Lemma~\ref{lem:shadasuniq} implies that 
$\beta'(\zeta_n)\to\si$. 
Thus $\si\in\beta'(\geo\Ga)$. 
\qed

\begin{rem}
\label{rem:bdmp}
(i) 
Since both embeddings $\beta$ and $\beta'$ 
are continuous and $\Ga$-equivariant, 
the relative position (see \cite{coco})
$\pos(\beta',\beta):\geo\Ga\to W$ 
is continuous (locally constant) and $\Ga$-periodic.
For nonelementary hyperbolic groups $\Gamma$ this map must be constant,
because the action of a nonelementary word hyperbolic group on its Gromov boundary
is minimal. 

(ii)
One can show that if $\Ga$ is Zariski dense in $G$ 
then $\beta=\beta'$, cf. \cite{GW}. 

(iii)
On the other hand, in general, 
$\beta$ and $\beta'$ can be different 
as Example \ref{ex:alnebe} shows. 
\end{rem}

We take up again the discussion of the coarse geometry of the orbit map.
Elaborating on part (ii) of Lemma~\ref{lem:qgeoconvinf},
we will see next that the images of quasirays, since they must stay close to euclidean Weyl chambers,
they are forced to move out to infinity at a linear rate.

\begin{dfn}[$\Theta$-regular quasigeodesic]
Let $\Theta\subset \si_{mod}$ be a compact subset. 
A discrete quasigeodesic $p:I\cap \Z\to X$ is {\em $(s,\Theta)$-regular} if for every $m,n\in I\cap \Z$ with $|m-n|\geq s$
the segment $p(m)p(n)$ is $\Theta$-regular.
\end{dfn}

\begin{lem}
\label{lem:qgeotoqgeo}
For every $(L,A)$-quasigeodesic $q:\Z\to\Gamma$, its image
$qx$ in $X$ is an $(s,\Theta)$-regular discrete quasigeodesic
with $s$, a compact subset $\Theta\subset int(\sigma_{mod})$ and quasi-isometry constants depending  on $\Ga,L,A,x$. 
\end{lem}
\proof  
Since $\Ga$ is a discrete subgroup of $G$, 
the distance between orbit points can be bounded from below 
in terms of the word metric, 
i.e.\ there is an estimate of the form 
$$d(\ga x,\ga' x)\geq f_x(d_{\Ga}(\ga,\ga'))$$ 
with $f_x(t)\to+\infty$ as $t\to+\infty$. 
As a consequence, 
for $\rho>0$ 
we have $d(q(m)x,q(m+n)x)\geq\rho$ 
for $n\geq n(f_x,L,A,\rho)$. 

As before, let $\zeta_{\pm}\in\geo\Ga$ denote the ideal endpoints of $q$.
We consider the 
nearest point projection $\bar q$ of $qx$ to 
$F(\beta(\zeta_{-}), \beta(\zeta_{+}))$. 
Choosing 
$\rho \gg r(\Ga,L,A,x),r'(\Ga,L,A,x)$ 
and invoking again the regularity of $\Ga$, 
cf.\ \eqref{ineq:asyreg}, 
we obtain 
as in the end of the proof of Lemma~\ref{lem:qgeoconvinf} 
that 
\begin{equation}
\label{eq:projinwhu}
\bar q(m\pm n)\in V(\bar q(m),\beta'(\zeta_{\pm}))
\end{equation}
for $n\geq n'(\Ga,L,A,x)$. 
It follows that along the coarsening $\bar q|_{n'\Z}$ of $\bar q$, 
the $\De$-distances between its points 
are additive in the sense that 
\begin{equation}
\label{eq:distadd}
d_{\De}(\bar q(m_1),\bar q(m_2))+d_{\De}(\bar q(m_2),\bar q(m_3))
=d_{\De}(\bar q(m_1),\bar q(m_3))
\end{equation}
if $m_3-m_2,m_2-m_1\geq n'$. 
In particular,
$\bar q$ and hence $qx$ is a $\Theta$-regular uniform quasigeodesic. 
\qed

We summarize the properties of asymptotically embedded subgroups $\Ga$ 
established in the previous two lemmas:

(i) 
For some (every) $x\in X$ there exists a constant $r$ 
such that for every discrete geodesic ray $q: \N\to \Ga$, 
the points $q(n)x$ belong to the $r$-neighborhood 
of the Weyl cone $V(q(0)(x), \st(\tau))$ 
for some $\tau\in\Flagt$ depending on $q$.

(ii) 
For all sequences $\ga_n\to\infty$ in $\Ga$,
all subsequential limits of the quantities  
$$
\frac{1}{|\ga_n|}d_\Delta(x, \ga_n x)
$$
are in the interior of the cone $V(0,\st(\tau_{mod}))$, 
where $|\ga_n|$ denotes the word length of $\ga_n\in\Ga$.

\medskip
As a consequence of Lemmata \ref{lem:qgeoconvinf}, \ref{lem:qgeotoqgeo}, we obtain:
 
\begin{thm}[Coarse geometric properties of boundary embedded regular subgroups]
\label{thm:qiemb}
Let $\Ga\subset G$ be 
a  $\sigma_{mod}$-boundary embedded discrete subgroup.
Then $\Ga$ is 
$\sigma_{mod}$-asymp\-to\-ti\-cal\-ly embedded, 
uniformly regular and the orbit maps $\Ga\to\Ga x\subset X$ 
are quasi-isometric embeddings.
\end{thm}
\proof
Since $\Ga$ is word hyperbolic, 
through any two elements $\ga,\ga'\in\Ga$ there
exists a complete quasigeodesic $q:\Z\to\Ga$ 
with uniform quasi-isometry constants $L,A$. 
The assertion then follows from Lemma~\ref{lem:qgeotoqgeo}. 
\qed

\subsubsection{Asymptotically embedded groups}
\label{sec:coarse_asymptotically_embedded}

In this section we prove a version of the results of the previous section for $\tau_{mod}$-regular  subgroups, under the stronger assumption of  asymptotic embeddedness.

Let $\Gamma\subset G$ be $\tau_{mod}$-asymptotically embedded,
and let $\alpha:\geo \Gamma\to\Lambda_{\tau_{mod}}(\Gamma)$ be the equivariant homeomorphism. 
Lemma~\ref{lem:qgeoflat} holds for $\alpha$. We will now obtain more specific information on the position
of the image $qx$ of the quasigeodesic $q$ under the orbit map,
generalizing Lemma~\ref{lem:qgeoconvinf}(ii).

\begin{lem}
\label{lem:qgeoconvinfgeneral}
For every $(L,A)$-quasigeodesic $q:\Z\to\Ga$ with ideal endpoints $\zeta_{\pm}$,
the point $q(m\pm n)x$ is contained in a tubular neighborhood 
of uniform radius $r''=r''(\Ga,L,A,x)$ 
of the Weyl cone $V(q(m)x, \st(\alpha(\zeta_{\pm})))$ 
for $m,n\in\N$. 
\end{lem}

\proof
Due to the continuity at infinity of orbit maps (Lemma~\ref{lem:continuity_at_infinity})
we have that $q(n)\tauto \tau_+$ as $n\to+\infty$,
where we abbreviate $\tau_\pm=\alpha(\zeta_\pm)$.

According to Lemma~\ref{lem:qgeoflat},
$qx$ is contained in a uniform tubular neighborhood of the parallel set $P(\tau_-,\tau_+)$. 
Therefore, the $\tau_{mod}$-regular sequence $(q(n))$ in $\Ga$ has a shadow sequence $(\si_n)$ in $\D_FX$ 
consisting of chambers $\si_n\subset\geo P(\tau_-,\tau_+)$. 
Taking their type $\tau_{mod}$ faces, one obtains a shadow sequence $(\tau_n)$ in $\Flagt$ 
consisting of simplices $\tau_n\subset\geo P(\tau_-,\tau_+)$. 
More precisely,
$q(n)x$ has uniformly bounded distance from $V(q(0)x,\si_n)\subset V(q(0)x,\st(\tau_n))$.
The asymptotic uniqueness of shadows (Lemma~\ref{lem:shadasuniq})
implies that $\tau_n\to \tau_+ $.

We use now that $\tau_+$ is
isolated among the type $\tau_{mod}$ simplices occurring in  
$\geo P(\tau_-, \tau_+)$, see Lemma \ref{lem:poleisol}(i). 
It follows that $\tau_n=\tau_+$ for sufficiently large $n$.
This means that the sequence $(q(n)x)$ enters a uniform tubular neighborhood of the Weyl cone 
$V(q(0)x,\st(\tau_+))$.

It remains to show that the entry time is uniform. 
We will do this by backtracking,
based on the fact that $V(q(0)x,\ost(\tau_+))$ is an open subset of $P(\tau_-,\tau_+)$.
The latter follows from the fact that $\ost(\tau_+)$ is an open subset of $\geo P(\tau_-,\tau_+)$
with respect to the visual topology, 
see Lemma~\ref{lem:poleisol}(ii). 

Let $\bar p_n$ denote the nearest point projection of  $q(n)x$ to the parallel set $P(\tau_-, \tau_+)$.
Since $\Gamma$ is $\tau_{mod}$-regular,    the segment $\bar p_0\bar p_n$ is $\tau_{mod}$-regular
for $n\geq n(\Ga,L,A,x)$.
More precisely, given $D>0$, the distance of $\bar p_n$ 
from $\partial V( \bar p_0, \st( \tau_n))= V( \bar p_0, \partial \st( \tau_n))$
is  at least $ D$ provided that $n\geq n_0(D)$, with $n_0(D)$ independent of $q$.
We choose $D$ sufficiently large depending on $L,A,x$ and $r$ (the constant $r$ from Lemma~\ref{lem:qgeoflat}) 
so that $D>>r$ and  $d(\bar p_n,\bar p_{n+1})< D$ for all $n$ and $q$.
Hence the sequence $\bar p_n$ cannot enter the cone $V(\bar p_0,\st(\tau_+))$ after time $n_0(D)$,
i.e.\ if $\bar p_{n+1}$ belongs to the cone for $n\geq n_0(D)$, then $\bar p_n$ itself belongs to the cone.
Therefore
$\bar p_n\in V(\bar p_0,\st(\tau_+))$ for $n\geq n_0(D)$.

The assertion follows by suitably enlarging the uniform radius of the tubular neighborhood, e.g.\ by $n_0(D)\cdot D$.
\qed

\medskip
The following is an analogue of Lemma~\ref{lem:qgeotoqgeo}.

\begin{lem}
\label{lem:qgeotoqgeotau}
For every $(L,A)$-quasigeodesic $q:\Z\to\Gamma$, its image
$qx$ in $X$ is an $(s,\Theta)$-regular discrete quasigeodesic
with $s$, a compact subset $\Theta\subset \ost(\tau_{mod})$ and the quasi-isometry constants depending  on $\Ga,L,A,x$. 
\end{lem}

\proof
We continue the argument in the previous proof, keeping the notation.

There exists a compact Weyl convex subset $\Theta\subset \ost(\tau_{mod})$
depending on $\Ga,L,A,x$ such that the segments $\bar p_0 \bar p_{n_0}$ are  $\Theta$-regular for all $q$. 
This is because there are only finitely many elements $q(0)^{-1}q(n_0)\in\Ga$.
Moreover, there is a lower bound $d(\bar p_0,\bar p_{n_0})\geq d_0=d_0(\Ga,L,A,x)>0$.

We consider the nearest point projection $V(\bar p_0,\st(\tau_+))\to V(\bar p_0,\tau_+)$, which can be regarded 
as the restriction of the projection from $P(\tau_-,\tau_+)$ to its Euclidean de Rham factor.
Let  $\bar{\bar p}_n\in V(\bar p_0,\tau_+)$ denote the projection of $\bar p_n$.
Since ${\bar p}_{m+n_0}\in V({\bar p}_m,\st(\tau_+))$, it follows that
$\bar{\bar p}_{m+n_0}\in V(\bar{\bar p}_m,\tau_+)$. 
Since $\Theta\subset \interior ( B(\xi,\pihalf))$ for any $\xi\in \interior (\tau_{mod})$ (see Lemma \ref{lem:half}), 
we have that
$d(\bar{\bar p}_m,\bar{\bar p}_{m+n_0})\geq c(\Theta) d({\bar p}_m,{\bar p}_{m+n_0}) $ 
with a constant $c(\Theta)>0$.
Inductively, we obtain that 
$d(\bar{\bar p}_m,\bar{\bar p}_{m+k n_0})\geq k c(\tau_{mod} )c(\Theta) d_0$
for $k\in\N$ and some constant $c(\tau_{mod})>0$,
compare the proof of Lemma~\ref{lem:qgeotoqgeo}.
This establishes that $q x$ is a uniform quasigeodesic.

The inclusion ${\bar p}_{m+n_0}\in V({\bar p}_m,\st(\tau_+))$ and the $\Theta$-regularity of $\bar p_m \bar p_{m+n_0}$ imply
that  ${\bar p}_{m+n_0}\in V({\bar p}_m,\st_\Theta(\tau_+))$. Induction and the convexity of $\Theta$-cones 
(Proposition~\ref{prop:wconeconv})
yield
${\bar p}_{m+ k n_0}\in V({\bar p}_m,\st_\Theta(\tau_+))$ for all $0<k\in\N$. 
After slightly enlarging  $\Theta$ and choosing $s$ sufficiently large (both depending on $\Ga,L,A,x$) we obtain that $qx$ is $(s,\Theta)$-regular.
\qed

\medskip
As a consequence we obtain, analogously to Theorem~\ref{thm:qiemb}:
\begin{thm}[Coarse geometric properties of asymptotically embedded subgroups]
\label{thm:qiembrtaureg}
Let $\Ga\subset G$ be a $\tau_{mod}$-asymptotically embedded discrete subgroup 
with boundary embedding $\alpha$.
Suppose furthermore that $\Ga$ is non-elementary as a word hyperbolic group, $|\geo\Ga|\geq3$.
Then:

(i)
$\Ga$ is uniformly $\tau_{mod}$-regular.

(ii)
The orbit maps $\Ga\to\Ga x\subset X$ are quasi-isometric embeddings.

(iii)
The action $\Gamma\acts X$ is strongly conical in the following sense: 
For every $\zeta\in\geo\Gamma$ and quasiray $q: \N\to \Gamma$ asymptotic to $\zeta$,
the image quasiray $qx$ lies in a uniform tubular neighborhood of the Weyl cone $V(q(0)x,\st (\alpha(\zeta)))$.
\end{thm}

The following example shows that being undistorted 
without further restrictions is a very weak concept in higher rank, unlike in rank one.
 
\begin{example}
[Infinitely presented undistorted subgroups]
\label{ex:stallings}
Consider the group $F_2\times F_2$ where $F_2$ is the free group of rank $2$. 
Let $\phi: F_2\to \Z$ be the homomorphism which sends both free generators of $F_2$ to the generator of $\Z$. Let 
$\Ga\triangleleft F_2\times F_2$ denote the normal subgroup
$$
\Ga=\{(h_1, h_2): \phi(h_1)=\phi(h_2)\}. 
$$
Then $\Ga$ is finitely generated but $\Ga$ 
is not finitely presentable (see \cite{Baumslag-Roseblade}). We claim that 
the subgroup $\Ga$ is undistorted  in $F_2\times F_2$. Indeed, let $w$ be a path in the Cayley graph of $F_2\times F_2$ connecting the unit element $e$ to an element 
$\ga\in \Ga$. We will equip $\Z$ with the presentation which is the quotient of the standard presentation of $F_2\times F_2$. Then $\phi(w)$ is a loop in $\Z$, which, therefore, bounds a van Kampen diagram $D$ in the Cayley complex of $\Z$. Since $\Z$ has 
a linear Dehn function, the diagram $D$ has combinatorial area at most $C\ell(w)$, where $\ell(w)$ is the length of 
$w$ (which is the same as the length of $\phi(w)$). Lifting $D$ to the Cayley complex of $F_2\times F_2$ results in a van Kampen diagram $\tilde{D}$ 
bounding a bigon one of whose sides is $w$ and the other is a path $u$ in the Cayley graph of $\Ga$. 
Since    $\tilde{D}$ has combinatorial area at most a constant times the combinatorial area of  $D$, we conclude that 
$$
\ell(u)\le C' \ell(w).  
$$
Thus, $\Ga$ is indeed undistorted in $F_2\times F_2$. 

Realizing $F_2$ as a convex cocompact subgroup of $\Isom(\H^2)$, we obtain a discrete quasi-isometric embedding $F_2\times F_2\to \Isom(X), X=\H^2\times \H^2$. Then the subgroup $\Ga\subset\Isom(X)$ 
is undistorted and not finitely presentable. On the other hand, since $\Ga$ is not finitely presented, there is no coarse Lipschitz retraction $X\to\Ga x$. 

Note that the group $\Ga$ in this example is not weakly regular.
\end{example}

Theorem~\ref{thm:qiembrtaureg}, in particular part 3, 
can be regarded as a higher rank version of the Morse Lemma for quasigeodesics in hyperbolic spaces.
We will study in section~\ref{sec:morse} quasigeodesics with such a Morse property and Morse actions,
whose orbit maps send uniform quasigeodesics to Morse quasigeodesics. We will show that the class of
Morse actions coincides with the class of actions having the strong asymptotic properties discussed in this section
(weakly RCA, asymptotically embedded and expanding at infinity).

\subsection{The Anosov condition}

\subsubsection{Anosov representations}\label{sec:anosov}

A notion of Anosov representations of surface groups into $PSL(n,\R)$ was introduced by Labourie in \cite{Labourie},
and generalized to a notion of $(P_+,P_-)$-Anosov representations $\Ga\to G$ of word hyperbolic groups into semisimple Lie groups by 
Guichard and Wienhard in \cite{GW}. The goal of this section is to review this definition of Anosov representations $\Ga\to G$ using the language of {\em expanding and contracting flows} and then present a closely related and equivalent definition which avoids the language of flows. 

Let $\Gamma$ be a non-elementary (i.e.\ not virtually cyclic) word hyperbolic group with a fixed word metric $d_\Ga$ and Cayley graph $C_\Ga$.  Consider a {\em geodesic flow} $\widehat\Ga$ of $\Ga$;  such a  flow was originally constructed by Gromov \cite{Gromov} and then improved by Champetier \cite{Champetier} and Mineyev \cite{Mineyev}, resulting in definitions with different properties. 
We note that the exponential convergence of asymptotic geodesic rays will not be used in our discussion; 
as we will see, it is also irrelevant whether the trajectories of the geodesic flow are geodesics or uniform quasigeodesics in $\widehat\Ga$. In particular, it will be irrelevant for us which definition of $\widehat\Ga$ is used. 
Only the following properties of $\widehat\Ga$ will be used in the sequel:

1. $\widehat\Ga$ is a proper metric space. 

2. There exists a properly discontinuous isometric action 
$\Ga\acts\widehat\Ga$. 

3. There exists a $\Ga$-equivariant quasi-isometry $\pi: \widehat\Ga\to \Ga$; in particular, the fibers of $\pi$ are relatively compact.

4. There exists a continuous action $\R\acts\widehat\Ga$, denoted $\phi_t$ and called the {\em geodesic flow}, 
whose trajectories are uniform quasigeodesics in $\widehat\Ga$, i.e. for each $\hat m\in \widehat\Ga$ the 
{\em flow line} 
$$
t\to \hat m_t := \phi_t(\hat m)
$$ 
is a uniform quasi-isometric embedding $\R\to \widehat\Ga$. 

5. The flow $\phi_t$ commutes with the action of $\Ga$. 

6. Each $\hat m\in\hat\Ga$ defines a uniform  
quasigeodesic $m:t\mapsto m_t$ in $\Ga$ by the formula: 
$$
m_t= \pi(\hat m_t) 
$$
The natural map 
$$
e=(e_-,e_+):\hat\Ga\to\geo\Ga\times\geo\Ga-\Diag
$$
assigning to $\hat m$ the pair of ideal endpoints $(m_{-\infty},m_{+\infty})$ 
of $m$ is continuous and surjective. 
In particular,
every uniform quasigeodesic in $\hat\Ga$ 
is uniformly Hausdorff close to a flow line.

\medskip 
The reader can think of the elements of $\widehat\Ga$ {\em as parameterized geodesics in $C_\Ga$, so that $\phi_t$ acts on geodesics via reparameterization}. This was Gromov's original viewpoint, 
although not the one in \cite{Mineyev}.  

We say that $\hat m\in \widehat\Ga$ is {\em normalized} if $\pi(\hat m)=1\in\Ga$. 
Similarly, maps $q: \Z\to \Ga$, and $q: \N\to \Ga$ will be called {\em normalized} if $q(0)=1$. 
It is clear that every $\hat m\in \widehat\Ga$ can be sent to a normalized element of $\widehat\Ga$ 
via the action of $m_0^{-1}\in \Ga$.

Since trajectories of $\phi_t$ are uniform quasigeodesics, for 
each normalized $\hat m\in \widehat\Ga$ we have
\begin{equation}\label{eq:distance}
d_\Ga(1, m_t)  \approx t
\end{equation}
in the sense that 
$$
C_1^{-1} t - C_2\le d_\Ga(1, m_t) \le C_1 t + C_2
$$
for some positive constants $C_1, C_2$. 

Let ${\mathrm F}^\pm=\Flag(\tau_{mod}^{\pm})$ be a pair of {\em opposite} partial flag manifolds 
associated to the Lie group $G$, i.e.\ they are quotient manifolds of the form ${\mathrm F}^\pm=G/P_\pm$, 
where $P_{\pm}$ are opposite (up to conjugation) parabolic subgroups. 
The conjugacy classes of $P_{\pm}$ correspond to 
faces $\tau_{mod}^\pm$ of the model spherical Weyl chamber $\si_{mod}$
related by $\iota(\tau_{mod}^\pm)=\tau_{mod}^\mp$.
As usual, we will regard elements of ${\mathrm F}^\pm$ as simplices of type 
$\tau_{mod}^\pm$ in the Tits boundary of $X$.

Define the trivial bundles  
$$
E^\pm =\widehat \Ga \times {\mathrm F}^\pm \to \widehat \Ga. 
$$
For every representation $\rho: \Ga\to G$, the group $\Ga$ acts on both bundles via its natural action on $\widehat\Ga$ and via the representation $\rho$ on ${\mathrm F}^\pm$. 
Put a $\Ga$-invariant background Riemannian metric on the fibers of theses bundles, which varies  continuously with respect to $\hat m\in \widehat\Ga$. We will use the notation ${\mathrm F}^\pm_{\hat m}$ for the fiber above the point $\hat m$ equipped with this Riemannian metric. Since the subspace of $\widehat\Ga$ consisting of normalized elements is 
compact, it follows that for normalized $\hat m, \hat m'$ the identity map
$$
{\mathrm F}^\pm_{\hat m} \to {\mathrm F}^\pm_{\hat m'}
$$
is uniformly bilipschitz (with bilipschitz constant independent of $\hat m, \hat m'$). 
We will identify $\Ga$-equivariant (continuous) sections of the bundles $E^\pm$ with equivariant 
maps $s_\pm: \widehat\Ga\to {\mathrm F}^\pm$. 
These sections are said to be 
{\em parallel along flow lines}
if 
$$
s_\pm(\hat m)= s_\pm(\hat m_t)
$$
for all $t\in \R$ and $\hat{m} \in \widehat\Ga$. 

\begin{defn}
Parallel sections $s_\pm$ are called 
{\em strongly parallel along flow lines} if for any two flow lines 
$\hat m,\hat m'$ with the same ideal endpoints, we have
$s_\pm(\hat m)=s_\pm(\hat m')$.
\end{defn}

Note that this property is automatic 
for the geodesic flows constructed by Champetier and Mineyev
since (for their flows) 
any two flow lines which are at finite distance from each other 
are actually equal. 
Strongly parallel sections define $\Ga$-equivariant {\em boundary maps}
\begin{equation*}
\beta_\pm:\geo\Ga\to {\mathrm F}^\pm
\end{equation*}
from the Gromov boundary $\geo \Ga$ of the word hyperbolic group $\Ga$ by:
\begin{equation}
\label{eq:bdmapsect}
\beta_{\pm}\circ e_{\pm} = s_{\pm}
\end{equation}
\begin{lem}
$\beta_{\pm}$ is continuous.
\end{lem}
\proof
Let $(\xi^n_-,\xi^n_+)\to(\xi_-,\xi_+)$
be a converging sequence in $\geo\Ga\times\geo\Ga-\Diag$.
There exists a bounded sequence $(\hat m^n)$ in $\hat\Ga$ 
so that $e_{\pm}(\hat m^n)=\xi^n_{\pm}$.
It subconverges, $\hat m^n\to\hat m$.
The continuity of $s_{\pm}$ implies that 
$\beta_{\pm}(\xi^n_{\pm})=s_{\pm}(\hat m^n)\to s_{\pm}(\hat m)=\beta_{\pm}(\xi_{\pm})$.
This shows that no subsequence of $(\beta_{\pm}(\xi^n_{\pm}))$ can have a limit $\neq\beta_{\pm}(\xi_{\pm})$,
and the assertion follows because ${\mathrm F}^{\pm}$ is compact.
\qed

\medskip
Conversely, equivariant continuous maps $\beta_\pm$
define $\Ga$-equivariant sections strongly parallel along flow lines 
just by (\ref{eq:bdmapsect}).

Consider the canonical ``identity" maps
$$
\Phi_{\hat m,t}: {\mathrm F}^\pm_{\hat m} \to {\mathrm F}^\pm_{\phi_t\hat m}. 
$$
These maps distort the Riemannian metric on the fibers. 
Using Definition~\ref{def:expansion factor},
we define the {\em expansion factor} of the flow $\phi(t)$ on the fiber  ${\mathrm F}^\pm_{\hat m}$ at the point $s_\pm(\hat m)$ as
$$
\eps_\pm(\hat m, t) :=
\eps(\Phi_{\hat m,t}, s_\pm(\hat m)),
$$
see Definition \ref{def:expansion factor} for the definition of the expansion factor of a diffeomorphism. 

\begin{definition}\label{defn:uniformly expanding}
The geodesic flow $\phi_t$ is said to be {\em uniformly exponentially expanding} on the bundles $E^\pm$ with respect to the sections $s_\pm$ if there exist constants $a, c>0$ such that
$$
\eps_{\pm}(\hat m,\pm t) \ge a e^{ct}
$$
for all $\hat m\in \widehat \Ga$ and $t\geq0$.
\end{definition}

Our next goal is to give an alternative interpretation for the uniform expansion in this definition. 
First of all, since the metrics on the fibers are $\Ga$-invariant, it suffices to verify uniform exponential expansion only for normalized elements of $\widehat\Ga$. Consider a normalized element $\hat m\in \widehat\Ga$ and for $t\in\R$ 
the composition
$$
m_t^{-1} \circ \Phi_{\hat m,t}:  {\mathrm F}^\pm_{\hat m} \to {\mathrm F}^\pm_{m_t^{-1}\hat m_t}.
$$
Note that $\pi(m_t^{-1}\hat m_t)=m_t^{-1}m_t=1$,
i.e.\ both $\hat m$ and $m_t^{-1}\hat m_t$ are normalized.  
Since the group $\Ga$ acts isometrically on the fibers of the bundles $E^\pm$, the metric distortion of the above compositions is exactly the same as the distortion of $\Phi_{\hat m,t}$. Furthermore, since, as we noted above, the metrics on 
${\mathrm F}^\pm_{\hat m}$ and ${\mathrm F}^\pm_{m_t^{-1}\hat m_t}$ 
are uniformly bilipschitz to each other (via the ``identity'' map), the rate of expansion for the above composition is (up to a uniform multiplicative error) the same as the expansion rate for the map 
$$
\rho(m_t^{-1}): {\mathrm F}^\pm\to {\mathrm F}^\pm .
$$
(Here we are using fixed background Riemannian metrics on ${\mathrm F}^\pm$.) Thus, we get the estimate 
$$
C_3^{-1} \eps( \rho(m_t^{-1}), \beta_{\pm}(m_{\pm\infty}))\le    \eps_\pm(\hat m, t) \le 
C_3 \eps( \rho(m_t^{-1}), \beta_{\pm}(m_{\pm\infty})) 
$$
for some uniform constant $C_3>1$. By taking into account the equation \eqref{eq:distance}, we  obtain the following equivalent reformulation of Definition \ref{defn:uniformly expanding}:

\begin{lemma}\label{lem:expansion}
The geodesic flow is uniformly exponentially expanding with respect to the sections $s_\pm$ 
if and only if for every normalized uniform quasigeodesic 
$\ga: \Z\to \Ga$, which is asymptotic to points $\xi_\pm=\ga(\pm \infty)\in \geo \Ga$, the elements $\rho(\ga(\pm n))^{-1}$ act on $T_{\beta_\pm(\xi_\pm)} {\mathrm F}^\pm$ with uniform exponential expansion rate, i.e.
$$
\eps( \rho(\ga(\pm n))^{-1}, \beta_\pm(\xi_\pm))\ge Ae^{Cn}
$$
for all $\hat m\in \widehat \Ga$ and $n\ge 0$ 
with some fixed constants $A,C>0$.
\end{lemma}
\proof
There exists a normalized flow line $\hat m$ uniformly close to $\ga$, 
i.e.\ $\ga(n)$ is uniformly close to $m_{t_n}$ with $t_n\approx n$ uniformly. 
Then $m_{\pm\infty}=\xi_\pm$, and
$\eps(\rho(\ga(\pm n))^{-1}, \beta_\pm(\xi_\pm))$
equals $\eps( \rho(m_{t_{\pm n}}^{-1}), \beta_{\pm}(m_{\pm\infty}))$
up to a uniform multiplicative error,
and hence also 
$\eps_\pm(\hat m, t_{\pm n})$.
\qed

\medskip
Since every uniformly quasigeodesic ray $\ga: \N\to \Ga$ extends to a uniform complete quasigeodesic  $\ga: \Z\to \Ga$,  and in view of Morse lemma for hyperbolic groups, in the above definition it suffices to consider only normalized integer geodesic rays $\ga: \N\to \Ga$.   

We can now  give the original and an alternative definition of Anosov representations. 
\begin{definition}
A pair of continuous maps $\beta_\pm: \geo \Ga \to {\mathrm F}^\pm$ is said to be {\em antipodal} if it satisfies the following conditions (called {\em compatibility} in \cite{GW}): 

(i) 
For every pair of distinct ideal points $\zeta, \zeta'\in \geo \Ga$, the simplices $\beta_+(\zeta)$, $\beta_-(\zeta')$ in the Tits boundary of $X$ are antipodal, equivalently, the corresponding parabolic subgroups of $G$ are opposite.  (In \cite{GW} this property is called {\em transversality}.) 

(ii) 
For every $\zeta\in \geo \Ga$, the simplices $\beta_+(\zeta), \beta_-(\zeta)$ belong to the same spherical Weyl chamber, i.e.\ the intersection of the corresponding parabolic subgroups of $G$ contains a Borel subgroup.
\end{definition}

Note that, as a consequence, the maps $\beta_{\pm}$ are {\em embeddings}, 
because antipodal simplices cannot be faces of the same chamber. 

\begin{definition}[\cite{GW}]
A representation $\rho: \Ga\to G$ is said to be $(P_+, P_-)$-Anosov if there exists an antipodal pair of continuous $\rho$-equivariant 
maps $\beta_\pm: \geo \Ga \to {\mathrm F}^\pm$ such that the geodesic flow on the associated bundles $E^\pm$ satisfies the uniform expansion property with respect to the sections $s_\pm$  associated to the maps $\beta_\pm$. 
\end{definition}

The pair of maps $(\beta_+,\beta_-)$ in this definition is called {\em compatible} with the Anosov representation $\rho$.
Note that 
a $(P_+,P_-)$-Anosov representation admits a unique compatible pair of maps. 
Indeed, 
the fixed points of infinite order elements $\ga\in\Ga$ 
are dense in $\geo\Ga$. 
The maps $\beta_{\pm}$ send the attractive and repulsive fixed points of $\ga$
to fixed points of $\rho(\ga)$ with contracting and expanding 
differentials,
and these fixed points are unique.
In particular, if $P_+$ is conjugate to $P_-$ 
then $\beta_-=\beta_+$. 

We note that Guichard and Wienhard in  \cite{GW} use in their definition the uniform contraction property of the reverse flow $\phi_{-t}$ instead of the expansion property used above, but the two are clearly equivalent. Note also that  in the definition, it suffices to verify the uniform exponential expansion property only for the bundle $E_+$.
We thus obtain, as a corollary of Lemma \ref{lem:expansion}, the following alternative definition of Anosov representations:

\begin{prop}
[Alternative definition of Anosov representations] 
A representation $\rho: \Ga\to G$ is $(P_+, P_-)$-Anosov if and only if there exists a pair of antipodal continuous $\rho$-equivariant maps $\beta_\pm: \geo \Ga \to {\mathrm F}^\pm$ such that for every normalized geodesic ray 
(equivalently, for every uniformly quasigeodesic ray) $\ga: \N\to \Ga$ asymptotic to $\xi\in \geo \Ga$, the elements $\rho(\ga(n))^{-1}$ act on $T_{\beta_+(\xi)} {\mathrm F}_+$ with uniform exponential expansion rate, i.e.
\begin{equation}\label{eq:ano}
\eps(\rho(\ga(n))^{-1}, \beta_+(\xi))\ge A e^{Cn}
\end{equation}
for $n\geq0$
with  constants $A,C>0$ which are independent of $\ga$. 
\end{prop}

\medskip 
We now restrict to the case that the parabolic subgroups $P_{\pm}$ are conjugate to each other, 
i.e.\ the simplices $\tau^\pm_{mod}$ are equal to 
an $\iota$-invariant face $\tau_{mod}$ of $\si_{mod}$. The $(P_+, P_-)$-Anosov representations will in this case be called simply {\em $P$-Anosov}, where $P=P_+$, or $\tau_{mod}$-Anosov. 
Note that the study of general $(P_+,P_-)$-Anosov representations quickly reduces 
to the case of $P$-Anosov representations by intersecting parabolic subgroups,
cf.\ \cite[Lemma 3.18]{GW}.
Now, 
$$
{\mathrm F}^{\pm}={\mathrm F}=G/P=\Flagt
$$ 
and 
$$
\beta_{\pm}=\beta:\geo\Ga\to {\mathrm F} 
$$ 
is a single continuous embedding.
The compatibility condition 
reduces to the {\em antipodality} condition: 
For any two distinct ideal points $\zeta,\zeta'\in\geo \Ga$
the simplices $\beta(\zeta)$ and $\beta(\zeta')$ are antipodal to each other. 
In other words, $\beta$ is a {\em boundary embedding} in the sense of Definition~\ref{dfn:nicastau}.

We thus arrive to our definition:
\begin{defn}[Anosov representation]
\label{defn:our anosov}
Let $P$ be a parabolic subgroup which is conjugate to its opposite parabolic subgroups, 
and let $\tau_{mod}\subset\si_{mod}$ the corresponding face type. 
We call a representation $\rho: \Ga\to G$ {\em $P$-Anosov} or $\tau_{mod}$-Anosov
if it is $\tau_{mod}$-boundary embedded (cf.\ Definition~\ref{dfn:nicastau})
with boundary embedding $\beta: \geo \Ga \to {\mathrm F}=G/P$ 
such that for every normalized geodesic ray $q: \N\to \Ga$ asymptotic to $\zeta\in \geo \Ga$, 
the elements $\rho(q(n))^{-1}$ act on $T_{\beta(\zeta)} {\mathrm F}$ with uniform exponential expansion rate, i.e.\ 
$$
\eps(\rho(q(n))^{-1}, \beta(\zeta))\ge A e^{Cn}
$$ 
for $n\geq 0$ with constants $A, C>0$  independent of $q$. 
\end{defn}

We will refer to $\rho(\Ga)$ as a $\tau_{mod}$-Anosov subgroup of $G$.

\subsubsection{Non-uniformly expanding Anosov representations}\label{sec:weak Anosov}

In this section we discuss a further weakening of the Anosov condition which leads, however, to the same 
class of group actions. We restrict our discussion to the case of $\iota$-invariant model simplices $\tau_{mod}\subset\si_{mod}$ and the corresponding parabolic subgroups $P\subset G$ (conjugate to their opposites), even though, with minor modifications, the same  proofs go through for arbitrary pairs of opposite (up to conjugation) parabolic subgroups $P_\pm\subset G$.  

We note that in the definition of Anosov representation (both the original definition and the alternative one) the constants $a$ and $c$ (respectively, $A$ and $C$) were required to be uniform for the entire group. The main goal of this section is to show that the requirement of uniform exponential expansion can be relaxed and that the weakened notion is still equivalent to the Anosov condition as well as to the concept of asymptotic embedding for discrete subgroups.

\begin{defn}
[Non-uniformly Anosov representation] 
We call a representation $\rho: \Ga\to G$ of a word hyperbolic group $\Ga$ 
{\em non-uniformly $\tau_{mod}$-Anosov} 
if it is $\tau_{mod}$-boundary embedded 
with boundary embedding $\beta: \geo \Ga \to\Flagt$
such that for every normalized discrete geodesic ray (equivalently, normalized  uniform quasiray) 
$q: \N\to \Ga$ asymptotic to $\zeta\in \geo \Ga$, 
the elements $\rho(q(n))^{-1}$ act on $T_{\beta(\zeta)} {\mathrm F}$ with unbounded expansion rate: 
\begin{equation}
\label{eq:ano2}
\sup_{n\geq0}\, \eps(\rho(q(n))^{-1}, \beta(\zeta))=+\infty. 
\end{equation}
\end{defn}

Note that this definition does not even have the requirement that the expansion rate of $\rho(q(n))^{-1}$ at $\tau$ diverges to infinity. Other weakenings of the Anosov condition appear in \cite[sec.\ 6.1]{Labourie} and \cite[Prop.\ 3.16]{GW}.
They assume uniform, not necessarily exponential, divergence to infinity of
the expansion factors along the trajectories of the geodesic flow. 

Non-uniformly Anosov representations clearly have finite kernel (and discrete image).  
Therefore we will consider from now on only the case when $\Ga$ is a subgroup of $G$. 

\begin{thm}[Non-uniformly Anosov implies asymptotically embedded]
\label{thm:weak anosov->as embedded}
Every non-uniformly $\tau_{mod}$-Anosov subgroup 
is $\tau_{mod}$-asymptotically embedded.
\end{thm}
\proof We first establish a weak form of continuity at infinity and conicality for the boundary embedding 
of a non-uniformly Anosov subgroups:

\begin{lem}
\label{lem:weak anosov->weak as embedded}
Suppose that $\Ga\subset G$ is non-uniformly $\tau_{mod}$-Anosov 
with boundary embedding $\beta$.
Then for every discrete geodesic ray $q: \N\to \Ga$, 
the sequence $(q(n))$ in $\Ga$ contains a $\tau_{mod}$-regular subsequence $(q(n_i))$ such that 
$q(n_i)\tauto\beta(q(+\infty))$ conically.
In particular, 
$\beta(\geo\Ga)\subset\LatGa$.
\end{lem}
\proof We fix a point $x\in X$. 
Since by definition, $\Ga$ is $\tau_{mod}$-boundary embedded,
the discussion in section~\ref{sec:bdembgp} applies.
By Lemma \ref{lem:qgeoflat}, 
the image under the orbit map $\Ga\to \Ga x \subset X$ 
of every discrete geodesic 
$q:\Z\to\Ga$ is uniformly close to the parallel set 
$P(\beta(q(-\infty)),\beta(q(+\infty)))$, 
i.e.\ there is a constant $r>0$ independent of $q$ such that 
\begin{equation}\label{eq:bdd distance}
d(q(n),P(\beta(q(-\infty)),\beta(q(+\infty)))) < r
\end{equation}
for all $q$ and $n\in\Z$.

We will now establish a weak analogue of  
Lemma \ref{lem:qgeoconvinfgeneral},
taking into account the infinitesimal expansion property (\ref{eq:ano2}).
To make use of the expansion property, 
we need a version, in particular, a converse of Lemma~\ref{lem:expconlim}
for sequences close to parallel sets
and their infinitesimal contraction at infinity:

\begin{sublem}
\label{sublem:expconliminfini}
Let $(g_n)$ be a sequence in $G$ whose orbit sequence $(g_nx)$ for a point $x\in X$ 
is contained in a tubular neighborhood of the parallel set $P(\tau_-,\tau_+)$. 
Then the following are equivalent:

(i) The sequence $(g_n)$ is $\tau_{mod}$-regular and $g_n\tauto\tau_+$ conically.

(ii) The differentials $(dg_n^{-1})_{\tau_+}$ expand arbitrarily strongly, i.e.\
\begin{equation}
\label{eq:exptoinf}
\lim_{n\to+\infty}\eps(g_n^{-1},\tau_+)=+\infty
\end{equation}
with respect to a fixed background metric on $\Flagt$.
\end{sublem}
\proof
We proceed as in the proof of Corollary~\ref{cor:expand}.  
Let $\bar x\in P(\tau_-,\tau_+)$ denote the nearest point projection of $x$.
We write the $g_n$ as products $g_n=t_nb_n$ 
of transvections $t_n\in G$ along geodesics $l_n\subset P(\tau_-,\tau_+)$ 
through $\bar x$ 
and bounded isometries $b_n\in G$,
e.g.\ such that $d(b_n\bar x,\bar x)<2r$.
Then the $t_n$ fix $\tau_+$ on $\Flagt$, 
and the expansion factors $\eps(g_n^{-1},\tau_+)$
are the same as $\eps(t_n^{-1},\tau_+)$
up to bounded multiplicative error.
In view of $\eps(t_n^{-1},\tau_+)=\|(dt_n)_{\tau_+}\|^{-1}$,
condition (\ref{eq:exptoinf}) translates to the infinitesimal contraction condition
$$
\lim_{n\to+\infty}\|(dt_n)_{\tau_+}\|=0.
$$
According to Theorem~\ref{thm:infcontrtrans}, the latter is equivalent to 
$t_n\bar x$ being contained in $V(\bar x,\st(\tau_+))$ for large $n$ and 
$$
\lim_{n\to+\infty} d(t_n\bar x,\D V(\bar x,\st(\tau_+))) = \lim_{n\to+\infty} d(t_n^{-1}\bar x,\D V(\bar x,\st(\tau_-))) =+\infty.
$$
Since $d(g_n\bar x,t_n\bar x)=d(b_n\bar x,\bar x)$ is bounded,
this is in turn equivalent to $(g_n)$ being $\tau_{mod}$-regular and the constant sequence $(\tau_+)$ 
being a shadow of $(g_n)$ in $\Flagt$, i.e.\ $g_n\tauto\tau_+$ conically.
\qed

\medskip
The lemma follows by applying the sublemma to $q$.
Since $\Ga$ is word hyperbolic, 
every discrete geodesic ray $\N\to\Ga$ 
is contained in a uniform tubular neighborhood of a discrete geodesic
$\Z\to\Ga$. 
So, we may assume that $q$ extends to a discrete geodesic $q:\Z\to\Ga$ and, moreover, 
that it is normalized.
Since 
$$\lim_{i\to+\infty}\eps(q(n_i)^{-1},\beta(q(+\infty)))=+\infty
$$ 
for some sequence of indices $n_i\to+\infty$ in $\N$ by property (\ref{eq:ano2}),
the sublemma yields that the subsequence $(q(n_i))_{i\in\N}$ is $\tau_{mod}$-regular 
and $q(n_i)\tauto\beta(q(+\infty))$ conically as $i\to+\infty$. 
\qed

\medskip
At this stage we do not yet know that the entire subgroup $\Ga$ 
is $\tau_{mod}$-regular, nor do we know that $\be(\geo\Ga)=\LatGa$. 
The problem is that we have no uniform control yet on the distance 
of $q(n)x$ from the Weyl cone $V(q(0)x,\st(\beta(q(+\infty))))$. 
This will be our next aim. 

In the sequel, 
we will denote the nearest point projections 
of points $y$ in $X$ 
to the parallel set $P(\beta(q(-\infty)),\beta(q(+\infty)))$
by $\bar y$.

The proof of Sublemma~\ref{sublem:expconliminfini}
yields the following additional information.
Since $d(g_nx,t_n\bar x)\leq d(g_nx,g_n\bar x)+d(g_n\bar x,t_n\bar x)<3r$, 
it follows (for unnormalized $q$) 
that 
$$
\ol{q(n_i)x}\in V(\ol{q(0)x},\st(\beta(q(+\infty))))
$$
for all $i$, 
and 
$$
\lim_{i\to+\infty}d(\ol{q(n_i)x},\D V(\ol{q(0)x},\st(\beta(q(+\infty))))) =+\infty
$$

We fix a constant $d>>0$ and define for a discrete geodesic $q:\Z\to\Ga$
the {\em entry time} 
$T(q)\in\N$ 
as the smallest natural number $\geq1$ for which 
$$
\ol{q(T)x}\in V(\ol{q(0)x},\st(\beta(q(+\infty))))
$$
and 
$$
d(\ol{q(T)x},\D V(\ol{q(0)x},\st(\beta(q(+\infty)))))>d.
$$

\begin{lem}[Bounded entry time]\label{lem:bounded entry}
$T(q)$ is bounded above independently of $q$. 
\end{lem}
\proof
We first observe that the function $T$ on the space 
${\mathcal G}(\Ga)$
of discrete geodesics 
in $\Ga$, equipped with the topology of pointwise convergence, 
is upper semicontinuous,
because $\beta$ is continuous.
Since $T$ is $\Ga$-periodic and the $\Ga$-action on ${\mathcal G}(\Ga)$
is cocompact, the claim follows.
\qed

\medskip
Now we can strengthen Lemma~\ref{lem:weak anosov->weak as embedded} and 
show that the orbit maps of rays in $\Ga$ 
stay uniformly close to the Weyl cones associated to them by the boundary embedding, 
and that their nearest point projections to these Weyl cones move away from the boundaries of the cones
at a uniform linear rate:

\begin{lem}
\label{lem:cltcon}
(i) The distance of $q(n)x$  from the Weyl cone $V(q(0)x,\st(\beta(q(+\infty))))$ 
is uniformly bounded for all $n\geq 0$ independently of $q$.

(ii) The distance of $q(n)x$ from the boundary 
$$\D V(q(0)x,\st(\beta(q(+\infty))))$$
of the Weyl cone 
is at least $Cn-A$ for all $n\geq0$ 
with constants $C,A>0$ independent of $q$.
\end{lem}
\proof 
We inductively define a sequence of $i$-th entry times $T_i(q)$ 
for $i\geq0$ by $T_0(q):=0$ and
$$
T_{i+1}(q)-T_i(q) := T(q_i)
$$ 
for $i\geq1$,
where $T(q_i)$ is the entry time of the shifted discrete geodesic 
$q_i(n)=q(T_i(q)+n)$
with the same ideal endpoints. 
Then $T(q_i)$ and the increments $T_{i+1}(q)-T_i(q)$ 
are bounded above independently of $q$, 
and $T_i(q)\to+\infty$ monotonically as $i\to+\infty$, because $T\geq1$.
By taking into account that $q(\pm\infty)= q_i(\pm \infty)$, Ê
the definition of the $T(q_i)$ implies:
$$
\ol{q(T_{i+1})x}\in V(\ol{q(T_i)x},\st(\beta(q(+\infty)))) -N_d(\D V(\ol{q(T_i)x},\st(\beta(q(+\infty))))) 
$$
Corollary~\ref{cor:nestcone} yields that the cones are nested,
$$
V(\ol{q(T_{i+1})x},\st(\beta(q(+\infty)))) \subset V(\ol{q(T_i)x},\st(\beta(q(+\infty)))),
$$
and it hence follows that 
$$
d(\ol{q(T_i)x},\D V(\ol{q(0)x},\st(\beta(q(+\infty)))))>id.
$$
Note that $T_i$ grows uniformly linearly with $i$, i.e.\
$$
c^{-1}i Ê\le T_i\le c i
$$
for $i\geq1$ 
with a constant $c\ge 1$ independent of $q$.
Moreover, the increments $T_{i+1}-T_i$ are uniformly bounded. 
These two observations imply that:

(i) The distance of $\ol{q(n)x}$  from $V(\ol{q(0)x},\st(\beta(q(+\infty))))$ 
is uniformly bounded for all $n\geq 0$ independently of $q$;
and 

(ii) the distance of $\ol{q(n)x}$ from the complement 
$$P(\beta(q(-\infty)),\beta(q(+\infty)))-V(\ol{q(0)x},\st(\beta(q(+\infty))))$$
is $\geq Cn-A$ for all $n\geq0$ 
with constants $C,A>0$ independent of $q$.

The claim follows, because the Hausdorff distance between the cones 
$V(q(0)x,\st(\beta(q(+\infty))))$ and $V(\ol{q(0)x},\st(\beta(q(+\infty))))$
is $\leq d(q(0)x,\ol{q(0)x})<r$,
as well as $d(q(n)x,\ol{q(n)x})<r$.
\qed

\medskip
It follows that the group $\Ga$ is uniformly $\tau_{mod}$-regular. 

We will finally deduce from the last lemma that $\beta(\geo\Ga)=\LatGa$,
i.e.\ that $\beta$ maps {\em onto} the $\tau_{mod}$-limit set.
We proceed as in the proof of Proposition~\ref{prop:homeochlim}.

Suppose that $\tau\in\LatGa$ 
and let $(\ga_n)$ be a sequence in $\Ga$ with 
\begin{equation}
\label{eq:gaconvtau}
\ga_n\tauto\tau.
\end{equation}
Since $\Ga$ is word hyperbolic, 
there exists a sequence of discrete geodesic rays $q_n:\N\to\Ga$ 
initiating in $q_n(0)=1_{\Ga}$ and passing at uniformly bounded distance from $\ga_n$,
i.e.\ $(q_n(+\infty))$ is a shadow sequence for $(\ga_n)$ in $\geo\Ga$.
After passing to a subsequence, 
the $q_n$ converge (pointwise) to a ray $q:\N\to\Ga$.
Then $q_n(+\infty)\to q(+\infty)$ and 
\begin{equation}
\label{eq:endptsconv}
\beta(q_n(+\infty))\to\beta(q(+\infty)).
\end{equation} 
As a consequence of Lemma~\ref{lem:cltcon},
$(\beta(q_n(+\infty)))$ is a shadow sequence for $(\ga_n)$ in $\Flagt$.
Therefore, 
(\ref{eq:gaconvtau}) and (\ref{eq:endptsconv})
imply that 
$\beta(q_n(+\infty))\tauto\tau$,
cf.\ Lemma~\ref{lem:shadasuniq}. 
Hence $\beta(q(+\infty))=\tau$, 
i.e.\ $\tau$ is in the image of $\beta$.
\qed

\medskip
Thus, $\Ga$ is $\tau_{mod}$-asymptotically embedded. 
This concludes the proof of Theorem \ref{thm:weak anosov->as embedded}. 
\qed 

\medskip
We are now ready to prove the equivalence of three concepts,
namely of the (non-uniformly) Anosov property and asymptotic embeddedness.
Note that equivalences of the Anosov condition with other weakened forms of it
appear in \cite[sec.\ 6.1]{Labourie} and \cite[Prop.\ 3.16]{GW}.

\begin{thm}\label{thm:anosov iff morse}
For a discrete subgroup $\Ga\subset G$ which is non-elementary word hyperbolic, the following are equivalent:

1. $\Ga$ is $\tau_{mod}$-Anosov. 

2. $\Ga$ is non-uniformly $\tau_{mod}$-Anosov.

3. $\Ga$ is  $\tau_{mod}$-asymptotically embedded.
\end{thm}
\proof 
The implication 1$\Rightarrow$2 is immediate; implication 
2$\Rightarrow$3 is established in Theorem \ref{thm:weak anosov->as embedded}. 
It remains to prove that 3$\Ra$1.  

Suppose that $\Ga\subset G$ is $\tau_{mod}$-asymptotically embedded with boundary embedding $\be$. 
We use our results on the coarse geometry of asymptotically embedded subgroups, 
see Theorem~\ref{thm:qiembrtaureg}.
For a discrete geodesic ray $q:\N\to\Ga$ we know that its image $qx:\N\to X$ 
is a uniformly $\tau_{mod}$-regular discrete quasigeodesic ray
uniformly close to the Weyl cone $V(q(0)x, \st(\beta(q(+\infty))))$.
In particular, 
it moves away from the boundary $\D V(q(0)x, \st(\beta(q(+\infty))))$ of the cone with uniform linear speed. 
Therefore, Corollary~\ref{cor:expand} yields that the elements $q(n)^{-1}$ 
have uniform exponential expansion rate at 
$\beta(q(+\infty))\in \Flagt$,
$$
\eps(q(n)^{-1}, \beta(q(+\infty)))\ge Ae^{Cn}, 
$$ 
with constants $A,C>0$ independent of $q$. 
Hence, the subgroup $\Ga\subset G$ is $\tau_{mod}$-Anosov. 
\qed 

\medskip
As a corollary we obtain: 
\begin{cor}
Let $\Ga$ be a non-elementary word hyperbolic group. 
An isometric action $\rho:\Ga\acts X$ is $\tau_{mod}$-Anosov 
if and only if the kernel of $\rho$ is finite and $\rho(\Ga)\subset G$ 
is a $\tau_{mod}$-asymp\-to\-ti\-cally embedded discrete subgroup. 
\end{cor}

\section{Morse maps and actions}
\label{sec:morse}

In section~\ref{sec:coarse_asymptotically_embedded} we saw that asymptotically embedded subgroups satisfy an important coarse geometric property,
namely the orbit maps sends uniform quasigeodesic rays in $\Gamma$ to quasigeodesic rays $x_n=q(n) x$ in $X$ which stay within
bounded distance from Weyl cones.
We also saw that this property has important consequences: nondistortion of $\Gamma$ and uniform weak regularity.
We will call this property of sequences $(x_n)$ the \emph{Morse property}. The corresponding class of group actions will
be called \emph{Morse actions}. 
It is not hard to verify that, conversely,   Morse actions satisfy  asymptotic embeddedness. 

In this section we further investigate the Morse property of sequences and group actions. 
The main aim of this section is to establish 
a local criterion for being Morse.
To do so we introduce a local notion of  \emph{straightness} for sequences of points in $X$. 
Morse sequences are in general not straight, but they become straight after suitable modification, namely by sufficiently
coarsifying them
and then passing to the sequence of successive midpoints.
Conversely, the key result is that sufficiently spaced  straight sequences are Morse. We conclude that there is a local-to-global
implication for the Morse property.

As a consequence of the local-to-global criterion we establish that the Morse property for isometric group actions is an open condition.
Furthermore, for two nearby Morse actions, the actions on their $\tau_{mod}$-limit sets are also close, i.e.\ conjugate 
by an equivariant homeomorphism close to identity. In view of the equivalence of Morse property with the asymptotic
properties discussed earlier,
this implies structural stability for asymptotically embedded groups.
Another corollary of the local-to-global result is the algorithmic recognizability of Morse actions.

We conclude the section by illustrating our technique by constructing Morse-Schottky actions of free groups on higher rank symmetric spaces.
Unlike all previously known constructions, our proof does not rely on ping-pong, but
is purely geometric and proceeds by constructing  equivariant quasi-isometric embeddings of trees.

\bigskip

For the rest of this section we fix the following notation and conventions: 

Let $\tau_{mod}\subseteq\si_{mod}$ be an  
$\iota$-invariant face type.

We fix as auxiliary datum 
a $\iota$-invariant type $\zeta=\zeta_{mod}\in\interior(\tau_{mod})$. (We will omit the subscript in $\zeta_{mod}$ 
in order to avoid cumbersome notation for $\zeta$-angles.) 
This allows us to the define the {\em $\zeta$-angle} $\angle^\zeta$ 
and $\zeta$-Tits angle $\tangle^\zeta$,  see  equations 
\eqref{eq:zeta-angle} and \eqref{eq:zeta-tangle}. 
For a simplex $\tau\subset\geo X$ of type $\tau_{mod}$ 
we define $\zeta(\tau)\in\tau$ as the ideal point of type $\zeta_{mod}$. 
For a $\tau_{mod}$-regular unit tangent vector $v\in TX$ 
we denote by $\tau(v)\subset\geo X$ the unique simplex of type $\tau_{mod}$ 
such that ray $\rho_v$ with initial direction $v$ 
represents an ideal point in $\ost(\tau(v))$. 
We put $\zeta(v)=\zeta(\tau(v))$. 
Note that $\zeta(v)$ depends continuously on $v$.

In this section 
$\Theta,\Theta'\subset\ost(\tau_{mod})$ will denote 
$\iota$-invariant $\tau_{mod}$-convex compact subsets 
such that $\Theta\subset\interior(\Theta')$. 
The constants $L,A,D,\eps, \de,l,a,s,S$ 
are meant to be always strictly positive.

\subsection{A Morse Lemma for straight sequences}

In the following,
we consider finite or infinite sequences $(x_n)$ of points in $X$. 
\begin{dfn}[Straight and spaced sequence]
We call a sequence $(x_n)$ 
{\em $(\Theta,\eps)$-straight} 
if the segments $x_nx_{n+1}$ are $\Theta$-regular 
and 
\begin{equation*}
\angle_{x_n}^{\zeta}(x_{n-1},x_{n+1})\geq\pi-\eps
\end{equation*}
for all $n$. 
We call it {\em $l$-spaced} if the segments $x_nx_{n+1}$ 
have length $\geq l$.  
\end{dfn}
Note that every straight sequence 
can be extended to a biinfinite straight sequence. 

Straightness is a local condition.
The goal of this section is to prove the following 
local to global result 
asserting that sufficiently straight and spaced sequences 
satisfy a higher rank version of the Morse Lemma 
(for quasigeodesics in hyperbolic space). 

\begin{thm}[Morse Lemma for straight spaced sequences]
\label{thm:locstrimplcoastrseq}
For $\Theta,\Theta',\de$ there exist $l,\eps$ such that:

Every $(\Theta,\eps)$-straight $l$-spaced sequence $(x_n)$
is $\de$-close to a parallel set 
$P(\tau_-,\tau_+)$ with simplices $\tau_{\pm}$ of type $\tau_{mod}$, 
and it moves from $\tau_-$ to $\tau_+$ in the sense that 
its nearest point projection $\bar x_n$ to $P(\tau_-,\tau_+)$ satisfies 
\begin{equation}
\label{eq:movstr}
\bar x_{n\pm m}\in V(\bar x_n,\st_{\Theta'}(\tau_{\pm}))
\end{equation}
for all $n$ and $m\geq1$. 
\end{thm}
\begin{rem}[Global spacing]
\label{rem:globsp}
1. As a corollary of this theorem, we will show that straight spaced sequences are quasigeodesic: 
\begin{equation*}
d(x_n,x_{n+m})\geq clm-2\de
\end{equation*}
with a constant $c=c(\Theta')>0$. See Corollary \ref{cor:retractions}. 

2. Theorem \ref{thm:locstrimplcoastrseq} is a higher-rank generalization of two familiar facts from geometry of Gromov-hyperbolic geodesic metric spaces:  The fact that local quasigeodesics (with suitable parameters) are global quasigeodesics and the Morse lemma stating that 
quasigeodesics stay uniformly close to geodesics. In the higher rank, quasigeodesics, of course, need not be close to geodesics, but, instead (under the straightness assumption), are close to  parallel sets. 
\end{rem}

In order to prove the theorem, 
we start by considering one-sided infinite sequences 
and prove that they keep moving away from an ideal simplex 
of type $\tau_{mod}$
if they do so initially. 

\begin{dfn}[Moving away from an ideal simplex]
Given a face $\tau\subset \tits X$ of type $\tau_{mod}$ and distinct points $x, y\in X$, define the angle 
$$
\angle_{x}^{\zeta}(\tau, y):= \angle_{x}(z, y)
$$ 
where $z$ is a point (distinct from $x$) on the geodesic ray $x\xi$, where $\xi\in \tau$ is the point of type $\zeta$. 

We say that a sequence $(x_n)$ 
{\em moves $\eps$-away} 
from a simplex $\tau$ of type $\tau_{mod}$ if  
\begin{equation*}
\angle_{x_n}^{\zeta}(\tau,x_{n+1})\geq\pi-\eps
\end{equation*}
for all $n$. 
\end{dfn}

\begin{lem}[Moving away from ideal simplices]
\label{lem:keepsmoving}
For small $\eps$ and large $l$, 
$\eps\leq\eps_0$ and $l\geq l(\eps,\Theta)$, 
the following holds:

If the sequence $(x_n)_{n\geq0}$
is $(\Theta,\eps)$-straight $l$-spaced 
and if 
\begin{equation*}
\angle_{x_0}^{\zeta}(\tau,x_1)\geq\pi-2\eps,
\end{equation*}
then $(x_n)$ moves $\eps$-away from $\tau$. 
\end{lem}
\proof
By Lemma~\ref{lem:distangdec}(ii), 
the unit speed geodesic segment $c:[0,t_1]\to X$ from $p(0)$ to $p(1)$ 
moves $\eps(d(2\eps))$-away from $\tau$ at all times, 
and $\eps'(2\eps,\Theta,l)$-away at times $\geq l$, 
which includes the final time $t_1$. 
For $l(\eps,\Theta)$ sufficiently large, 
we have $\eps'(2\eps,\Theta,l)\leq\eps$. 
Then $c$ moves $\eps$-away from $\tau$ at time $t_1$, 
which means that 
$\angle_{x_1}^{\zeta}(\tau,x_0)\leq\eps$. 
Straightness at $x_1$ and the triangle inequality yield that again
$\angle_{x_1}^{\zeta}(\tau,x_2)\geq\pi-2\eps$. 
One proceeds by induction. 
\qed

\medskip
Note that there do exist simplices $\tau$ 
satisfying the hypothesis of the previous lemma. 
For instance, 
one can extend the initial segment $x_0x_1$ backwards to infinity 
and choose $\tau=\tau(x_1x_0)$. 

\medskip
Now we look at {\em biinfinite} sequences. 

We assume in the following that $(x_n)_{n\in\Z}$
is $(\Theta,\eps)$-straight $l$-spaced
for small $\eps$ and large $l$. 
As a first step, 
we study the asymptotics of such sequences 
and use the argument for Lemma~\ref{lem:keepsmoving} 
to find a pair of opposite ideal simplices $\tau_{\pm}$ 
such that $(x_n)$ moves from $\tau_-$ towards $\tau_+$. 

\begin{lem}[Moving towards ideal simplices]
\label{lem:locstrcpclpar}
For small $\eps$ and large $l$, 
$\eps\leq\eps_0$ and $l\geq l(\eps,\Theta)$, 
the following holds: 

There exists a pair of opposite simplices $\tau_{\pm}$ 
of type $\tau_{mod}$ 
such that the inequality
\begin{equation}
\label{ineq:tauminusplusatnb}
\angle_{x_n}^{\zeta}(\tau_{\mp},x_{n\pm1})\geq\pi-2\eps
\end{equation}
holds for all $n$.
\end{lem}
\proof
1. For every $n$ define a compact set 
$C^\mp_n\subset\Flag(\tau_{mod})$
$$
C^\pm_n=\{\tau_\pm : \angle_{x_n}^{\zeta}(\tau_{\pm},x_{n\mp 1})\geq\pi-2\eps\}. 
$$
As in the proof of Lemma~\ref{lem:keepsmoving}, 
straightness at $x_{n+1}$ implies that 
$C^-_n\subset C^-_{n+1}$. 
Hence the family $\{C^-_n\}_{n\in \Z}$ form a nested sequence of nonempty compact subsets
and therefore have nonempty intersection containing a simplex $\tau_-$.  
Analogously, there exists a simplex $\tau_+$ which belongs to $C^+_n$ for all $n$.  

2. It remains to show that the simplices $\tau_-, \tau_+$ are antipodal. Using straightness and the triangle inequality, we see that 
\begin{equation*}
\label{ineq:tauminusplusatn}
\angle_{x_n}^\zeta(\tau_-,\tau_+)\geq\pi-5\eps 
\end{equation*}
for all $n$. Hence, if $5\eps<\eps(\zeta)$, then the simplices $\tau_-, \tau_+$ are antipodal in view of Remark \ref{rem:antip}. \qed

\medskip
The pair of opposite simplices $(\tau_-, \tau_+)$ which we found determines a parallel set in $X$. 
The second step is to show that $(x_n)$ is uniformly close to it. 

\begin{lem}[Close to parallel set]
\label{lem:closeparsetb}
For small $\eps$ and large $l$, 
$\eps\leq\eps(\de)$ and $l\geq l(\Theta,\de)$, 
the sequence $(x_n)$ is $\de$-close to $P(\tau_-,\tau_+)$. 
\end{lem}
\proof The statement follows from the combination of 
the inequality \eqref{ineq:tauminusplusatn} (in the second part of the proof of Lemma \ref{lem:locstrcpclpar}) and Lemma~\ref{lem:distangcontr}.
\qed

\medskip
The third and final step 
is to show that the nearest point projection $(\bar x_n)$ of $(x_n)$ 
to $P(\tau_-,\tau_+)$ 
moves from $\tau_-$ towards $\tau_+$. 
\begin{lem}[Projection moves towards ideal simplices]
\label{lem:projpathstrb}
For small $\eps$ and large $l$, 
$\eps\leq\eps_0$ and $l\geq l(\eps,\Theta,\Theta')$, 
the segments $\bar x_n\bar x_{n+1}$ are $\Theta'$-regular
and 
\begin{equation*}
\label{ineq:tauminusplusatnproj}
\angle_{\bar x_n}^{\zeta}(\tau_-,\bar x_{n+1})=\pi
\end{equation*}
for all $n$. 
\end{lem}
\proof
By the previous lemma, 
$(x_n)$ is $\delta_0$-close to $P(\tau_-,\tau_+)$
if $\eps_0$ is sufficiently small and $l$ is sufficiently large. 
Since $x_nx_{n+1}$ is $\Theta$-regular, 
the triangle inequality for $\De$-lengths yields that 
the segment $\bar x_n\bar x_{n+1}$ is $\Theta'$-regular, 
again if $l$ is sufficiently large.

Let $\xi_+$ denote the ideal endpoint of the ray extending this segment, 
i.e.\ $\bar x_{n+1}\in\bar x_n\xi_+$. 
Then $x_{n+1}$ is $2\delta_0$-close to the ray $x_n\xi_+$. 
We obtain that 
\begin{equation*}
\tangle^{\zeta}(\tau_-,\xi_+)\geq
\angle^{\zeta}_{x_n}(\tau_-,\xi_+)\simeq
\angle^{\zeta}_{x_n}(\tau_-,x_{n+1})\simeq\pi
\end{equation*}
where the last step follows from inequality (\ref{ineq:tauminusplusatnb}). 
The discreteness of Tits distances between ideal points 
of fixed type $\zeta$ implies that in fact
\begin{equation*}
\tangle^{\zeta}(\tau_-,\xi_+)=\pi, 
\end{equation*}
i.e.\ the ideal points $\zeta(\tau_-)$ and $\zeta(\xi_+)$ 
are antipodal. 
But the only simplex opposite to $\tau_-$ in 
$\geo P(\tau_-,\tau_+)$ is $\tau_+$, 
so $\tau(\xi_+)=\tau_+$ and 
\begin{equation*}
\angle_{\bar x_n}^{\zeta}(\tau_-,\bar x_{n+1})=
\angle_{\bar x_n}^{\zeta}(\tau_-,\xi_+)
=\pi,
\end{equation*}
as claimed.
\qed

\medskip
{\em Proof of Theorem~\ref{thm:locstrimplcoastrseq}.}
It suffices to consider biinfinite sequences. 

The conclusion of 
Lemma~\ref{lem:projpathstrb} 
is equivalent to 
$\bar x_{n+1}\in V(\bar x_n,\st_{\Theta'}(\tau_+))$. 
Combining Lemmas~\ref{lem:closeparsetb} and~\ref{lem:projpathstrb}, 
we thus obtain the theorem for $m=1$. 

The convexity of $\Theta'$-cones, 
cf.\ Proposition~\ref{prop:thconeconv}, 
implies that 
\begin{equation*}
V(\bar x_{n+1},\st_{\Theta'}(\tau_+))\subset 
V(\bar x_n,\st_{\Theta'}(\tau_+)),
\end{equation*}
and the assertion follows for all $m\geq1$ by induction.  
\qed

\begin{rem}
\label{rem:xinotconv}
The conclusion of the theorem implies that 
$x_{\pm n}\tauto\tau_{\pm}$ as $n\to+\infty$. 
However, 
the $x_n$ do in general not converge at infinity,
but accumulate at a compact subset of $\st_{\Theta'}(\tau_{\pm})$. 
\end{rem}

\subsection{Lipschitz retractions}

Consider a (possibly infinite) closed interval $J$ in $\R$; we will assume that $J$ has integer or infinite bounds.  
Suppose that $p: J\cap \Z\to P=P(\tau_-, \tau_+)\subset X$ is an $l$-separated, $\lambda$-Lipschitz, $(\Theta,0)$-straight coarse sequence 
pointing away from $\tau_-$ and towards $\tau_+$. We extend $p$ to a piecewise-geodesic map $p: J\to P$ by sending intervals $[n, n+1]$ to geodesic segments $p(n)p(n+1)$ via affine maps. We retain the name $p$ for the extension. 

\begin{lem}
There exists $L=L(l,\lambda,\Theta)$ and an  $L$-Lipschitz retraction of $X$ to $p$, i.e., an $L$-Lipschitz map $r: X\to J$ so that $r\circ p=Id$. 
In particular, $p: J\cap \Z\to X$ is a $(\bar L, \bar A)$-quasigeodesic, where $\bar L, \bar A$ depend only on $l, \lambda, \Theta$.  
\end{lem}
\proof It suffices to prove existence of a retraction. Since $P$ is convex in $X$, 
it suffices to construct a map $P\to J$. Pick a generic point $\xi=\xi_+\in \tau_+$ and let $b_\xi: P\to \R$  denote the Busemann function normalized so that $b_\xi(p(z))=0$ for some $z\in J\cap \Z$.  Then the $\Theta$-regularity assumption on $p$ implies that the slope of the piecewise-linear function $b_\xi\circ p: J\to \R$ is strictly positive, bounded away from $0$. The assumption that $p$ is $l$-separated $\la$-Lipschitz implies that
$$
l\le |p'(t)|\le \la  
$$
for each $t$ (where the derivative exists). The straightness assumption on $p$ implies that the function $h:= b_\xi\circ p: J\to \R$ is strictly increasing. By combining these observations, we conclude that $h$ is an $L$-biLipschitz homeomorphism for some 
$L=L(l,\la,\Theta)$. Lastly, we define 
$$
r: P\to J, \quad r=h^{-1}\circ b_\xi. 
$$ 
Since $b_\xi$ is $1$-Lipschitz, the map $r$ is $L$-Lipschitz. By construction, $r\circ p=Id$. \qed 

\begin{cor}\label{cor:retractions}
Suppose that $p: J\cap \Z\to X$ is a $l$-separated, $\lambda$-Lipschitz, $(\Theta,\eps)$-regular straight coarse sequence. Pick some $\Theta'$ such that 
$\Theta\subset int(\Theta')$ and let $\delta=\delta(l, \Theta, \Theta', \eps)$ be the constant as in Theorem \ref{thm:locstrimplcoastrseq}. Then for 
$L=L(l-2\delta,\la+2\delta,\Theta')$  we have: 

1. There exists an $(L,2\delta)$-coarse Lipschitz retraction $X\to J$. 

2. The map $p$ is an $(L',A')$-quasigeodesic for $L', A'$ depending only on  $l, \Theta, \Theta', \eps$.  
\end{cor}
\proof The statement immediately follows the above lemma combined with Theorem \ref{thm:locstrimplcoastrseq}.  \qed

\subsection{Morse quasigeodesics}
\label{sec:morseqg}

According to Theorem~\ref{thm:locstrimplcoastrseq}, 
sufficiently spaced straight sequences 
satisfy a Morse Lemma. 
If the spacing is also bounded above, 
then these sequences are quasigeodesics (Corollary \ref{cor:retractions}). 
This motivates considering quasigeodesics 
which satisfy a higher rank version of the Morse Lemma 
as it appears in the conclusion of Theorem~\ref{thm:locstrimplcoastrseq}. 
\begin{dfn}[Morse quasigeodesic]
\label{dfn:mqg}
An $(L,A,\Theta,D)$-{\em Morse quasigeodesic} in $X$ 
is an $(L,A)$-quasigeodesic $p:I\to X$ 
such that for all $t_1,t_2\in I$ 
the subpath $p|_{[t_1,t_2]}$ is $D$-close to a $\Theta$-diamond 
$\diamo_{\Theta}(x_1,x_2)$
with $d(x_i,p(t_i))\leq D$.
\end{dfn}

Note that each quasigeodesic satisfying Theorem~\ref{thm:locstrimplcoastrseq} is Morse; the converse is also true 
as we will see in Lemma \ref{lem:morsecon}.  

We will now prove that, conversely, 
the Morse property implies straightness in a suitable sense, 
namely that for sufficiently spaced quadruples 
the associated midpoint triples are arbitrarily straight. 
(For the quadruples themselves this is in general not true.)

\begin{dfn}[Quadruple condition]
\label{dfn:quad}
For points $x, y\in X$ we let $\midp(x,y)$ denote the midpoint of the geodesic segment $xy$. 
A map $p:I\to X$ satisfies the 
{\em $(\Theta,\eps,l,s)$-quadruple condition}
if for all $t_1,t_2,t_3,t_4\in I$
with $t_2-t_1,t_3-t_2,t_4-t_3\geq s$ 
the triple of midpoints
\begin{equation*}
(\midp(t_1,t_2),\midp(t_2,t_3),\midp(t_3,t_4))
\end{equation*}
is $(\Theta,\eps)$-straight and $l$-spaced. 
\end{dfn}
\begin{prop}[Morse implies quadruple condition]
\label{prop:morseimplquad}
For $L,A,\Theta,\Theta',D,\eps,l$ 
exists a scale $s=s(L,A,\Theta,\Theta',D,\eps,l)$ 
such that every $(L,A,\Theta,D)$-Morse quasigeodesic 
satisfies the $(\Theta',\eps,l,s')$-quadruple condition for every $s'\ge s$.  
\end{prop}
\proof
Let $p:I\to X$ be an $(L,A,\Theta,D)$-Morse quasigeodesic, 
and let $t_1,\dots,t_4\in I$
such that $t_2-t_1,t_3-t_2,t_4-t_3\geq s$. 
We abbreviate $p_i:=p(t_i)$ and $m_i=\midp(p_i,p_{i+1})$. 

Regarding straightness, 
it suffices to show that 
the segment $m_2m_1$ is $\Theta'$-regular 
and that $\angle^{\zeta}_{m_2}(p_2,m_1)\leq\frac{\eps}{2}$
provided that $s$ is sufficiently large 
in terms of the given data. 

By the Morse property,
there exists a diamond
$\diamo_{\Theta}(x_1,x_3)$
such that $d(x_1,p_1),d(x_3,p_3)\leq D$ 
and $p_2\in N_D(\diamo_{\Theta}(x_1,x_3))$.
The diamond spans a unique parallel set $P(\tau_-,\tau_+)$. 
(Necessarily, 
$x_3\in V(x_1,\st_{\Theta}(\tau_+))$ and 
$x_1\in V(x_3,\st_{\Theta}(\tau_-))$.)

We denote by $\bar p_i$ and $\bar m_i$
the projections of $p_i$ and $m_i$ to the parallel set. 

We first observe 
that $m_2$ (and $m_3$) is arbitrarily close to the parallel set 
if $s$ is large enough. 
If this were not true,
a limiting argument would produce a geodesic line 
at strictly positive finite Hausdorff distance 
$\in(0,D]$ 
from $P(\tau_-,\tau_+)$ 
and asymptotic to ideal points in 
$\st_{\Theta}(\tau_{\pm})$. 
However,
all lines asymptotic to ideal points in
$\st_{\Theta}(\tau_{\pm})$ 
are contained in $P(\tau_-,\tau_+)$. 

Next,
we look at the directions of the segments $\bar m_2\bar m_1$ 
and $\bar m_2\bar p_2$
and show that they have the same $\tau$-direction. 
Since $\bar p_2$ is $2D$-close to 
$V(\bar p_1,\st_{\Theta}(\tau_+))$, 
we have that the point 
$\bar p_1$ is $2D$-close to 
$V(\bar p_2,\st_{\Theta}(\tau_-))$, 
and hence also 
$\bar m_1$ is $2D$-close to 
$V(\bar p_2,\st_{\Theta}(\tau_-))$.
Therefore, 
$\bar p_1,\bar m_1\in V(\bar p_2,\st_{\Theta'}(\tau_-))$ 
if $s$ is large enough. 
Similarly, 
$\bar m_2\in V(\bar p_2,\st_{\Theta'}(\tau_+))$
and hence 
$\bar p_2\in V(\bar m_2,\st_{\Theta'}(\tau_-))$. 
The convexity of $\Theta'$-cones, 
see Proposition~\ref{prop:thconeconv}, 
implies that also 
$\bar m_1\in V(\bar m_2,\st_{\Theta'}(\tau_-))$. 
In particular,
$\angle^{\zeta}_{\bar m_2}(\bar p_2,\bar m_1)=0$ 
if $s$ is sufficiently large. 

Since $m_2$ is arbitrarily close to the parallel set 
if $s$ is sufficiently large, 
it follows by another limiting argument that 
$\angle^{\zeta}_{m_2}(p_2,m_1)\leq\frac{\eps}{2}$
if $s$ is sufficiently large. 

Regarding the spacing, 
we use that 
$\bar m_1\in V(\bar p_2,\st_{\Theta'}(\tau_-))$
and 
$\bar m_2\in V(\bar p_2,\st_{\Theta'}(\tau_+))$.
It follows that 
\begin{equation*}
d(\bar m_1,\bar m_2) \geq c\cdot (d(\bar m_1,\bar p_2)+d(\bar p_2,\bar m_2))
\end{equation*}
with a constant $c=c(\Theta')>0$, 
and hence that 
$d(m_1,m_2)\geq l$ 
if $s$ is sufficiently large. 
\qed

\medskip
Theorem~\ref{thm:locstrimplcoastrseq} and Proposition~\ref{prop:morseimplquad} 
tell that the Morse property for quasigeodesics is equivalent 
to straightness (of associated spaced sequences of points). 
Since straightness is a local condition, 
this leads to a local to global result for Morse quasigeodesics, 
namely that the Morse property holds globally 
if it holds locally up to a sufficiently large scale. 
\begin{dfn}[Local Morse quasigeodesic]
An $(L,A,\Theta,D,S)$-{\em local Morse quasigeo\-de\-sic} in $X$ 
is a map $p:I\to X$ 
such that for all $t_0$ the subpath 
$p|_{[t_0,t_0+S]}$ is an $(L,A,\Theta,D)$-Morse quasigeodesic. 
\end{dfn}
Note that local Morse quasigeodesics are uniformly coarse Lipschitz. 
\begin{thm}[Local to global for Morse quasigeodesics]
\label{thm:locglobmqg}
For $L,A,\Theta,\Theta',D$ exist $S,L',A',D'$ such that 
every $(L,A,\Theta,D,S)$-local Morse quasigeo\-de\-sic in $X$ 
is an $(L',A',\Theta',D')$-Morse quasigeodesic. 
\end{thm}
\proof
We choose an auxiliary Weyl convex subset $\Theta''$ 
depending on $\Theta,\Theta'$ 
such that 
$\Theta\subset\interior(\Theta'')$ and 
$\Theta''\subset\interior(\Theta')$. 

Let $p:I\to X$ be an $(L,A,\Theta,D,S)$-local Morse quasigeo\-de\-sic. 
We consider its coarsification on a (large) scale $s$ 
and the associated midpoint sequence, 
i.e.\ we put $p^s_n=p(ns)$ and $m^s_n=\midp(p^s_n,p^s_{n+1})$. 
Whereas the coarsification itself 
does in general not become arbitrarily straight 
as the scale $s$ increases, 
this is true for its midpoint sequence
due to Proposition~\ref{prop:morseimplquad}. 
We want it to be sufficiently straight and spaced 
so that we can apply to it the Morse Lemma from 
Theorem~\ref{thm:locstrimplcoastrseq}. 
Therefore we first fix an auxiliary constant $\de$, 
and further auxiliary constants $l,\eps$ 
as determined by Theorem~\ref{thm:locstrimplcoastrseq} 
in terms of $\Theta',\Theta''$ and $\de$. 
Then 
Proposition~\ref{prop:morseimplquad} applied to the 
$(L,A,\Theta,D)$-Morse quasigeodesics $p|_{[t_0,t_0+S]}$ 
yields that 
$(m^s_n)$ is $(\Theta'',\eps)$-straight and $l$-spaced 
if $S\geq 3s$ and the scale $s$ is large enough 
depending on $L,A,\Theta,\Theta'',D,\eps,l$. 

Now we can apply Theorem~\ref{thm:locstrimplcoastrseq} to $(m^s_n)$. 
It yields 
a nearby sequence $(\bar m^s_n)$, 
$d(\bar m^s_n,m^s_n)\leq\de$, 
with the following property:
For all $n_1<n_2<n_3$ 
the segments $\bar m^s_{n_1}\bar m^s_{n_3}$ 
are uniformly regular 
and the points $m^s_{n_2}$ are $\de$-close to the diamonds
$\diamo_{\Theta'}(\bar m^s_{n_1},\bar m^s_{n_3})$. 

Since the subpaths $p|_{[ns,(n+1)s]}$ filling in $(p^s_n)$ 
are $(L,A)$-quasigeodesics (because $S\geq s$), 
and it follows that 
for all $t_1,t_2\in I$ 
the subpaths $p|_{[t_1,t_2]}$ are 
$D'$-close to $\Theta'$-diamonds 
with $D'$ depending on $L,A,s$. 

The conclusion of Theorem~\ref{thm:locstrimplcoastrseq} 
also implies a global spacing for the sequence $(m^s_n)$, 
compare Remark~\ref{rem:globsp}, i.e.\ 
$d(m^s_n,m^s_{n'})\geq c\cdot |n-n'|$
with a positive constant $c$ depending on $\Theta',l$. 
Hence $p$ is a global $(L',A')$-quasigeodesic 
with $L',A'$ depending on 
$L,A,s,c$. 

Combining this information, 
we obtain that $p$ is an $(L',A',\Theta',D')$-Morse quasigeodesic 
for certain constants $L',A'$ and $D'$ 
depending on $L,A,\Theta,\Theta'$ and $D$, 
provided that the scale $S$ is sufficiently large 
in terms of the same data. 
\qed

\medskip
We discuss now the asymptotics of Morse quasigeodesics. 

There is much freedom for the asymptotic behavior 
of arbitrary quasigeodesics in euclidean spaces,
and therefore also in symmetric spaces of higher rank.  
However, 
the asymptotic behavior of Morse quasigeodesics 
is as restricted as for quasigeodesics in rank one symmetric spaces. 

Morse quasirays do in general not converge at infinity, 
but they $\tau_{mod}$-converge at infinity, 
compare Remark~\ref{rem:xinotconv}. 
This is a consequence of:
\begin{lem}[Conicality]
\label{lem:morsecon}
Every Morse quasiray $p:[0,\infty)\to X$ 
is uniformly Hausdorff close to a cone $V(p(0),\st_{\Theta}(\tau))$
for a unique simplex $\tau$ of type $\tau_{mod}$. 
\end{lem}
\proof
The subpaths $p|_{[0,t_0]}$ are uniformly Hausdorff close 
to $\Theta$-diamonds. 
These subconverge to a cone $V(x,\st_{\Theta}(\tau))$ 
$x$ uniformly close to $p(0)$ and $\tau$ a simplex of type $\tau_{mod}$. 
This establishes the existence. 
Since $p(n)\tauto\tau$, 
the uniqueness of $\tau$ 
may be deduced from the uniqueness of $\tau_{mod}$-limits, 
cf.\ Lemma~\ref{lem:shadasuniq}.
\qed

\begin{dfn}[End of Morse quasiray]
We call the unique simplex given by the previous lemma 
the {\em end} of the Morse quasiray $p:[0,\infty)\to X$ 
and denote it by 
\begin{equation*}
p(+\infty)\in\Flag(\tau_{mod}).
\end{equation*}
\end{dfn}
Hausdorff close Morse quasirays have the same end, 
again by Lemma~\ref{lem:shadasuniq}, 
and this lemma also implies the continuous dependence of the end 
on the Morse quasiray:
\begin{lem}[Continuity of end]
\label{lem:cont-dependence}
The assignment $p\mapsto p(+\infty)$ 
is a continuous map from the space of 
Morse quasirays $[0,+\infty)\to X$ with fixed data, 
equipped with the topology of  pointwise convergence (equivalently, uniform convergence on compacts), 
to $\Flag(\tau_{mod})$. 
\end{lem}

\subsection{Morse embeddings}

We consider maps into $X$ from metric spaces 
which are coarsely geodesic in the sense that their points can be connected 
by uniform quasigeodesics. 
\begin{dfn}[Quasigeodesic metric space]
A metric space is called {\em $(l,a)$-quasigeodesic} 
if pairs of points can be connected by $(l,a)$-quasigeodesics. 
It is called {\em quasigeodesic} 
if it is $(l,a)$-quasigeodesic for some parameters $l,a$.
\end{dfn}
The quasigeodesic spaces considered in this paper 
are discrete groups equipped with word metrics. 
\begin{dfn}[Morse embedding]
\label{dfn:morseemb}
A {\em Morse embedding} 
from a quasigeodesic space $Z$ into $X$ is a map 
$f:Z\to X$ 
which sends uniform quasigeodesics in $Z$ 
to uniform Morse quasigeodesics in $X$. 
We call it a {\em $\Theta$-Morse embedding} 
if it sends uniform quasigeodesics 
to uniform $\Theta$-Morse quasigeodesics. 
\end{dfn}
Thus, to be a Morse embedding means that for any parameters $l,a$ 
the $(l,a)$-quasigeodesics in $Z$ are mapped to 
$(L,A,\Theta,D)$-Morse quasigeodesics in $X$ 
with the parameters $L,A,\Theta,D$ depending on $l,a$. 

Note that Morse embeddings are quasi-isometric embeddings. 

Our definition is chosen so that 
it depends only on the quasi-isometry class of $Z$ 
whether a map $f:Z\to X$ is a ($\Theta$-)Morse embedding,
i.e.\ 
the precomposition of a ($\Theta$-)Morse embedding with a quasi-isometry 
is again ($\Theta$-)Morse. 
For this to be true 
is why we require control on the images of quasigeodesics 
of arbitrarily bad quality. 

However, as we observe next, 
in the case of maps from Gromov hyperbolic spaces 
control on the images of quasigeodesics of a fixed quality suffices. 
This is due to the 
Morse Lemma for quasigeodesics in Gromov hyperbolic spaces. 
We recall that it asserts that quasigeodesics with the same endpoints 
are uniformly close to each other, 
the closeness depending on the quasi-isometry and hyperbolicity constants.
\begin{lem}
\label{lem:morsefromhyp}
Let $f:Z\to X$ be a map from a Gromov hyperbolic space $Z$ into $X$. 
If $Z$ is $(l,a)$-quasigeodesic
and if $f$ sends $(l,a)$-quasigeodesics 
to uniform $\Theta$-Morse quasigeodesics, 
then $f$ is a $\Theta$-Morse embedding. 
\end{lem}
\proof
This is a consequence of the definition of Morse quasigeodesics, 
see Definition~\ref{dfn:mqg} 
and the Morse Lemma applied to $Z$. 
\qed

\medskip
We now deduce from our local to global result for Morse quasigeodesics, 
see Theorem~\ref{thm:locglobmqg}, 
a local to global result for Morse embeddings. 

Since we need to fix one scale of localness, 
we can expect a local to global control 
for the $f$-images of quasigeodesics in $Z$ 
only if they have a certain fixed quality. 
This is why we need to restrict to maps from Gromov hyperbolic spaces. 
\begin{dfn}[Local Morse embedding]
We call a map $f:Z\to X$ from a quasigeodesic space $Z$ into $X$ 
an {\em $(l,a,L,A,\Theta,D,S)$-local Morse embedding} 
if $Z$ is $(l,a)$-quasigeodesic
and if for any $(l,a)$-quasigeodesic $q:I\to Z$ 
defined on an interval $I$ of length $\leq S$ 
the image path $f\circ q$ is an 
$(L,A,\Theta,D)$-Morse quasigeodesic in $X$. 
\end{dfn}

\begin{thm}[Local-to-global for Morse embeddings of Gromov hyperbolic spaces]
\label{thm:locglobmqiembhyp}
For $l,a,L,A,\Theta,\Theta',D$ exists a scale $S$ 
such that every 
$(l,a,L,A,\Theta,D,S)$-local Morse embedding 
from a quasigeodesic  Gromov hyperbolic space into $X$ 
is a $\Theta'$-Morse embedding. 
\end{thm}
\proof
Let $f:Z\to X$ denote the local Morse embedding. 
It sends every $(l,a)$-quasigeodesic $q:I\to Z$ 
to a $(L,A,\Theta,D,S)$-local Morse quasigeodesic $p=f\circ q$ in $X$. 
By Theorem~\ref{thm:locglobmqg}, 
$p$ is $(L',A',\Theta',D')$-Morse 
if $S$ is sufficiently large,
where $L',A',D'$ depend on the given data. 
Lemma~\ref{lem:morsefromhyp} implies 
that $f$ is a $\Theta'$-Morse embedding. 
\qed

\medskip
We now discuss the asymptotics of Morse embeddings. 

Morse embeddings $f:Z\to X$ map sufficiently spaced pairs of points 
to uniformly $\tau_{mod}$-regular pairs of points. 
Therefore, 
their images accumulate in the $\tau_{mod}$-regular part of $\geo X$ 
and there is a well-defined {\em flag limit set}
\begin{equation}
\label{eq:limsmemb}
\La_{\tau_{mod}}(f)\subset\Flag(\tau_{mod}) .
\end{equation} 
For Morse embeddings $f:Z\to X$ of Gromov hyperbolic spaces $Z$ 
we obtain, 
by applying our discussion of the asymptotics of Morse quasirays, 
a well-defined continuous boundary map at infinity 
\begin{equation}
\label{eq:bdmapmemb}
\geo f:\geo Z\to\Flag(\tau_{mod}) 
\end{equation}
which sends the ideal endpoint $q(+\infty)$ 
of a quasiray $q:[0,+\infty)\to Z$ 
to the end $(f\circ q)(+\infty)$ of the image Morse quasiray. 
The continuity of $\geo f$ follows from Lemma~\ref{lem:cont-dependence}. 
\begin{prop}[Asymptotic properties of Morse embeddings]
\label{prop:asymorseemb}
Let $f:Z\to X$ be a Morse embedding from a 
proper, quasigeodesic Gromov hyperbolic space $Z$.

(i) $f$ is ${\tau_{mod}}$-conical in the sense that $f$ sends each quasigeodesic ray $q$ in $Z$ to a quasigeodesic ray $p$ in $X$ so that the sequence 
$(p(n))$ converges conically to the simplex $\geo f(q(\infty))$.  

(ii) 
$\geo f$ is antipodal, 
i.e.\ it maps distinct ideal points to antipodal simplices. 

(iii)
$f$ has nice asymptotics, 
i.e.\ $\geo f$ is a homeomorphism onto $\La_{\tau_{mod}}(f)$.
\end{prop}
\proof
(i)
This is a consequence of 
Lemma~\ref{lem:morsecon}
applied to the $f$-images of quasirays in $Z$.

(ii)
By hyperbolicity and the quasigeodesic properties of $Z$,   
any two distinct points in $\geo Z$
can be connected by a quasigeodesic $q:\R\to Z$. 
Then $f\circ q$ is a Morse quasigeodesic. 
Since Morse quasigeodesics are Hausdorff close to 
biinfinite straight spaced sequences 
(Proposition~\ref{prop:morseimplquad}), 
and such sequences satisfy a Morse Lemma 
(Theorem~\ref{thm:locstrimplcoastrseq})
and therefore $\tau_{mod}$-converge to a pair of opposite simplices, 
the assertion follows. 

(iii)
Our argument follows the end of the proof of 
Theorem~\ref{thm:weak anosov->as embedded}.
By construction of the boundary map $\geo f$,
its image is contained in the limit set $\La_{\tau_{mod}}(f)$. 
It is injective by antipodality. 
To prove surjectivity, 
let $\tau\in\La_{\tau_{mod}}(f)$, 
and let $z_n\to\infty$ be a sequence in $Z$ 
so that $f(z_n)\to\tau$ in the sense of $\tau_{mod}$-convergence. 
Since $Z$ is quasigeodesic, 
there exists a sequence of uniform quasigeodesic segments 
$q_n:[0,l_n]\to Z$ 
connecting a base point to $z_n$.
After passing to a subsequence, 
the $q_n$ converge to a quasiray $q_{\infty}:[0,+\infty)\to Z$ 
because $Z$ is proper. 
We need to verify that 
$(f\circ q_{\infty})(+\infty))=\tau$. 

Let $t_n\to+\infty$ be a sequence of times $t_n\leq l_n$.
Then the Morse property of $f\circ q_n$ 
and Lemma~\ref{lem:shadasuniq} 
imply that 
$(f\circ q_n)(t_n)\to\tau$. 
We can choose the $t_n$ so that $q_n(t_n)$ is uniformly close to 
$q_{\infty}$. 
The Morse property of $f\circ q_{\infty}$ implies 
that $(f\circ q_n)(t_n)$ is uniformly close to a $\Theta$-cone 
$V(x,\st_{\Theta}((f\circ q_{\infty})(+\infty)))$. 
Applying Lemma~\ref{lem:shadasuniq} again,
we conclude that 
$\tau=(f\circ q_{\infty})(+\infty)=\geo f(q_{\infty}(+\infty))$, 
i.e.\ $\tau$ is in the image of $\geo f$.
\qed

\subsection{Morse actions}

We consider isometric actions $\Ga\acts X$ 
of finitely generated groups. 
\begin{dfn}[Morse action]
\label{dfn:morseact}
We call an action $\Ga\acts X$ {\em $\Theta$-Morse} 
if one (any) orbit map $\Ga\to\Ga x\subset X$ 
is a $\Theta$-Morse embedding 
with respect to a(ny) word metric on $\Ga$. 
We call an action $\Ga\acts X$ {\em $\tau_{mod}$-Morse} if it is $\Theta$-Morse 
for some $\tau_{mod}$-Weyl convex compact subset $\Theta\subset\ost(\tau_{mod})$. 
\end{dfn}

\begin{rem}[Morse actions are weakly regular and undistorted]
\label{rem:morseacundist}
(i) It follows immediately from the definition of Morse quasigeodesics that $\Theta$-Morse actions 
are $\tau_{mod}$-regular for the simplex type $\tau_{mod}$ determined by $\Theta$. 

(ii) Morse actions are {\em undistorted} in the sense 
that the orbit maps are quasi-isometric embeddings. 
In particular, 
they are properly discontinuous. 
\end{rem}
We denote by $\Hom_M(\Ga,G)\subset\Hom(\Ga,G)$ 
the subset of Morse actions $\Ga\acts X$.

By analogy with {\em local Morse quasigeodesics}, we define {\em local Morse group actions} $\rho: \Ga\acts X$ of a hyperbolic group 
(with fixed generating set):  

\begin{dfn}
An action $\rho$ is called {\em $(L,A,\Theta, D, S)$-locally Morse}, or 
{\em $(l, a, L,A,\Theta, D)$-locally Morse on the scale $S$},  if the 
orbit map $\Ga\to \Ga\cdot x\subset X$ is an $(l,a,L,A,\Theta, D,S)$-local Morse embedding. 
\end{dfn}

According to our local to global result for Morse embeddings, 
see Theorem~\ref{thm:locglobmqiembhyp}, 
an action of a word hyperbolic group is Morse 
if and only if it is local Morse on a sufficiently large scale. 
Since this is a finite condition, 
it follows that the Morse property is stable under perturbation of the action:
\begin{thm}[Morse is open for word hyperbolic groups]
\label{thm:morsestab}
For any word hyperbolic group $\Ga$ 
the subset $\Hom_M(\Ga,G)$ is open in $\Hom(\Ga,G)$. 
\end{thm}
\proof
Let $\rho:\Ga\acts X$ be a Morse action. 
We fix a word metric on $\Ga$ and a base point $x\in X$. 
Then there exist data $(L, A, \Theta, D)$ 
such that the orbit map $\Ga\to\Ga x\subset X$ sends is an $(L,A)$-quasi-isometric embedding, which sends 
(discrete) geodesics to $(L,A,\Theta,D)$-Morse quasigeodesics. 

We relax the Morse parameters slightly, 
i.e.\ we consider $(L,A,\Theta,D)$-Morse quasigeodesics 
as $(L,A+1,\Theta,D+1)$-Morse quasigeodesics
satisfying strict inequalities. 
For every scale $S$, 
the orbit map is, in particular,  
an $(L,A+1,\Theta,D+1,S)$-{\em local} Morse embedding. 
Due to $\Ga$-equivariance, 
this is a finite condition in the sense 
that it is equivalent to a condition involving 
only finitely many orbit points. 
Since we relaxed the Morse parameters, 
the same condition is satisfied by all actions sufficiently close to $\rho$.

Theorem~\ref{thm:locglobmqiembhyp} provides a scale $S$ 
such that $(L,A+1,\Theta,D+1,S)$-local Morse embeddings 
are global Morse. 
(More precisely, they are uniform $\Theta'$-Morse embeddings.)
It follows that actions sufficiently close to $\rho$ are ($\Theta'$-)Morse. 
\qed

\begin{cor}
For every hyperbolic group $\Ga$ the space of {\em faithful} Morse representations 
$$\Hom_{f,M}(\Ga,G)$$ is open in $\Hom_M(\Ga,G)$. 
\end{cor}
\proof Every hyperbolic group $\Ga$ has the unique maximal finite normal subgroup 
$F\triangleleft \Ga$ (if $\Ga$ is nonelementary then $F$ is the kernel of the action of $\Ga$ on $\geo \Ga$). Since Morse actions are properly discontinuous, kernel of every Morse representation $\Ga\to G$ is contained in $F$. Since $\Hom(F, G)/G$ is finite, it follows that the set of faithful Morse representations is open in $\Hom_M(\Ga,G)$. \qed

\medskip
We now turn to asymptotic properties of Morse actions. 
We apply our earlier discussion of the asymptotics of Morse embeddings. 

For a Morse action $\rho:\Ga\acts X$,  
the limit sets (\ref{eq:limsmemb}) of the orbit maps 
coincide with the limit set $\La_{\tau_{mod}}(\rho(\Ga))$ 
of the subgroup $\rho(\Ga)\subset G$. 
If $\Ga$ is word hyperbolic group, 
then the boundary map (\ref{eq:bdmapmemb})
induced by the orbit maps 
yield a well-defined $\Ga$-equivariant homeomorphism at infinity 
\begin{equation*}
\geo\rho:\geo\Ga\buildrel\cong\over\to
\La_{\tau_{mod}}(\Ga)\subset\Flag(\tau_{mod})
\end{equation*}
which does not depend on the $\Ga$-orbit, 
cf.\ Proposition~\ref{prop:asymorseemb}. 
The same proposition implies  
together with the fact that 
Morse embeddings are uniformly regular:

\begin{thm}[Asymptotic characterization of Morse actions]
\label{thm:asymorse}
An action of a word hyperbolic group 
is $\tau_{mod}$-Morse 
if and only if it is $\tau_{mod}$-asymptotically embedded.
\end{thm}
\proof
One direction follows from the discussion above, 
the converse from Theorem~\ref{thm:qiembrtaureg}.
\qed

\medskip
Our result on the openness of the Morse condition 
for actions of word hyperbolic groups, 
cf.\ Theorem~\ref{thm:morsestab}, 
can be strengthened in the sense that 
the asymptotics of Morse actions vary continuously:

\begin{thm}[Morse actions are structurally stable]
\label{thm:strcstab}
The boundary map at infinity of a Morse action 
depends continuously on the action. 
\end{thm}
\proof
Nearby actions are uniformly Morse, 
see the proof of Theorem~\ref{thm:morsestab}. 
The assertion therefore follows from the fact 
that the ends of Morse quasirays
vary continuously, 
cf.\ Lemma~\ref{lem:cont-dependence}.
\qed

\begin{rem}
(i)
Note that since the boundary maps at infinity are embeddings, 
the $\Ga$-actions on the $\tau_{mod}$-limit sets 
are topologically conjugate to each other and, for nearby actions, 
by a homeomorphism close to the identity.

(ii)
In rank one, our argument yields a different proof for 
Sullivan's Structural Stability Theorem \cite{Sullivan}
for convex cocompact group actions. 
\end{rem}

\subsection{Schottky actions}
\label{sec:schottky actions}

In this section we apply our local-to-global result for straight sequences 
(Theorem~\ref{thm:locstrimplcoastrseq})
to construct Morse actions of free groups, 
generalizing and sharpening\footnote{In the sense that we obtain free subgroups which are not only embedded, but also asymptotically embedded in $G$.} 
Tits's ping-pong construction.

\medskip 
We consider two oriented $\tau_{mod}$-regular geodesic lines 
$a,b$ in $X$. 
Let $\tau_{\pm a},\tau_{\pm b}\in\Flag(\tau_{mod})$ 
denote the simplices 
which they are $\tau$-asymptotic to, 
and let $\theta_{\pm a},\theta_{\pm b}\in\si_{mod}$ 
denote the types of their forward/backward ideal endpoints in $\geo X$.
(Note that $\theta_{-a}=\iota(\theta_a)$
and $\theta_{-b}=\iota(\theta_b)$.) Let $\Theta$ be a compact convex subset of $\ost(\tau_{mod})\subset \si_{mod}$, which is invariant under $\iota$.

\begin{dfn}[Generic pair of geodesics]
We call the pair of geodesics $(a,b)$ {\em generic} 
if the four simplices $\tau_{\pm a},\tau_{\pm b}$ 
are pairwise opposite. 
\end{dfn}
Let $\al,\beta\in G$ be axial isometries 
with axes $a$ and $b$ respectively 
and translating in the positive direction along these geodesics. 
Then $\tau_{\pm a}$ and $\tau_{\pm b}$ 
are the attractive/repulsive fixed points of $\al$ and $\beta$ 
on $\Flag(\tau_{mod})$. 

For every pair of numbers $m,n\in\N$ 
we consider the representation of the free group in two generators
\begin{equation*}
\rho_{m,n}: F_2=\<A,B\> \to G
\end{equation*}
sending the generator $A$ to $\al^m$ and $B$ to $\beta^n$. 
We regard it as an isometric action 
$\rho_{m,n}:F_2\acts X$.

\begin{definition}[Schottky subgroup]
A {\em $\tau_{mod}$-Schottky subgroup} of $G$ is a free $\tau_{mod}$-asymp\-to\-ti\-cally embedded subgroup of $G$. 
\end{definition}

If $G$ has rank one, this definition amounts to the requirement that $\Ga$ is convex cocompact and free. Equivalently, this is a discrete finitely generated subgroup of $G$ which contains  no nontrivial elliptic and parabolic elements and has totally disconnected limit set (see see \cite{Kapovich2007}). We note that this definition essentially agrees with the standard definition of Schottky groups in rank 1 Lie groups, provided one allows fundamental domains at infinity for such groups to be bounded by pairwise disjoint compact submanifolds which need not be topological spheres, see \cite{Kapovich2007} for the detailed discussion.

\begin{thm}[Morse Schottky actions]
\label{thm:mschott}
If the pair of geodesics $(a,b)$ is generic 
and if $\theta_{\pm a},\theta_{\pm b}\in\interior(\Theta)$, 
then the action $\rho_{m,n}$ is $\Theta$-Morse 
for sufficiently large $m,n$. Thus, such $\rho_{m,n}$ is injective and its image is a $\tau_{mod}$-Schottky subgroup 
of $G$. 
\end{thm}
\begin{rem}
In particular,
these actions are faithful and undistorted, 
compare Remark~\ref{rem:morseacundist}. 
\end{rem}
\proof
Let $S=\{A^{\pm1},B^{\pm1}\}$ be the standard generating set. 
We consider the sequences $(\ga_k)$ in $F_2$
with the property that $\ga_k^{-1}\ga_{k+1}\in S$ and 
$\ga_{k+1}\neq\ga_{k-1}$ for all $k$.  
They correspond to the geodesic segments in the Cayley tree of $F_2$ 
associated to $S$ which connect vertices. 

Let $x\in X$ be a base point. 
In view of Lemma~\ref{lem:morsefromhyp} 
we must show that the corresponding sequences $(\ga_kx)$ 
in the orbit $F_2\cdot x$
are uniformly $\Theta$-Morse. 
(Meaning e.g.\ that the maps $\R\to X$ 
sending the intervals $[k,k+1)$ to the points $\ga_kx$ are 
uniform $\Theta$-Morse quasigeodesics.)
As in the proof of Theorem~\ref{thm:locglobmqg}
we will obtain this by applying 
our local to global result for straight spaced sequences 
(Theorem~\ref{thm:locstrimplcoastrseq}) 
to the associated midpoint sequences. 
Note that the sequences $(\ga_kx)$ themselves cannot expected to be straight. 

Taking into account the $\Ga$-action, 
the uniform straightness of all midpoint sequences 
depends on the geometry of a finite configuration in the orbit. 
It is a consequence of the following fact.
Consider the midpoints $y_{\pm m}$ of the segments $x \al^{\pm m}(x)$ 
and $z_{\pm n}$ of the segments 
$x \beta^{\pm n}(x)$. 
\begin{lem}
\label{lem:geomquad}
For sufficiently large $m,n$ 
the quadruple $\{y_{\pm m},z_{\pm n}\}$ is arbitrarily separated 
and $\Theta$-regular.
Moreover, 
for any of the four points, 
the segments connecting it to the other three points 
have arbitrarily small $\zeta$-angles with the segment connecting it to $x$. 
\end{lem}
\proof
The four points are arbitrarily separated from each other and from $x$ 
because the axes $a$ and $b$ diverge from each other 
due to our genericity assumption. 

By symmetry, it suffices to verify the rest of the assertion 
for the point $y_m$, 
i.e.\ we show that 
the segments $y_my_{-m}$ and $y_mz_n$ are $\Theta$-regular 
for large $m,n$
and that 
$\lim_{m\to\infty} \angle^\zeta_{y_m}(x, y_{-m})=0$ 
and 
$\lim_{n,m\to\infty} \angle^\zeta_{y_m}(x, z_n)=0$. 

The orbit points $\al^{\pm m}x$ and the midpoints $y_{\pm m}$
are contained in a tubular neighborhood of the axis $a$. 
Therefore, 
the segments $y_mx$ and $y_my_{-m}$ 
are $\Theta$-regular for large $m$ 
and $\angle_{y_m}(x, y_{-m})\to0$. 
This implies that also $\angle^{\zeta}_{y_m}(x, y_{-m})\to0$. 

To verify the assertion for $(y_m,z_n)$ 
we use that, due to genericity, 
the simplices $\tau_a$ and $\tau_b$ are opposite 
and we consider the parallel set $P=P(\tau_a,\tau_b)$.
Since the geodesics $a$ and $b$ are forward asymptotic to $P$, 
it follows that the points $x,y_m,z_n$ have 
uniformly bounded distance from $P$. 
We denote their projections to $P$ by $\bar x,\bar y_m,\bar z_n$. 

Let $\Theta''\subset\interior(\Theta)$ be an auxiliary Weyl convex subset
such that 
$\theta_{\pm a},\theta_{\pm b}\in\interior(\Theta'')$. 
We have that 
$\bar y_m\in V(\bar x,\st_{\Theta''}(\tau_a))$ for large $m$ 
because the points $y_m$ lie in a tubular neighborhood 
of the ray with initial point $\bar x$ and asymptotic to $a$. 
Similarly, 
$\bar z_n\in V(\bar x,\st_{\Theta''}(\tau_b))$ for large $n$. 
It follows that 
$\bar x\in V(\bar y_m,\st_{\Theta''}(\tau_b))$ and,
using the convexity of $\Theta$-cones (Proposition~\ref{prop:thconeconv}), 
that $\bar z_n\in V(\bar y_m,\st_{\Theta''}(\tau_b))$. 

The cone $V(y_m,\st_{\Theta''}(\tau_b))$
is uniformly Hausdorff close to the cone 
$V(\bar y_m,\st_{\Theta''}(\tau_b))$
because the Hausdorff distance of the cones 
is bounded by the distance $d(y_m,\bar y_m)$ of their tips. 
Hence there exist points $x',z'_n\in V(y_m,\st_{\Theta''}(\tau_b))$ 
uniformly close to $x,z_n$. 
Since $d(y_m,x'),d(y_m,z'_n)\to\infty$ as $m,n\to\infty$, 
it follows that 
the segments $y_mx$ and $y_mz_n$ are $\Theta$-regular for large $m,n$. 
Furthermore, 
since $\angle^\zeta_{y_m}(x',z'_n)=0$ 
and 
$\angle_{y_m}(x,x')\to0$
as well as
$\angle_{y_m}(z_n,z'_n)\to0$, 
it follows that $\angle^\zeta_{y_m}(x,z_n)\to0$. 
\qed

\medskip
{\em Proof of Theorem concluded.} 
The lemma implies that for any given $l,\eps$ 
the midpoint triples of the four point sequences $(\ga_kx)$
are $(\Theta,\eps)$-straight and $l$-spaced 
if $m,n$ are sufficiently large, 
compare the quadruple condition (Definition~\ref{dfn:quad}). 
This means that the midpoint sequences 
of all sequences $(\ga_kx)$ are $(\Theta,\eps)$-straight and $l$-spaced 
for large $m,n$. 
Theorem~\ref{thm:locstrimplcoastrseq} then implies 
that the sequences $(\ga_kx)$ are uniformly $\Theta$-Morse. 
\qed 

\begin{rem}
Generalizing the above argument to free groups with finitely many generators, 
one can construct Morse Schottky subgroups 
for which the set $\theta(\La)\subset\si_{mod}$ of types of limit points 
is arbitrarily Hausdorff close to a given $\iota$-invariant 
Weyl convex subset $\Theta$. 
This provides an alternative approach 
to the second main theorem in \cite{Benoist}
using geometric arguments. 
\end{rem}

\subsection{Algorithmic recognition of Morse actions}
\label{sec:algrec}

In this section, we describe an algorithm which has an isometric action $\rho: \Ga\acts X$ and a point $x\in X$ 
as its input and terminates if and only if the action $\rho$ is Morse (otherwise, the algorithm runs forever). 

We begin by describing briefly a {\em J{\o}rgensen's algorithm} accomplishing a similar task, namely, detecting 
geometrically finite actions on $X=\H^3$. 
Suppose that we are given a finite (symmetric) set of generators $g_1=1,\ldots, g_m$ 
of a subgroup  $\Ga\subset PO(3,1)$ and a base-point $x\in X=\H^n$. The idea of J{\o}rgensen's algorithm is to construct a finite sided Dirichlet fundamental domain $D$ for $\Ga$ (with the center at $x$): Every geometrically finite subgroup of 
$PO(3,1)$ admits such a domain. (The latter is false for geometrically finite subgroups of $PO(n,1)$, $n\ge 4$, but is, nevertheless true for convex cocompact subgroups.) Given  a finite sided convex fundamental domain, one concludes that $\Ga$ is geometrically finite. 
Here is how the algorithm works: For each $k$ define the subset $S_k\subset \Ga$ represented by words of length $\le k$ in  the letters $g_1,\ldots, g_m$. For each $g\in S_k$ consider the half-space $Bis(x, g(x))\subset X$ bounded by the bisector of the segment $x g(x)$ and containing the point $x$. Then compute the intersection
$$
D_k=\bigcap_{g\in S_k} Bis(x, g(x)).
$$ 
Check if $D_k$ satisfies the conditions of the {\em Poincar\'e's Fundamental Domain theorem}. 
If it does, then $D=D_k$ is a finite sided fundamental domain of $\Ga$. If not, increase $k$ by $1$ and repeat the process. Clearly, this process terminates if and only if $\Ga$ is geometrically finite. 

One can enhance the algorithm in order to detect if a geometrically finite group is convex cocompact. Namely, after a Dirichlet domain $D$ is constructed, one checks for the following:

1. If the ideal boundary of a Dirichlet domain $D$  has isolated ideal points (they would correspond to 
rank two cusps which are not allowed in convex cocompact groups). 

2. If the ideal boundary of $D$  contains tangent circular arcs with points of tangency fixed by parabolic elements (coming from the ``ideal vertex cycles''). Such points correspond to rank 1 cusps, which again are not allowed  
 in convex cocompact groups. 
 
Checking 1 and 2 is a finite process; after its completion, one concludes that $\Ga$ is convex cocompact.  
 
\medskip
We now consider group actions on general symmetric spaces. Let $\Ga$ be a hyperbolic group with a fixed finite (symmetric) generating set; we equip the group $\Ga$ with the word metric determined by this generating set. 

For each $n$, let ${\mathcal L}_n$ denote the set of maps $q: [0, 3n]\cap \Z \to \Gamma$ which are restrictions  of geodesics $\tilde{q}: \Z\to \Ga$, so that $q(0)=1\in \Ga$. In view of the geodesic automatic structure on $\Ga$ (see e.g. \cite[Theorem 3.4.5]{Epstein}), the set  ${\mathcal L}_n$ can be described via a finite state automaton. 

Suppose that $\rho: \Gamma\acts X$ is an isometric action on a symmetric space $X$; we fix a base-point $x\in X$ 
and the corresponding orbit map $f: \Ga\to \Ga x\subset X$. We also fix an $\iota$-invariant face $\tau_{mod}$ of the model spherical simplex $\si_{mod}$ of $X$. The algorithm that we are about to describe will detect that the action $\rho$ is $\tau_{mod}$-Morse. 

\begin{rem}
If the face $\tau_{mod}$ is not fixed in advance, we would run algorithms for each face $\tau_{mod}$ in parallel. 
\end{rem}

For the algorithm we will be using a special (countable) increasing family of Weyl-convex compact subsets $\Theta=\Theta_i\subset \ost(\tau_{mod})\subset \si_{mod}$ which exhausts $\ost(\tau_{mod}$; in particular, 
every compact $\iota$-invariant convex subset of $\ost(\tau_{mod})\subset \si_{mod}$ is contained in some $\Theta_i$: 
\begin{equation}\label{Theta_i}
\Theta_i:=\{v\in \si: \min_{\alpha\in \Phi_{\tau_{mod}}} \alpha(v)\ge \frac{1}{i} \},
\end{equation}
where $\Phi_{\tau_{mod}}$ is the subset of the set of simple roots $\Phi$ (with respect to $\si_{mod}$) which vanish on the face $\tau_{mod}$. Clearly, the sets $\Theta_i$ satisfy the required properties. Furthermore, we consider only those $L$ and $D$ which are natural numbers. 

Next, consider the sequence
$$
(L_i, \Theta_{i}, D_i)= (i, \Theta_i, D_i), i\in \N. 
$$

In order to detect $\tau_{mod}$-Morse actions we will use the local characterization of Morse quasigeodesics given by 
Theorem  \ref{thm:locstrimplcoastrseq} and Proposition \ref{prop:morseimplquad}. Due to the discrete nature of quasigeodesics that we will be considering, it suffices to assume that the additive quasi-isometry constant $A$ is zero. 

Consider the functions 
$$
l(\Theta, \Theta', \delta), \eps(\Theta, \Theta', \delta)
$$
as in Theorem \ref{thm:locstrimplcoastrseq}. Using these functions,  for the sets $\Theta=\Theta_{i}, \Theta'=\Theta_{i+1}$ and the constant $\delta=1$ we define the numbers 
$$
l_i=l(\Theta, \Theta', \delta), \eps_i= \eps(\Theta, \Theta', \delta). 
$$

Next, for the numbers $L=L_i, D=D_i$ and the sets  $\Theta=\Theta_{i}, \Theta'=\Theta_{i+1}$, consider the numbers 
$$
s_i=s(L_i,0, \Theta_i, \Theta_{i+1}, D_i, \eps_{i+1}, l_{i+1})
$$ 
as in Proposition \ref{prop:morseimplquad}. According to this proposition, every $(L_i,0,\Theta_i,D_i)$-Morse quasigeodesic satisfies the $(\Theta_{i+1},\eps_{i+1},l_{i+1}, s)$-quadruple condition for all $s\ge s_i$. We note that, a priori, the sequence $s_i$ need not  be increasing. We set $S_1=s_1$ and 
define a monotonic sequence $S_i$ recursively by 
$$
S_{i+1}= \max(S_i, s_{i+1}). 
$$
Then every $(L_i, 0, \Theta_i, D_i)$-Morse quasigeodesic also satisfies the 
$(\Theta_{i+1}, \eps_{i+1}, l_{i+1}, S_{i+1})$-quadruple condition. 

We are now ready to describe the algorithm. For each $i\in \N$ we compute the numbers $l_i, \eps_i$ and, then, $S_i$, as above.  
We then consider finite discrete paths in $\Ga$, $q\in {\mathcal L}_{S_i}$, and the corresponding discrete paths 
in $X$, $p(t)=q(t)x$, $t\in [0, 3S_i]\cap \Z$. The number of paths $q$ (and, hence, $p$) for each $i$ is finite, bounded by the growth function of the group 
$\Ga$. 

For each discrete path $p$ we check the $(\Theta_{i}, \eps_{i}, l_{i}, S_{i})$-quadruple condition. If for some $i=i_*$, all paths $p$ satisfy this condition, 
the algorithm terminates: It follows from Theorem \ref{thm:locstrimplcoastrseq} that the map $f$ sends all normalized discrete biinfinite geodesics in $\Ga$ to Morse quasigeodesics in $X$. Hence, the action $\Ga\acts X$ is Morse in this case. Conversely, suppose that the action of $\Ga$ is $(L, 0, \Theta, D)$-Morse. Then $f$ sends all isomeric embeddings $\tilde{q}: \Z\to \Ga$ to $(L, 0, \Theta, D)$-Morse quasigeodesics $\tilde{p}$ in $X$. In view of the properties of the sequence 
$$
(L_i, \Theta_i, D_i),
$$
it follows that for some $i$, 
$$
(L, \Theta, D)\le (L_i, \Theta_i, D_i), 
$$
i.e., $L\le L_i, \Theta\subset \Theta_i, D\le D_i$; hence, all the biinfinite discrete paths $\tilde p$ are $(L_i, 0, \Theta_i, D_i)$-Morse quasigeodesic. 
By the definition of the numbers $l_i, \eps_i, S_i$, it then follows that all the discrete paths $p=f\circ q, q\in {\mathcal L}_{S_i}$ satisfy the  
$(\Theta_{i+1}, \eps_{i+1}, l_{i+1}, S_{i+1})$-quadruple condition. Thus, the algorithm will terminate at the step $i+1$ in this case. 

Therefore, the algorithm terminates if and only if the action is Morse (for some parameters). If the action is not Morse, the algorithm will run forever. \qed

Addresses:

\noindent M.K.: Department of Mathematics, \\
University of California, Davis\\
CA 95616, USA\\
email: kapovich@math.ucdavis.edu

\noindent B.L.: Mathematisches Institut\\
Universit\"at M\"unchen \\
Theresienstr. 39\\ 
D-80333, M\"unchen, Germany\\ 
email: b.l@lmu.de

\noindent J.P.: Departament de Matem\`atiques,\\
 Universitat Aut\`onoma de Barcelona,\\ 
 08193 Bellaterra, Spain\\
email: porti@mat.uab.cat

\end{document}